\documentclass[reqno,final]{amsart}
\usepackage{natbib}  
\usepackage{fancyhdr} 
\usepackage{color} 
\usepackage{graphicx} 

\usepackage{amssymb}

\usepackage{yfonts}
\usepackage{euscript}

\usepackage[applemac]{inputenc}
\usepackage{amsmath,amsthm,dsfont}
\usepackage{graphics}
\usepackage[T1]{fontenc}
\usepackage{epsfig,psfrag}

\RequirePackage[colorlinks,citecolor=blue,urlcolor=blue]{hyperref}
\newcommand{\gr}[1]{{\color{black}{#1}}}

\newcommand{\alex}[1]{{\color{black}{#1}}}
\newcommand{\bl}[1]{{\color{black}{#1}}}

\definecolor{aleacolor}{rgb}{0.16,0.59,0.78}

\hypersetup{
breaklinks,
colorlinks=true,
linkcolor=aleacolor,
urlcolor=aleacolor,
citecolor=aleacolor}

\renewcommand{\headheight}{11pt}
\renewcommand{\cite}{\citet}

\theoremstyle{plain}
\newtheorem{theorem}{Theorem}[section]
\newtheorem{proposition}[theorem]{Proposition}
\newtheorem{lemma}[theorem]{Lemma}

\theoremstyle{definition}

\theoremstyle{remark}
\newtheorem{remark}[theorem]{Remark}

\makeatletter \@addtoreset{equation}{section} \makeatother
\renewcommand\theequation{\thesection.\arabic{equation}}
\renewcommand\thefigure{\thesection.\arabic{figure}}
\renewcommand\thetable{\thesection.\arabic{table}}

\newcommand{\aleaIndex}[1]{\href{http://alea.impa.br/english/index_v#1.htm}{\bf #1}}
\eheader{}{\aleaIndex{}}{2016}{1}{61}

\newcommand{\N}{\mathbb{Z}_{+}}
\newcommand{\Z}{\mathbb{Z}}

\newcommand{\R}{\mathbb{R}}

\newtheorem{fact}[theorem]{Fact}

\newcommand{\egloi}{\stackrel{\textgoth{L}}{=}}
\newcommand{\cvloiS}{\stackrel{\textgoth{L}_{S}}{\rightarrow}}
\newcommand{\cvloi}{\stackrel{\textgoth{L}_{}}{\rightarrow}}

\def\un{\mathds{1}}
\def\a{\alpha}

\def\d{\delta}
\def\e{\varepsilon}

\def\E{\mathbb{E}}
\def\CK{\mathcal{C}_{\kappa}}

\def\Dt{\mathcal{D}}
\def\B{\mathcal{A}}

\def\o{\omega}

\def\H{\mathcal{H}}

\def\N{\mathbb{N}}
\def\P{\mathbb{P}}
\def\Pw{\P^{W_{\kappa}}}

\def\Ymi{\mathcal{Y}_2^{-1}}
\def\Ymj{\mathcal{Y}_1^{\natural}}

\def\BP{{W_{\kappa}^{\uparrow}}}

\def\Cct{C_0}

\def\R{\mathbb{R}}

\def\BQ{Q}

\def\tt{\theta^-}

\def\Z{\mathbb{Z}}

\def\tV{\tilde{V}^{(j)}}

\def\loX{\mathcal{L}}

\def\mV{\mathcal{V}}
\def\tA2{\tilde{A}({\tilde L_2})}
\def \tts2{\tilde \tau^+_2(  h_t/2 ))}
\def \tts1{\tilde \tau_1(  h_t/2 ))}

\def\L{{\ell}}

\def\Ip{{\mathcal{I}}^+}
\def\Im{{\mathcal{I}}^-}

\def\tt{\tilde{\tau}}

\def\Lt{\tilde{L}}
\def\mt{\tilde{m}}

\def\mf{{m}}
\def\Mf{{M}}

\def\do{E_{\omega}}

\def\tm{\tilde m_j}
\def\tL{\tilde L_j}

\def\lo{\mathcal L}

\def\siX{\sigma}

\def\bU{ {\bf U} }
\def\k{\kappa}
\def\dd{\textnormal{d}}
\def\wk{W_{\k}}

\title[Renewal structure and local time for diffusions in random environment]
{Renewal structure and local time for diffusions in random environment}

\def\lx{\loX}

\author{Pierre Andreoletti}
\address{Universit\'e d'Orl\'eans, Laboratoire MAPMO - Fédération Denis Poisson, B\^atiment de math\'ematiques - Rue de Chartres
B.P. 6759 - 45067 Orl\'eans cedex 2,
France}
\email{Pierre.Andreoletti@univ-orleans.fr}

\author{Alexis Devulder}
\address{
Laboratoire de Math\'ematiques de Versailles, UVSQ, CNRS, Universit\'e Paris-Saclay, 78035 Versailles, France.}
\email{devulder@math.uvsq.fr }

\author{Grégoire V\'echambre}
\address{Universit\'e d'Orl\'eans, Laboratoire MAPMO - Fédération Denis Poisson, B\^atiment de math\'ematiques - Rue de Chartres
B.P. 6759 - 45067 Orl\'eans cedex 2,
France}

\email{Gregoire.Vechambre@ens-rennes.fr}


\date{\today}
\thanks{This research was partially supported by the french ANR project MEMEMO2 2010 BLAN 0125.}

\begin{document}

\maketitle

\begin{abstract}
We study a one-dimensional diffusion $X$ in a drifted Brownian potential  $W_\kappa$, with $ 0<\kappa<1$,  and focus on the behavior of the local times $(\loX(t,x),x)$ of $X$ before time $t>0$.
In particular we characterize the limit law of the supremum of the local time, as well as the position of the favorite site. These limits can be written explicitly from a two dimensional stable Lévy process. Our analysis is based on the study of an extension of the renewal structure which is deeply involved in the asymptotic behavior of $X$.

\end{abstract}


\section{Introduction}
\subsection{Presentation of the model}
Let $(X(t),\ t\geq 0)$ be a diffusion in a random c\`adl\`ag potential $(V(x),\ x\in\R)$, defined informally by $X(0)=0$ and
$$
    \text{d}X(t)=\text{d}\beta(t)-\frac{1}{2}V'(X(t))\text{d}t,$$
where $(\beta(s),\ s\geq 0)$ is a Brownian motion independent of $V$. Rigorously, $X$ is defined by its conditional generator given $V$,
\begin{align*} & \frac{1}{2}e^{V(x)}\frac{\text{d}}{\text{d} x}\left(e^{-V(x)}\frac{\text{d}}{\text{d} x}\right).
\end{align*}
We put ourselves in the case where $V$ is a negatively drifted Brownian motion:
$V(x)=\wk(x):=W(x)-\frac{\k}{2}x,\ x\in\R$, with $0<\kappa<1$ and $(W(x), \ x\in\R)$ is a two sided Brownian motion. We explain at the end of Section \ref{Results} what should be done to extend our results to a more general Lévy potential.

We denote by $P$ the probability measure associated to $\wk(.)$. The probability conditionally on the potential $W_{\kappa}$ is denoted by $\P^{W_{\kappa}}$ and is called the {\it quenched probability}.
We also define the {\it annealed probability} as
$$
    \P(.)
:=
    \int \P^{W_{\kappa}}(.) P(\wk\in \dd\o).
$$
We denote respectively by $\E^{W_{\kappa}}$, $\E$, and $E$ the expectations with regard to $\P^{W_{\kappa}}$, $\P$ and $P$.
In particular, $X$ is a Markov process under $\P^{W_{\kappa}}$ but not under $\P$.

This diffusion $X$ has been introduced by \cite{Schumacher}. It is generally considered as a continuous time analogue of
random walks in random environment (RWRE). We refer e.g. to \cite{Zeitouni_St_Flour} for general properties of RWRE.

In our case, since $\k>0$, the diffusion $X$ is a.s. transient and its asymptotic behavior was first studied by Kawazu and Tanaka:
if $H(r)$ is the hitting time of $r\in\R$ by $X$,
\begin{align}
H(r):=\inf\{s>0,\ X(s)=r\} \label{DefH},
\end{align}
\cite{KawazuTanaka} proved that, for $0<\k<1$
under the annealed probability $\P$,  $H(r)/r^{1/\k}$ converges in law as $r\to+\infty$
to a $\kappa$-stable distribution (see also \cite{HuShYo}, and
\cite{Tanaka2}). Here we are interested in the local time of $X$,
which is the jointly continuous process $(\loX(t,x),\, t>0,\, x \in \R)$ satisfying, for any positive measurable function $f$,
$$
    \int_0^t f(X(s)) \dd s
=
    \int_{-\infty}^{+\infty} f(x) \loX(t,x)\dd x,
\qquad
    t>0.
$$
One quantity of particular interest is the supremum of the local time of $X$ at time $t$, defined as
$$
    \lo^*(t)
:=
    \sup_{x \in \R} \loX(t,x),
\qquad
    t>0.
$$

For Brox's diffusion, that is, for the diffusion $X$ in the recurrent case $\kappa=0$, it is proved in  \cite{AndDiel}
that the local time process until time $t$ re-centered at the localization coordinate $b_t$ (see \cite{Brox})
and renormalized by $t$ converges in law under the annealed probability $\P$.
This allows the authors of \cite{AndDiel} to derive the limit law of the supremum of the local time at time $t$
as $t\to+\infty$.
We recall their result below in order to compare it with the results of the present paper.
To this aim, we introduce for every $\kappa\geq 0$,
\begin{align}
\mathcal{R}_{\kappa} :=\int_0^{+ \infty} e^{- W_{\kappa}^{\uparrow}(x)} \dd x
+ \int_0^{+ \infty} e^{- \widetilde{W}_{\kappa}^{\uparrow}(x)} \dd x,
\label{Rk}
\end{align}
where $\big(W_{\kappa}^{\uparrow}(x),\, x\geq 0\big)$ and $\big(\widetilde{W}_{\kappa}^{\uparrow},\, x\geq 0\big)$
are two independent copies of the process $({W}_{\kappa}(x),\, x\geq 0)$ Doob-conditioned to remain positive.

\begin{theorem}\label{cvloi} (\cite{AndDiel}) If $\kappa=0$, then
\begin{align*}
    \frac{\loX^{*}(t)}{t}
\cvloi
    \frac{1}{\mathcal{R}_{\kappa}},
\end{align*}
where $\cvloi$ denotes convergence in law under the annealed probability $\P$ as $t\to +\infty$.
\end{theorem}

Extending their approach, and following the results of \cite{Shi}, \cite{Diel}
obtains the non-trivial normalizations for the almost sure behavior of the $\limsup$ and the $\liminf$ of $\lo^*(t)$ as $t\to+\infty$
when $\kappa=0$.
Notice that corresponding results have been previously established
in \cite{DemGanPerShi} and \cite{GanPerShi}
for the discrete analogue of $X$ in the recurrent case $\kappa=0$, the recurrent RWRE generally called Sinai's random walk.

One of our aims in this paper is to extend the study of the local time of $X$ in the case $0< \kappa <1$, and deduce from that the weak asymptotic behavior of ${\loX^{*}(t)}$ suitably renormalized as $t\to+\infty$.

Before going any further, let us recall to the reader what is known for the slow transient cases.
For transient RWRE in the case $0<\kappa\leq 1$ (see \cite{KesKozSpi} for the seminal paper), a result of \cite{GanShi} states the almost sure behavior for the $\limsup$ of the supremum of the local time $\lo_S^*(n)$ of these random walks (denoted by $S$) at time $n$: there exists a constant $c>0$ such that
$\limsup_{n \rightarrow + \infty} \lo_S^*(n)/n=c>0$ $\P$ almost surely.
Contrarily to the recurrent case (\cite{GanPerShi}) their method,
based on a relationship between the RWRE $S$ and a branching process in random environment,
cannot be exploited to determine the limit law of $\lo_S^*(n)/n$.

For the transient diffusion $X$ considered here, the only paper dealing with $\loX^*(t)$ is \cite{Devulder}, in which it is proved,
among other results, that when $0<\kappa<1$,
$\limsup_{t \rightarrow + \infty} \loX^*(t)/t=+ \infty$ almost surely. But once again his method cannot be used to characterize the limit law of $\loX^*(t)/t$ in the case $0<\kappa<1$.

Our motivation here is twofold, first we prove that our approach enables to characterize the limit law of $\loX^*(t)/t$ and open a way to determine the correct almost sure behavior of $\loX^*(t)$ as was done for Brox's diffusion by  \cite{Shi} and \cite{Diel}. Second we  make a first step on a specific way to study the local time which  could be used in estimation problems in random environment, see
 \cite{AdeEnr}, \cite{Pierre4}, \cite{AndDiel2}, \cite{AndLouMat},
 \cite{Comets_etal2}, \cite{Comets_etal}, \cite{FaLoMa}.

The method we develop here is an improvement of the one used in \cite{AndDev} about the localization of $X(t)$ for large $t$.

Before recalling the main result of this paper \cite{AndDev}, we need to introduce some new objects. We start with the notion of $h$-extrema, with $h>0$, introduced by \cite{NevPit} and studied more specifically in our case of drifted Brownian motions by \cite{Faggionato}.
For $h>0$, we say that $x\in\R$ is an {\it $h$-minimum} for a given continuous function $f$, $\R\to\R$,
if there exist $u<x<v$ such that $f(y)\geq f(x)$ for all $y\in[u,v]$, $f(u)\geq f(x)+h$ and $f(v)\geq f(x)+h$.
Moreover, $x$ is an {\it $h$-maximum} for $f$ iff $x$ is an $h$-minimum for $-f$.
Finally, $x$ is an {\it $h$-extremum} for $f$ iff it is an $h$-maximum or an $h$-minimum for $f$.

As we are interested in the diffusion $X$ until time $t$ for large $t$, we only focus on the $h_t$-extrema of $\wk$, where
$$
    h_t
:=
    \log t-\phi(t),\quad  \textrm{ with }  0 < \phi(t) =o(\log t),\quad  \log \log t=o(\phi(t)),
$$
and $t\mapsto \phi(t)$ is an increasing function,
as in \cite{AndDev}.
It is known (see \cite{Faggionato}) that almost surely, the $h_t$-extrema of $\wk$ form a sequence indexed by
$\Z$, unbounded from below and above, and that the $h_t$-minima and $h_t$-maxima alternate.
We denote respectively by $(\mf_j,\ j \in \Z)$ and $(\Mf_j,\ j \in \Z)$ the increasing sequences of $h_t$-minima
and of $h_t$-maxima of $\wk$,
such that $\mf_{0}\leq 0<\mf_1$ and $\mf_j<\Mf_j<\mf_{j+1}$ for every  $j\in\Z$. Define
\begin{equation}\label{eqDefNt}
    N_t
:=
    \max\Big\{k\in\N,\ \sup_{ 0 \leq s \leq t}X(s) \geq \mf_{k}\Big\},
\end{equation}
the number of (positive) $h_t$-minima on $\R_+$ visited by $X$ until time $t$. We have the following result.

\medskip

\begin{theorem} \label{ththm} (\cite{AndDev})
Assume $ 0<\kappa<1$. There exists a constant $ \mathcal{C}_1 >0$, such that
\[
    \lim_{t \rightarrow + \infty }
    \P\Big( \big|X(t)- \mf_{N_{t}}\big| \leq \mathcal{C}_1 \phi(t) \Big)
=
    1. \]
\end{theorem}

\medskip

This result proves that before time $t$, the diffusion $X$  visits the $N_t$ leftmost positive $h_t$-minima, and then gets stuck
in a very small neighborhood of an ultimate one, which is $\mf_{N_{t}}$.
An analogous result was proved for transient RWRE
 in the zero speed regime $0<\kappa<1$
by \cite{EnSaZi}. This phenomenon is due to two facts: the first one is the appearance of a renewal structure which is composed of the times it takes  the process to move from one $h_t$-minimum to the following one. The second is the fact that like in Brox's case $\kappa=0$,
the process is trapped a significant amount of time in the neighborhood of the local minimum $m_{N_t}$.

It is the extension of this renewal structure to the sequence of local times at the $h_t$-minima that we study here. We now detail our results.

\subsection{Results\label{Results}}

Let us introduce some notation involved in the statement of our  results.
Assume that $0<\kappa<1$.
\noindent Denote by $\big(D\big([0,+ \infty),\R^2\big),J_1\big)$ the space of c\`adl\`ag functions $[0,+ \infty)\to\R^2$
with $J_1$-Skorokhod topology and denote by  $\cvloiS$ the convergence in law for this topology.
On this space, define a 2-dimensional L\'evy process $(\mathcal{Y}_1, \mathcal{Y}_2)$ taking values in $\R_+\times \R_+$, which is a pure positive jump process with $\kappa$-stable Lévy measure  $\nu$ given by
\begin{equation}\label{eqDefMesureNU}
    \forall x > 0, \, \forall y > 0,
\qquad
    \nu \big( [x, + \infty[ \times [y,+ \infty[ \big) = \bl{\frac{\mathcal{C}_2}{y^{\kappa}}} \mathbb{E} \left [ (\mathcal{R}_{\kappa})^{\kappa} \mathds{1}_{\mathcal{R}_{\kappa} \leq \frac{y}{x}} \right ] + \bl{\frac{\mathcal{C}_2}{x^{\kappa}}} \mathbb{P} \left ( \mathcal{R}_{\kappa} > \frac{y}{x} \right ),
\end{equation}
where $\mathcal{R}_{\kappa}$ is defined in \eqref{Rk} \bl{and $\mathcal{C}_2$ is a positive constant (see Lemma \ref{lemproba}).
The Laplace transform of $\mathcal{R}_{\kappa}$ is given by
\begin{equation}\label{eqTransfoLaplaceRk}
    E\big(e^{-\gamma \mathcal{R}_{\kappa}}\big)
=
    \left(\frac{ (2\gamma)^{{\k/2}}}{\k\Gamma(\k)I_\k(2\sqrt{2\gamma})}\right)^2
\qquad
    \gamma>0,
\end{equation}
as proved in Lemma \ref{LTR} below, where $I_\k$ is the modified Bessel function of the first kind of index $\kappa$.
Moreover, $\mathcal{R}_{\kappa}$ admits moments of any positive order (see also Lemma \ref{LTR}).
In particular $\mathbb{E}[ (\mathcal{R}_{\kappa})^{\kappa} ]$ is finite and $\nu$ is well defined. } \\
 For a given c\`adl\`ag function $f$ in $D([0,+ \infty),\R)$, define for any $s>0, a>0$:
$$
    f^{\natural}(s):= \sup_{0\leq r \leq s}(f(r)-f(r^-)),
\qquad
    f^{-1}(a):= \inf\{x \geq 0,\ f(x) > a\},
$$
where $f(r^-)$ denotes the left limit of $f$ at $r$. In words, $f^{\natural}(s)$ is the largest jump of $f$ before time $s$,
whereas $f^{-1}(a)$ is the first time $f$ is strictly larger than $a$.
We also introduce the couple of random variables $(\mathcal{I}_1, \mathcal{I}_2)$ as follows,
\begin{equation}
    \mathcal{I}_1
:=
    \mathcal{Y}_1^{\natural}\big(\mathcal{Y}_2^{-1}(1)^-\big),
\qquad
    \mathcal{I}_2
:=
    \left(1 - \mathcal{Y}_2\big(\mathcal{Y}_2^{-1}(1)^-\big)\right )
    \times
    \frac{\mathcal{Y}_1\big(\mathcal{Y}_2^{-1}(1)\big) - \mathcal{Y}_1\big(\mathcal{Y}_2^{-1}(1)^-\big)}
    {\mathcal{Y}_2\big(\mathcal{Y}_2^{-1}(1)\big) - \mathcal{Y}_2\big(\mathcal{Y}_2^{-1}(1)^-\big)}. \label{defI}
\end{equation}

\noindent
We recall that $\cvloi$ denotes convergence in law under the annealed probability $\P$ as $t\to+\infty$.
We are now ready to state our first result.
\begin{theorem} \label{ADV}  Assume $0<\kappa<1$. We have,
$$
    \frac{\loX^*(t)}{t} \cvloi \mathcal{I} =: \max(\mathcal{I}_1,\mathcal{I}_2).
$$
\end{theorem}

\alex{Contrary to the recurrent case $\k=0$, we have no scaling property for the potential, and the diffusion $X$ cannot be localized in a single valley
as we can see in Theorem \ref{ththm}. However in the transient case we can make appear and use a renewal structure.
}
\\
\noindent
We now give an intuitive interpretation of this theorem, explaining the appearance of the Lévy process $(\mathcal{Y}_1, \mathcal{Y}_2)$.
\\
\noindent
First for any $s>0$, $\mathcal{Y}_1(s)$ is the limit of the sum of the first $\lfloor s e^{\kappa \phi(t)} \rfloor $ normalized (by $t$) local times taken specifically at the  $\lfloor s e^{\kappa \phi(t)} \rfloor$ first $h_t$-minima (see Proposition \ref{CVY1Y2} below).
Similarly, $\mathcal{Y}_2(s)$ is the limit of the sum of the
exit times of the $\lfloor s e^{\kappa \phi(t)} \rfloor$ first $h_t$-valleys, normalized (by $t$),
where an $h_t$-valley is a large neighborhood of an $h_t$-minimum.
For a rigorous definition of these $h_t$-valleys,
see Section \ref{sec2.2} and Figure \ref{fig1}.\\
So, by definition, $\mathcal{I}_1$ is the largest jump of the process $\mathcal{Y}_1$ before  the first time $\mathcal{Y}_2$ is larger than 1.
It can be interpreted as the largest (re-normalized) local time  among the local times at the $h_t$-minima visited by $X$ until time $t$
and from which $X$ has already escaped. That is to say, $\mathcal{I}_1$ is the limit of the random variable $\sup_{k \leq N_t-1} \loX(m_{k},t)/t$.
\\
$\mathcal{I}_2$ is a product of two factors: the first one, $\left (1 - \mathcal{Y}_2(\mathcal{Y}_2^{-1}(1)^-)\right )$,
corresponds to the (re-normalized) amount of time left to the diffusion $X$ before time $t$ after it has reached the ultimate visited  $h_t$-minimum $m_{N_t}$, \alex{that is, to $(t-H(m_{N_t}))/t$}. The second factor  corresponds to the local time of $X$ at this ultimate $h_t$-minimum $m_{N_t}$, that is to say $\mathcal{I}_2$ is the limit of $\lo(t,m_{N_t})/t$.  Intuitively $\mathcal{Y}_2$ is built from $\mathcal{Y}_1$ by multiplying each of its jumps by an independent copy of the variable $\mathcal{R}_{\kappa}$. Therefore this second factor can be seen as an independent copy of $1/\mathcal{R_{\kappa}}$ taken at the instant of the overshoot of $\mathcal{Y}_2$ which makes it larger than 1. Notice that this variable $\mathcal{R}_{\kappa}$ plays a similar role as $\mathcal{R}_{0}$ of Theorem \ref{cvloi}. Indeed as in the case $\kappa=0$, the diffusion $X$ is prisoner in the neighborhood of the last $h_t$-minimum  visited before time $t$.

We prove Theorem \ref{ADV} by showing first that portions of the trajectory of $X$ re-centered at the local $h_t$-minima, until time $t$,  are made (in probability) with independent parts. This has been partially proved in \cite{AndDev} but we have to improve their results and add  simultaneously the study of the local time.

Second, we prove that the supremum of the local time is, mainly,  a function of the sum of theses independent parts, which converges to a Lévy process. We now provide some details about this.

Recall that
$(\BP(s),\ s\geq 0)$ is defined as a continuous process, taking values in $\R_+$, with  infinitesimal generator given for every $x>0$ by
\begin{align*}
    \frac{1}{2}\frac{\dd^2}{\dd x^2}+\frac{\kappa}{2} \coth\left(\frac{\kappa}{2} x\right)\frac{\dd}{\dd x}.
\end{align*}
This process $\BP$ can be thought of as a
$(-\kappa/2)$-drifted Brownian motion $W_{\kappa}$ Doob-conditioned to stay positive, with the terminology of
\cite{Bertoin}, which is called Doob conditioned to reach $+\infty$ before $0$ in \cite{Faggionato}
(for more details, see Section 2.1 in \cite{AndDev}, where $\BP$ is denoted by $R$).
We call $\textnormal{BES}(3,\k/2)$ the law of $(\BP(s),\ s\geq 0)$.
That is, $(\BP(s),\ s\geq 0)$ is a $3$-dimensional $(\kappa/2)$-drifted Bessel process starting from $0$.
For any process $(U(t),\ \ t \in \R_+)$, we denote by
$$
\tau^{U}(a)
:=
\inf\{t>0,\  U_{}(t)=a\},
$$
the first time this process hits $a$, with the convention $\inf\emptyset=+\infty$.
For $a<b$,
$
\big(W_{\kappa}^{b} (s),\ 0 \leq s \leq \tau^{W_{\kappa}^b}(a)\big)
$
is defined as a $(-\kappa/2)$-drifted Brownian motion starting from $b$
and killed when it first hits $a$.
We now introduce some functionals of $W_{\kappa}$ and $\BP$, which already appeared in (\cite{AndDev}, Section 4.1):
\begin{eqnarray}
\label{eqDefF+G+}
    F^{\pm}(x)
& := &
    \int_0^{\tau^{\BP}(x)} \exp(\pm \BP_{}(s))\dd s,
\qquad
    x>0,
\\
    G^{\pm}(a,b)
& := &
    \int_0^{\tau^{W_{\kappa}^{b}}(a)} \exp\big(\pm W_{\kappa}^{b} (s)\big)\dd s,
\qquad a<b.
\label{eqDefF+G+Bis}
\end{eqnarray}
Let $0<\delta<1$, define
$$
    n_t:=\big\lfloor e^{\kappa \phi(t)(1+ \delta)} \big\rfloor,
\qquad
    t>0,
$$
which is, with large probability, an upper bound for $N_t$ as stated in Lemma \ref{lemtps}.

{Let $(S_j,R_j,{\bf e}_j, j \leq n_t)$ be a sequence} of i.i.d. random variables depending on $t$,
with $S_j$, $R_j$ and $\bf e_j$  independent, $S_1 \egloi F^{+}(h_t)+G^{+}(h_t/2, h_t)$, $ R_1 \egloi F^{-}(h_t/2)+\tilde F^{-}(h_t/2)$ and ${\bf e}_1\egloi \mathcal{E}(1/2)$ (an exponential random variable with parameter $1/2$),
where $\tilde F^{-}$ is an independent copy of $ F^{-}$ and $F^{+}$ is independent of $G^+$,
and $\egloi$ denotes equality in law.
Define $\ell_j:={\bf e}_j S_j$ and $\mathcal{H}_j:= \ell_j R_j $. Note that to simplify the notation, we do not make appear the dependence in $t$ in the sequel. Intuitively, $\ell_j$ plays the role of the local time at the $j$-th  positive $h_t$-minimum $m_j$ if $X$ escapes from the $j$-th $h_t$-valley before time $t$, that is, if $j<N_t$. Similarly, $\mathcal{H}_j$ plays the role of the time $X$ spends in
the $j$-th  $h_t$-valley before escaping from it. \\
\noindent Define the family of processes $(Y_1, Y_2)^t$ indexed by $t$, by
$$
    \forall s \geq 0,
\qquad
    (Y_1, Y_2)^t_s
=
    \big(Y_1^t(s),Y_2^t(s)\big)
:=
    \frac{1}{t} \sum_{j = 1}^{ \lfloor s e^{\kappa \phi(t)} \rfloor} (\ell_j, \mathcal{H}_j).
$$
Recall that $\cvloiS$ denotes convergence in law under $J_1$-Skorokhod topology.
Here is our next result.

\begin{proposition} \label{CVY1Y2}
Assume $0<\kappa<1$.
We have under $\P$, as $t\to+\infty$,
$$
    (Y_1, Y_2)^t\,
\cvloiS\,
    (\mathcal{Y}_1, \mathcal{Y}_2).
$$
\end{proposition}

Once this is proved, we check that
we can approximate, in law, the renormalized local time $\loX^*(t)/t$ by
a function of $(Y_1, Y_2)^t$. We obtain such  an expression in Proposition \ref{ProSupL}.
Then to obtain the limit claimed in Theorem \ref{ADV}, we prove the continuity (in $J_1$-topology) of the involved mapping  and apply a continuous mapping Theorem (see Section \ref{secCont}).

It appears that with this method we can also obtain some other asymptotics.
Indeed, we obtain in the following theorem the convergence in law of
the supremum of the local time of $X$ before $X$ hits the last $h_t$-minimum $m_{N_t}$ visited before time $t$,
of the supremum of the local time of $X$ before $X$ leaves
the last $h_t$-valley visited before time $t$ (the one around $m_{N_t}$)
approximately at time $H(m_{N_t+1})$, and of
the position of the favorite site.

\medskip

\begin{theorem} \label{ADV2}
Assume $0<\kappa<1$.
We have the following convergences in law under $\P$ as $t\to+\infty$,
\begin{eqnarray}
\label{Th15eq1}
    \frac{ \loX^{*}(H(m_{N_t + 1}))}{t}
& \cvloi &
    \Ymj\big(\Ymi(1)\big),
\\
    \frac{ \loX^{*}(H(m_{N_t}))}{t}
& \cvloi &
    \Ymj\big(\Ymi(1)^-\big)=\mathcal{I}_1.
\label{Th15eq2}
\end{eqnarray}
Let us call $F^*_t$ the position of the first favorite site, that is,
{$F^*_t:= \inf\{x \in \R,\, \loX(t,x)= \lo^*(t)\}$}. Then,
\begin{align}
    \frac{F^*_t}{X(t)}
\cvloi
    \mathcal{B}U_{[0,1]}+1-\mathcal{B},
\label{Th15eq3}
\end{align}
where $\mathcal{B}$ is a Bernoulli random variable with parameter  $\P(\mathcal{I}_1<\mathcal{I}_2)$,
and $U_{[0,1]}$ is a uniform random variable on $[0,1]$,
independent of $\mathcal{B}$.
\end{theorem}

\medskip

We remark that with probability one there is at most one point $x$ such that $\lo(t,x)= \lo^*(t)$ so $F_t^*$ is actually the favorite site.
Note that similar questions about favorite points for $X$ have been studied in the recurrent case $\k=0$ by \cite{Cheliotis}.

One question we may ask here is: what happens in the discrete case (that is, for RWRE),
or with a more general Lévy potential?

For RWRE, we expect a very similar behavior because the renewal structures which appear in both cases (RWRE and our diffusion $X$) are very similar (see \cite{EnSaZi}). The main difference comes essentially from the functional $\mathcal{R}_{\kappa}$, which should be replaced by a sum of exponentials of  simple random walks conditioned to remain positive (see \cite{EnSaZi2}, \cite{EnSaZi}).

For a more general Lévy potential, we have in mind for example a spectrally negative Lévy process
(diffusions in such potentials have been studied by \cite{Singh}). More work needs to be done, especially for the potential.
First, to obtain a specific decomposition of the Lévy's path (similar to what is done for the drifted Brownian motion in \cite{Faggionato}), and also to study the more complicated functional $\mathcal{R}_{\kappa}$ which is less known than in the Brownian case. This is a work in preparation by \cite{Vechambre}.

The rest of the paper is organized as follows.

\textit{In Section \ref{Section2}}, we recall the results of Faggionato on the path decomposition of the trajectories of $W_{\kappa}$. Also we recall from \cite{AndDev} the construction of specific $h_t$-minima which plays an important role in the appearance of independence, under $\P$, on the path of $X$ before time $t$.

\textit{In Section \ref{secontion3}}, we study the joint process of the hitting times of the $h_t$-minima $m_j$, $1\leq j\leq n_t$
and of local times at these $m_j$. We show that parts of the trajectory of $X$ are not important for our study, that is, we prove that the time spent outside the $h_t$-valleys,
and the supremum of the local time outside the $h_t$-valleys are negligible compared to $t$.
We then prove the main result of this section: Proposition \ref{Pro3.3}.
It shows that the joint process (exit times, local times) can be approximated in probability by $i.i.d$ random variables (which are the $\mathcal{H}_j$ and $\ell_j$).
This part makes use of some technical results inspired from  \cite{AndDev}, they are summarized in Section 6.

\textit{In Section \ref{Section4}}, we prove Proposition \ref{CVY1Y2}, and study the continuity of certain  functionals of $(\mathcal{Y}_1,\mathcal{Y}_2)$ which appear in the expression of the limit law $\mathcal{I}$. This section is independent of the other ones, we essentially prove a basic functional limit theorem and prepare to the application of continuous mapping theorem.

\textit{Section \ref{Section5}} is where we make appear the renewal structure in the problem we want to solve. In particular we show how the distribution of the supremum of the local time can be approximated by the distribution of some function of the couple $(Y_1,Y_2)^t$, the main step being Proposition \ref{ProSupL}.

\textit{Section \ref{Section6}} is a reminder of some key results and their extensions extracted from \cite{AndDev}.
For some of these results, sketch of proofs or complementary proofs are added in order for this paper to be more self-contained.

Finally, Section \ref{SectAppendix}
is a reminder of some estimates on Brownian motion, Bessel processes, and functionals of both of these processes.

\subsection{Notation}
In this section we introduce typical notation and tools for the study of diffusions in a random potential.


\noindent
 For any process $(U(t),\ \ t \in \R_+)$ we denote by $\lo_U$ a bicontinuous version of the local time of $U$ when it exists. Notice that for our main process $X$ we simply write $\lo$ for its local time. The inverse of the local time  for every $x \in \R$ is denoted by $\sigma_U(t,x):=\inf\{s>0,\ \lo_U(s,x) \geq t\}$ and in the same way $\sigma(t,x):=\sigma_X(t,x)$.
We also denote by $U^a$ the process $U$ starting from $a$, and by $P^a$ the law of $U^a$, with the notation $U=U^0$.
Now, let us introduce the following functional of $W_{\kappa}$,
$$
    A(r)
:=
    \int_0^r e^{\wk(x)}\dd x,\qquad r\in\R.
$$
We recall that since $\kappa>0$, $A_{\infty}:=\lim_{r\to+\infty}A(r)<\infty$ a.s.
As in \cite{Brox}, there exists a Brownian motion $(B(s),\ s\geq 0)$, independent of $W_{\kappa}$,
such that $X(t)=A^{-1}[B(T^{-1}(t))]$ for every $t\geq 0$, where
\begin{align}
T(r):=\int_0^r \exp\{-2\wk[A^{-1}(B(s))]\}\text{d}s,\qquad 0\leq r< \tau^B(A_{\infty}) \label{T}.
\end{align}
The local time of the diffusion $X$ at location $x$ and time $t$, simply denoted by $\loX(t,x)$, can be written as
(see \cite{Shi}, eq. (2.5))
\begin{equation}
    \loX(t,x)
=
    e^{- \wk(x)}\lo_B(T^{-1}(t),A(x)),
\qquad
    t>0,x\in\R.
\label{1.2}
\end{equation}
\noindent With this notation, we recall  the following expression of the hitting times of $X$,
\begin{equation}
    H(r)
=
    T\big[\tau^B(A(r))\big]
=
    \int_{-\infty}^r e^{-\wk(u)}\mathcal{L}_B[\tau^B(A(r)),A(u)]\dd u,
\qquad
    r\geq 0.\label{1.3}
\end{equation}



\section{Path decomposition and Valleys} \label{Section2}

\subsection{Path decomposition in the neighborhood of the $h_t$-minima $\mf_i$  }
We first recall some results for $h_t$-extrema of $\wk$.
Let
$$
V^{(i)}(x):=W_{\kappa}(x)-W_{\kappa}(\mf_i), \qquad x\in\R,\  i\in\N^*,
$$
which is the potential $\wk$ translated so that it is $0$ at the local minimum $\mf_i$. We also define
\begin{eqnarray}\label{eqDefTauj}
    \tau_i^-(h) & :=  & \sup \{s < \mf_i,\  V^{(i)}(s)=h\},
\qquad h>0,
\\
    \tau_i(h) & := & \inf \{s > \mf_i,\  V^{(i)}(s)=h\},
\qquad
     h>0.
\label{eqDefTaujBis}
\end{eqnarray}

The following result has been proved by \cite{Faggionato}
[for {\bf (i)} and {\bf (ii)}], and the last fact comes from the strong Markov property
(see also (\cite{AndDev}, Fact 2.1) and its proof).

\begin{fact}\label{Fact_Williams} (path decomposition  of $\wk$ around the $h_t$-minima $m_i$)\\
{\bf (i)} The truncated trajectories
$\big(V^{(i)}(\mf_i-s),\ 0 \leq s \leq \mf_i-\tau_i^-(h_t)\big)$,
$\big(V^{(i)}(\mf_i+s),\  0\leq s \leq \tau_i(h_t)-\mf_i\big)$, $i\geq 1$
are independent.
\\
{\bf (ii)} Let $(\BP(s),\ s\geq 0)$ be a process with law $BES(3,\k/2)$.
All the truncated trajectories
$\big(V^{(i)}(\mf_i-s),\ 0 \leq s \leq \mf_i-\tau_i^-(h_t)\big)$ for $i\geq 2$
and
$\big(V^{(j)}(\mf_j+s),\  0\leq s \leq \tau_j(h_t)-\mf_j\big)$ for $j\geq 1$
are equal in law to
$\big(\BP(s),\ 0 \leq s \leq \tau^{\BP}(h_t) \big)$.
\\
{\bf (iii)} For $i\geq 1$, the truncated trajectory
$\big(V^{(i)}(s+\tau_i(h_t)), \ s \geq 0\big)$
is independent of $\big(\wk(s),\ s\leq \tau_i(h_t)\big)$
and is equal in law to $\big(\wk^{h_t}(s),\ s\geq 0\big)$, that is, to a $(-\k/2)$-drifted Brownian motion starting from $h_t$.
\end{fact}

\medskip

\subsection{Definition of $h_t$-valleys and of standard $h_t$-minima \label{sec2.2} $\tm$, $j\in\N^*$ }~\\
We are interested in the potential around the $h_t$-minima $m_i$, $i\in\N^*$, in fact intervals containing at least $[\tau_i^-((1+\k)h_t), M_i]$. However, these valleys could intersect.
In order to define valleys which are well separated and i.i.d., we introduce the following notation.  This notation is used to define valleys of the potential around some $\tilde m_i$, which are thanks to Lemma \ref{CVs}  equal to the $m_i$ for $1\leq i \leq n_t$ with large probability.

\begin{figure}[h]
\begin{center}
\scalebox{1.2}{\input{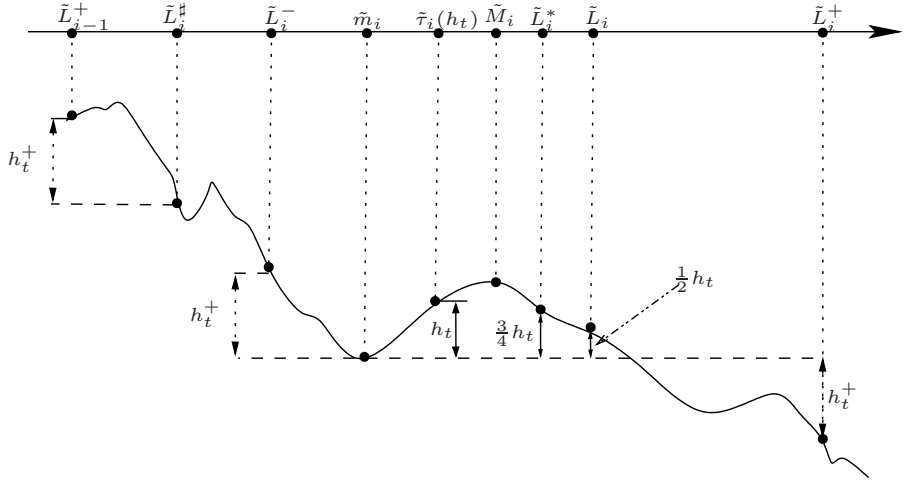}}
\caption{Schema of the potential between $\tilde L_{i-1}^+$ and $\tilde L_i^+$, in the case $\tilde L_i^\sharp <\tilde L_i^-$.}\label{fig1}
\end{center}
\end{figure}

\noindent Let
$$
    h_t^+
:=
    (1+ \kappa+ 2\delta) h_t.
$$
As in \cite{AndDev}, we define  $\tilde L_0^+:=0$, $\tilde m_0:=0$,
and recursively for $i \geq 1$ (see Figure \ref{fig1}),
\begin{align}
     \tilde L_i^{\sharp}
&:=
    \inf\{x>\tilde L_{i-1}^+,\ \wk(x)\leq\wk(\tilde L_{i-1}^+)-h_t^+   \},
\nonumber\\
    \tilde \tau_i(h_t)
& :=
    \inf \big\{x \geq   \tilde L_i^{\sharp},\ \wk(x)-\inf\nolimits_{[\tilde L_i^{\sharp},x]}\wk \geq h_t\big\},
\label{eqDefTaui1}
\\
    \tilde m_i
&:=
    \inf\big\{x \geq \tilde L_i^{\sharp},\ \wk(x)=\inf\nolimits_{[\tilde L_i^{\sharp},\tilde \tau_i(h_t)]}\wk\big\},
\nonumber\\
    \tilde L_i^{+}
&:=
    \inf\{x>\tilde \tau_i(h_t),\ \wk(x)\leq\wk(\tilde \tau_i(h_t))-h_t-h_t^+   \}.
\nonumber
\end{align}
\smallskip
We also introduce the following random variables for $i\in\N^*$:
\begin{align}
    \tilde M_i
&:=
    \inf\{ s>\tilde m_i,\ \wk(s)=\max\nolimits_{\tilde m_i \leq u\leq \tilde L_i^+} \wk(u)  \},
\nonumber\\
 \tilde L_i^* & :=
    \inf\{x>\tilde \tau_i(h_t),\ \wk(x)-\wk(\tilde m_i)= 3h_t/4  \}, \nonumber \\
    \tilde L_i
& :=
    \inf\{x>\tilde \tau_i(h_t),\ \wk(x)-\wk(\tilde m_i)= h_t/2  \},
\label{eqDefLiTilde}
\\
\tilde \tau_i(h)
& :=
    \inf \{s > \tilde m_i,\  \wk(x)-\wk(\tilde m_i)=h\}, \qquad h>0,
\label{eqDefTaui2}
\\
    \tilde \tau_i^-(h)
& :=
    \sup \{s < \tilde m_i,\  \wk(x)-\wk(\tilde m_i)=h\},
    \qquad h>0,
\label{eqDefTauiMoins}
\\
    \tilde L_i^{-}
&:=
    \tilde \tau_i^-(h_t^+). 
\nonumber
\end{align}


We stress that these random variables depend on $t$, which we do not write as a subscript to simplify the notation.
Notice also that $\tilde \tau_i(h_t)$ is the same in definitions \eqref{eqDefTaui1} and \eqref{eqDefTaui2} with $h=h_t$.
Moreover by continuity of $\wk$, $\wk(\tilde \tau_i(h_t))=\wk(\tilde m_i)+h_t$.
Thus, the $\tilde m_i$, $i\in\N^*$, are  $h_t$-minima,
since $\wk(\tilde m_i)=\inf_{[\tilde L_{i-1}^+,\tilde \tau_i(h_t)]}\wk$,
$\wk(\tilde \tau_i(h_t))=\wk(\tilde m_i)+h_t$ and
$\wk(\tilde L_{i-1}^+)\geq\wk(\tilde m_i)+h_t$.
In addition,
\begin{equation}\label{InegRvTilde1}
\tilde L_{i-1}^+<\tilde L_i^{\sharp}\leq \tilde m_i < \tilde \tau_i(h_t)<\tilde L_i^* <\tilde L_i <\tilde L_i^{+},
\qquad
    i\in\N^*,
\end{equation}
\begin{equation}\label{InegRvTilde2}
\tilde L_{i-1}^+\leq \tilde L_i^{-}< \tilde m_i < \tilde \tau_i(h_t)<\tilde M_i<\tilde L_i^{+},
\qquad
    i\in\N^*.
\end{equation}
Also by induction, the random variables $\tilde L_i^{\sharp}$, $\tilde \tau_i(h_t)$ and $\tilde L_i^{+}$, $i\in\N^*$ are stopping times
for the natural filtration of $(\wk(x), \ x\geq 0)$,
and so $\tilde L_i$, $\tilde L_i^*$,  $i\in\N^*$, are also stopping times.
Moreover by induction,
\begin{equation}\label{InegLiPremieresDescentes}
\begin{split}
    \wk(\tilde L_i^{\sharp})
=
    \inf_{[0, \tilde L_i^{\sharp}]} \wk,
\qquad
    \wk(\tilde m_i)
=
    \inf_{[0,\tilde \tau_i(h_t)]}\wk,
\\
    \wk(\tilde L_i^+)
=
    \inf_{[0, \tilde L_i^+]} \wk
=
    \wk(\tilde m_i)-h_t^+,
\end{split}
\end{equation}
for $i\in\N^*$.
We also introduce the analogue of $V^{(i)}$ for $\tilde m_i$ as follows:
$$
    \tilde{V}^{(i)}(x)
:=
    W_{\kappa}(x)-W_{\kappa}(\tilde m_i),
\qquad
    x\in\R,
    \ i\in\N^*.
$$

\noindent
We call $i\ th$ $h_t$-\textit{valley} the translated truncated potential
$\big(\tilde{V}^{(i)}(x),\ \tilde L_i^- \leq x \leq \tilde L_i \big)$, for $i\geq 1$.


The following lemma states that, with a very large probability, the first $n_t+1$ positive
$h_t$-minima $\mf_i$, $1 \leq i\leq n_t+1$,
coincide with the random variables $\tilde m_i$, $1 \leq i\leq n_t+1$.
We introduce the corresponding event
$\mV_t:=\cap_{i=1}^{n_t+1}\{\mf_i=\tilde m_i\}$.

\bigskip

\begin{lemma} \label{CVs}
Assume $0<\delta<1$.
There exists a constant $C_1>0$ such that
for $t$ large enough,
\[
    P\left(\overline{\mathcal V}_t\right)
\leq
    C_1 n_t e^{-\k h_t/2}
=
    e^{[-\k /2+o(1)]h_t}.
\]
Moreover, the sequence $\left(\big(\tilde{V}^{(i)}(x+\tilde {L}^+_{i-1}),\ 0\leq x \leq \tilde {L}^+_i-\tilde {L}^+_{i-1} \big),\ i \geq 1 \right)$, is i.i.d.
\end{lemma}

\noindent{\bf Proof:}
This lemma is proved in \cite{AndDev}: Lemma 2.3.
\hfill $\Box$

\bigskip

The following remark is used several times in the rest of the paper.

\medskip

\begin{remark} \label{RemEgaliteAvecouSansTilde}
On ${ \mV}_t$, we have for every $1\leq i\leq n_t$, $m_i=\tilde m_i$,
and as a consequence,
$\tilde{V}^{(i)}(x)=V^{(i)}(x)$, $x\in\R$,
$\tau_i^-(h)=\tilde \tau_i^-(h)$ and
$\tau_i(h)=\tilde \tau_i(h)$ for  $h>0$.
Moreover, $\tilde M_i=M_i$.
Indeed, $\tilde M_i$ is an $h_t$-maximum for $\wk$, which belongs to $[\tilde m_i,\tilde m_{i+1}]=[m_i, m_{i+1}]$ on
${ \mV}_t$, and there is exactly one $h_t$-maximum in this interval since the $h_t$-maxima and minima alternate,
which we defined as $M_i$, so $\tilde M_i=M_i$. So in the following, on ${ \mV}_t$, we can write these random variables with or without tilde.
\end{remark}

\section{Contributions for  hitting and local times \label{secontion3}}

\subsection{Negligible parts for hitting times\label{sectionEventB}}
\noindent \\
In the following lemma  we recall results of \cite{AndDev} which say, roughly speaking,
that the time spent by the diffusion $X$ outside the $h_t$-valleys is negligible compared to the amount of time spent
by $X$ inside the $h_t$-valleys.
This lemma also gives an upper bound for the number of $h_t$-valleys visited  before time $t$. Finally, it tells us that  with large probability, up to time $t$,
after first hitting the bottom $\tilde m_j$ of each $h_t$-valley $[\tilde L_j^-, \tilde L_j]$, $X$ leaves this $h_t$-valley on its right, that is on $\tilde L_j$, and that $X$ never backtracks in a previously visited $h_t$-valley. We define
$H_{x\to y}:= \inf\{s>H(x),\, X(s)=y\}-H(x)$ for any $x\geq 0$ and $y\geq 0$, which is equal to $H(y)-H(x)$ if $x<y$.
Let
$$
    U_0:=0,
\qquad
    U_i
:=
    H(\tilde L_i)-H(\tilde m_i)
=
    H_{\tilde m_i\to \tilde L_i},
\quad
    i \geq 1,
$$
$$
    \mathcal{B}_{1}(m)
:=
    \bigcap_{k=1}^{ m}\left\{0\leq H(\tilde m_{k})-\sum_{i=1}^{k-1}U_i< \tilde v_t \right\},
\qquad
    m\geq 1,
$$
where $\tilde v_t:=2t/ \log h_t$
and $\sum_{i=1}^0 U_i=0$ by convention. Finally, we introduce
$$
    \mathcal{B}_{2}(m)
:=
    \bigcap_{j=1}^{m} \Big\{
                H_{\tm \to \tL}<H_{\tm  \to  \tL^-},
                \ H_{\Lt_j\to \mt_{j+1}}<H_{\tL  \to \tilde{L}^{*}_j }
    \Big\},
\qquad
    m\geq 1.
$$

\smallskip

\begin{lemma} \label{lemtps} For any $\delta>0$ small enough, we have for all large $t$,
\begin{align}
    \P\big[H( \mt_1 ) \leq\tilde v_t \big]
\geq
    \P\big[\mathcal{B}_{1}(n_t)\big]
\geq
    1- C_2 v_t ,
    \label{TpsNeg}
\end{align}
with $v_t:=n_t \cdot (\log h_t) e^{- \phi(t)}=o(1)$
and $C_2>0$.
Moreover, there exists $C_3>0$ such that for large $t$,
\begin{eqnarray}
 \P\left(\mathcal{B}_{2}(n_t) \right)
 & \geq & 1-C_3 n_t e^{-\delta \kappa h_t},\label{HLLM} \\
 \P(N_t < n_t)
 & \geq & 1-  e^{- \phi(t) }. \label{Bas1}
\end{eqnarray}
\end{lemma}

\medskip

\noindent{\bf Proof:}
The first statement is Lemma 3.7 in  \cite{AndDev}.
The second one follows directly from Lemmata 3.2 and 3.3 in  \cite{AndDev}.
For the proof of \eqref{Bas1} see Lemma \ref{6.7}.
\hfill$\Box$

\subsection{Negligible parts for local times\label{sectionEventC}}
\noindent \\
We now provide estimates for the local time of $X$ at time $t$.
We first prove that the local time of $X$  outside the first $n_t$  $h_t$-valleys is negligible compared to $t$.
Second, we prove that for every $1\leq j\leq n_t$ the local time of $X$ inside the $h_t$-valley $[\tilde L_j^-, \tilde L_j]$
but outside a small neighborhood of $\tilde m_j$ is also negligible compared to $t$.

\subsubsection{Supremum of the local time outside the valleys}

\noindent \\ The aim of this subsection is to prove that at time $t$, the maximum of the local time outside the $h_t$-valleys
is negligible compared to $t$. More precisely,
let $f(t):= t   e^{[\k(1+3\d)-1]\phi(t)}$ and,
for $m\geq 1$,
\begin{align*}
\mathcal{B}_3^1(m)& :=\bigg\{\sup_{x\in[0, \tilde m_1]}\lx(H(\tilde m_1), x)\leq f(t)\bigg\}
\\
&\qquad\qquad
        \cap
        \bigcap_{j=1}^{m-1}
            \bigg\{\sup_{x\in[\tilde L_j, \tilde m_{j+1}]}\lx(H(\tilde m_{j+1}), x)\leq f(t)\bigg\},  \\
\mathcal{B}_3^2(m)& := \bigcap_{j=1}^{m-1}
            \left\{\sup\nolimits_{x\leq \tilde L_j}\left(\lx(H(\tilde m_{j+1}), x)-\lx\big(H\big(\tilde Lj\big), x\big)\right)\leq f(t)\right\}, \\
   \mathcal{B}_3(m)& :=   \mathcal{B}_3^1(m)\cap \mathcal{B}_3 ^2(m).
\end{align*}
This section is devoted to the proof of the following lemma.

\medskip

\begin{lemma} \label{LemmaProbaMaxLocHorsdeTOUTESVallees} \label{negloc1}
Assume that $\delta$ is small enough such that $\k(1+3\delta)<1$. There exists $C_5>0$ such that for any large $t$
\begin{align*}
    \P\left(\mathcal{B}_3(n_t)
    \right) \geq 1-C_5 w_t,
\end{align*}
with $w_t:=e^{-\k \delta \phi(t)}$.
\end{lemma}

\medskip

Its proof is based on Lemma \ref{LemmaMajorationMaxLocJusqueMontee1} below,
for which we introduce the following notation, depending only on the potential $\wk$:
\begin{eqnarray*}
    \tau_1^*(h)
& := &
    \inf\{u\geq 0,\ \wk(u)-\inf\nolimits_{[0,u]}\wk\geq h\},
\qquad
    h> 0,
\\
    m_1^*(h)
& := &
    \inf\{y\geq 0,\ \wk(y)=\inf\nolimits_{[0,\tau_1^*(h)]}\wk\},
\qquad
    h>0.
\end{eqnarray*}

Throughout the paper, $C_+$ (resp. $c_-$) denotes a positive constant that may grow
(resp. decrease) from line to line.

\medskip

\begin{lemma}\label{LemmaMajorationMaxLocJusqueMontee1} Assume that  $\k(1+3\delta)<1$. For large $t$,
\begin{equation}\label{eqLemmaTempsLocalExterieurVallees}
    \P\left(
        \sup\nolimits_{x\in[0, m_1^*(h_t)]}
        \loX[H(\tau_1^*(h_t)),x]
    >
         t   e^{[\k(1+3\d)-1]\phi(t)}
    \right)
\leq
    \frac{C_+ }{n_t e^{\k \delta \phi(t)}}.
\end{equation}
\end{lemma}

\medskip

\noindent{\bf Proof of Lemma \ref{LemmaMajorationMaxLocJusqueMontee1}:}
Thanks to \eqref{1.2} and \eqref{1.3} there exists a Brownian motion $(B(s),\ s\geq 0)$, independent of $\wk$, such that
\begin{equation}\label{eqLocalTimeUpToTau1}
    \loX[H(\tau_1^*(h_t)),x]
=
    e^{-\wk(x)}\mathcal{L}_B[\tau^B(A(\tau_1^*(h_t))),A(x)],  \qquad
  x\in\R.
\end{equation}
By the first Ray--Knight theorem (see e.g. \cite{RevYor}, chap.~XI),
for every $\a>0$,
there exists a Bessel processes $\BQ_2$
of dimension $2$
starting from $0$,
such that
$\lo_B(\tau^B(\a),x)$ is equal to $\BQ_2^2(\a-x)$ for every $x\in[0,\a]$.
Consequently, using \eqref{eqLocalTimeUpToTau1} and the independence of $B$ and $\wk$, there exists a $2$-dimensional Bessel process $\BQ_2$
such that
\begin{equation}\label{eqLocalTimeUpToTau1b}
    \loX[H(\tau_1^*(h_t)),x]
=
    e^{-\wk(x)}
    \BQ_2^2\big[A(\tau_1^*(h_t))-A(x)\big]
\qquad
    0\leq x\leq \tau_1^*(h_t).
\end{equation}
In order to evaluate this quantity, the idea is to say that loosely speaking, $\BQ_2^2$
grows almost linearly.
More formally, we consider the  functions $k(t):=e^{2\k^{-1}\phi(t)}$,
$a(t):=4\phi(t)$ and $b(t):=6 \k^{-1} \phi(t) e^{\k h_t}$,
and define the following events
\begin{eqnarray*}
    \B_0
& := &
   \bl{ \left\{A_\infty : =\int_0^{+ \infty} e^{\wk(u)}\dd u \leq k(t) \right\}, }
\\
    \B_1
& := &
    \big\{
        \forall u\in(0,k(t)],\ \BQ_2^2(u)
        \leq
        2 e  u \big[a(t)+4\log\log[e k(t)/u] \big]
    \big\},
\\
    \B_2
& := &
    \big\{
        \inf\nolimits_{[0,\tau_1^*(h_t)]}\wk\geq -b(t)
    \big\}.
\end{eqnarray*}
We know that $P(A_\infty\geq y)\leq C_+ y^{-\k}$ for $y>0$ since $2/A_\infty$ is a gamma variable of parameter $(\k,1)$
(see  \cite{Dufresne}, or \cite{Borodin} IV.48 p. 78), having a density equal to $e^{-x}x^{\k-1}{\mathds 1}_{\R_+}(x)/\Gamma(\k)$,
so
$
    P\big(\overline{\B_0}\big)
\leq
    C_+ k(t)^{-\k}
=
    C_+ e^{-2\phi(t)}
$.
Moreover,
$
    \P\big(\overline{\B}_1\big)
\leq
    C_+\exp[-a(t)/2]
=
    C_+e^{-2\phi(t)}
$
by Lemma \ref{lemmaTechniqueMajorationBessel2}.
Also we know that
$-\inf_{[0,\tau_1^*(h)]}\wk$, denoted by $-\beta$ in (\cite{Faggionato}, eq. (2.2)) is exponentially distributed
with mean $2\k^{-1} \sinh(\k h/2)e^{\k h/2}$ (\cite{Faggionato}, eq. (2.4)).
So for large $t$,
\begin{eqnarray*}
    P\big(\overline{\B}_2\big)
& = &
    P[-\inf\nolimits_{[0,\tau_1^*(h_t)]}\wk> b(t)]
\\
& = &
    \exp\big[-b(t) \k /(2\sinh(\k h_t/2)e^{\k h_t/2})\big]
\\
& \leq &
    e^{-2\phi(t)}.
\end{eqnarray*}
Now, assume we are on $\B_0\cap \B_1\cap \B_2$.
Due to \eqref{eqLocalTimeUpToTau1b}, we have for every $    0\leq x< \tau_1^*(h_t)$,
since $0< A(\tau_1^*(h_t))-A(x)\leq A_\infty\leq k(t)$,
\begin{align}
&
    \loX[H(\tau_1^*(h_t)),x]
\nonumber\\
& \leq
    e^{-\wk(x)}
    2 e [A(\tau_1^*(h_t))-A(x)]
    \big\{a(t)+4\log\log\big[e k(t)/[A(\tau_1^*(h_t))-A(x)]\big] \big\}.
\label{InegTempsLocalEntre0etTau1a}
\end{align}
We now introduce
$$
    f_i
:=
    \inf\{u\geq 0,\ \wk(u)\leq -i\}=\tau^{\wk}(-i),
\qquad
    i\in\N,
$$
and let $0\leq x< \tau_1^*(h_t)$. There exists $i\in\N$ such that $f_i\leq x < f_{i+1}$.
Moreover, we are on $\B_2$,
so $i\leq b(t)$.
Furthermore, $x < f_{i+1}$, so $\wk(x)\geq -(i+1)$ and then $e^{-\wk(x)}\leq e^{i+1}=e^{-\wk(f_i)+1}$.
All this leads to
\begin{eqnarray}\label{InegDifferenceA}
    e^{-\wk(x)}[A(\tau_1^*(h_t))-A(x)\big]
& = &
    e^{-\wk(x)} \int_x^{\tau_1^*(h_t)}e^{\wk(u)}\dd u
\nonumber\\
& \leq &
   e \int_{f_i}^{\tau_1^*(h_t)}e^{\wk(u)-\wk(f_i)}\dd u.
\end{eqnarray}
To bound this, we introduce the event
$$
    \B_3
:=
    \bigcap_{i=0}^{\lfloor b(t)\rfloor }
    \bigg\{
        \int_{f_i}^{\tau_1^*(h_t)}e^{\wk(u)-\wk(f_i)}\dd u
        \leq
        e^{(1-\k)h_t}
        b(t)n_t e^{\k \d \phi(t)}
    \bigg\}.
$$
We now consider
$\tau_1^*(u,h_t):=\inf\{y\geq u, \ \wk(y)-\inf_{[u,y]}\wk\geq h_t\}\geq \tau_1^*(h_t)$ for $u\geq 0$.
We have
$$
    E
    \bigg(
        \int_{f_i}^{\tau_1^*(h_t)}
        e^{\wk(u)-\wk(f_i)}\dd u
    \bigg)
\leq
    E
    \bigg(
        \int_{f_i}^{\tau_1^*(f_i,h_t)}
        e^{\wk(u)-\wk(f_i)}\dd u
    \bigg)
=
\beta_0(h_t),
$$
by the strong Markov property applied at stopping time $f_i$,
where we define
$
\beta_0(h):=
    E
    \left(
        \int_{0}^{\tau_1^*(h)}
        e^{\wk(u)}\dd u
    \right)
$.
By \eqref{6.1.1},
$\beta_0(h)\leq C_+ e^{(1-\k)h}$ for large $h$.
Hence for large $t$ by Markov inequality,
\begin{eqnarray*}
    P\big(\overline{\B}_3\big)
& \leq &
    \sum_{i=0}^{\lfloor b(t)\rfloor}
    P\bigg(
        \int_{f_i}^{\tau_1^*(h_t)}e^{\wk(u)-\wk(f_i)}\dd u
        >
        e^{(1-\k)h_t}
        b(t)n_t e^{\k \d \phi(t)}
    \bigg)
\\
& \leq &
    \frac{[b(t)+1]\beta_0(h_t)}{e^{(1-\k)h_t}b(t)n_t e^{\k \d \phi(t)}}
\leq
    \frac{C_+}{n_t e^{\k \d \phi(t)}}.
\end{eqnarray*}
Now, on $\cap_{j=0}^3 \B_j$,
\eqref{InegTempsLocalEntre0etTau1a} and \eqref{InegDifferenceA} lead to
\begin{eqnarray}
&&    \loX[H(\tau_1^*(h_t)),x]
\nonumber\\
& \leq &
    2e^{2+(1-\k)h_t}        b(t)n_t e^{\k \d \phi(t)}
    \big\{a(t)+4\log\log\big[e k(t)/[A(\tau_1^*(h_t))-A(x)]\big] \big\}.\hphantom{az}
\label{InegTempsLocalEntre0etTau1b}
\end{eqnarray}
We now consider only $0\leq x \leq m_1^*(h_t)$.
By definition of $\B_2$, $\inf_{[0, \tau_1^*(h_t)]}\wk  \geq -b(t)$, such that
\begin{eqnarray*}
    A(\tau_1^*(h_t))-A(x)
& = &
    \int_x^{\tau_1^*(h_t)} e^{\wk(u)}\dd u
\\
& \geq &
    \int_{m_1^*(h_t)}^{\tau_1^*(h_t)} e^{\wk(u)}\dd u
\\
& \geq &
    e^{-b(t)} [\tau_1^*(h_t)-m_1^*(h_t)]
\\
& \geq &
    e^{-b(t)}
\end{eqnarray*}
on the event $\cap_{i=0}^4 \B_i$ with
$
    \B_4
:=
    \{\tau_1^*(h_t)-m_1^*(h_t)\geq 1\}
$.
Since $m_1=m_1^*(h_t)$ and $\tau_1(h_t)=\tau_1^*(h_t)$ on $\{M_0\leq 0\}$ by definition of $h_t$-extrema, we have
\begin{eqnarray*}
    P\big(\overline{\B}_4\big)
& \leq &
    P(0<M_0<m_1)+ P[\tau_1(h_t)-m_1< 1]
\\
& \leq &
    C_+ h_t e^{-\k h_t}+P\big[\tau^{\BP}(h_t)-\tau^{\BP}(h_t/2)< 1\big]
\\
& \leq &
    C_+ h_t e^{-\k h_t}+C_+\exp[-(c_-)h_t^2]
\end{eqnarray*}
due to (\cite{AndDev}, eq. (2.8), coming from \cite{Faggionato}), Fact \ref{Fact_Williams} {\bf (ii)} and \eqref{3.10b}. \\
Now, we have
$
    e k(t)/[A(\tau_1^*(h_t))-A(x)]
\leq
    e k(t) e^{b(t)}
$
on $\cap_{i=0}^4 \B_i$,
and then, on this event, \eqref{InegTempsLocalEntre0etTau1b} leads to
\begin{eqnarray*}
    \loX[H(\tau_1^*(h_t)),x]
& \leq &
    2e^{2+(1-\k)h_t}b(t)n_t e^{\k \d \phi(t)}
    \big\{a(t)+4\log\log\big[e k(t) e^{b(t)}\big] \big\}.
\\
& \leq &
    C_+ t  \phi(t) e^{[\k(1+\d)-1]\phi(t)} e^{\k \d \phi(t)}
    h_t,
\end{eqnarray*}
since $\phi(t)=o(\log t)$, $h_t=\log t-\phi(t)$ and $n_t=\lfloor e^{\k(1+\d)\phi(t)}  \rfloor$.
We notice that for large $t$, $C_+ \phi(t) h_t\leq e^{\k \delta \phi(t)}$
since $\log \log t=o(\phi(t))$. Hence, for large $t$,
$$
    \loX[H(\tau_1^*(h_t)),x]
\leq
     t   e^{[\k(1+3\d)-1]\phi(t)},
$$
on $\cap_{i=0}^4 \B_i$
for every $0\leq x \leq m_1^*(h_t)$.
This gives for large $t$,
$$
    \P\left(
        \sup\nolimits_{x\in[0, m_1^*(h_t)]}
        \loX[H(\tau_1^*(h_t)),x]
    \leq
         t   e^{[\k(1+3\d)-1]\phi(t)}
    \right)
\geq
    \P\left(
        \cap_{i=0}^4 \B_i
    \right)
\geq
    1-\frac{C_+ }{n_t e^{\k \delta \phi(t)}},
$$
due to the previous bounds for $\P\big(\B_i \overline{}\big)$, $0\leq i \leq 4$.
This proves the lemma.
\hfill$\Box$



\bigskip

With the help of the previous lemma, we can now prove Lemma \ref{LemmaProbaMaxLocHorsdeTOUTESVallees}.

\bigskip

\noindent {\bf Proof of Lemma \ref{LemmaProbaMaxLocHorsdeTOUTESVallees}:}
The method is to do a coupling, similarly as in the proof of Lemma 3.7 of \cite{AndDev}.
Recall the definition of $\tilde L_i^*<\tilde L_i<\tilde L_{i+1}^\sharp$ just above \eqref{eqDefTaui2}. Also, let
\begin{align}
    \tilde \tau_{i+1}^*(h_t)
& :=
    \inf\big\{u\geq \tilde L_i^*, \ \wk(u)-\inf\nolimits_{[\tilde L_i^*,u]}\wk \geq h_t\big\}\leq \tilde \tau_{i+1}(h_t),
\qquad
    i\geq 1,
\nonumber\\
    \tilde m_{i+1}^*(h_t)
& :=
    \inf\big\{u\geq \tilde L_i^*, \ \wk(u)=\inf\nolimits_{[\tilde L_i^*, \tilde \tau_{i+1}^*(h_t)]}\wk \big\},
\qquad
    i\geq 1,
\nonumber\\
    \B_5
& :=
    \cap_{i=1}^{n_t-1}
    \big\{\tilde \tau_{i+1}^*(h_t)=\tilde \tau_{i+1}(h_t)\big\},
\nonumber\\
\label{eqDefXi}
    X_i(u)
& :=
    X\big(u+H\big(\tilde L_i\big)\big),
\qquad
    X_i^*(u)
:=
    X\big(u+H\big(\tilde L_i^*\big)\big),
\qquad
    u\geq 0,\, i\geq 1.
\end{align}
Let $i\geq 1$.
By the strong Markov property,
$X_i$ and $X_i^*$ are diffusions in the potential $\wk$, starting respectively from $\tilde L_i$ and $\tilde L_i^*$.
We denote respectively by $\mathcal{L}_{X_i}$, $\mathcal{L}_{X_i^*}$, $H_{X_i}$ and $H_{X_i^*}$ the local times and hitting times of $X_i$ and $X_i^*$.
We have for every $x\geq \tilde L_i^*$,
\begin{eqnarray*}
    \loX(H(\tilde m_{i+1}),x)-\loX(H(\tilde L_i),x)
& \leq &
    \loX(H(\tilde m_{i+1}),x)-\loX(H(\tilde L_i^*),x)
\\
& = &
    \mathcal{L}_{X_i^*}\big(H_{X_i^*}(\tilde m_{i+1}),x\big).
\end{eqnarray*}
Consequently, on
$\B_5\cap \B_6$
with
$
    \B_6
:=
    \cap_{j=1}^{n_t-1}\big\{H_{X_j}(\tilde m_{j+1})<H_{X_j}\big(\tilde L_j^*\big)\big\}
$, for $1\leq i\leq n_t-1$,
\begin{eqnarray}
&&
    \sup_{x\in\R}
    \left(\lx(H(\tilde m_{i+1}), x)-\lx\big(H\big(\tilde L_i\big), x\big)\right)
\nonumber\\
& = &
    \sup_{\tilde L_i^*\leq x \leq \tilde m_{i+1}}
    \left(\lx(H(\tilde m_{i+1}), x)-\lx\big(H\big(\tilde L_i\big), x\big)\right)
\nonumber
\\
& \leq &
    \sup_{\tilde L_i^*\leq x \leq \tilde m_{i+1}}
    \mathcal{L}_{X_i^*}\big(H_{X_i^*}(\tilde m_{i+1}),x\big)
\nonumber
\\
& \leq &
    \sup_{\tilde L_i^*\leq x \leq \tilde m_{i+1}^*}
    \mathcal{L}_{X_i^*}\big(H_{X_i^*}(\tilde \tau_{i+1}^*(h_t)),x\big),
\label{RHSTempsLocalDeXiEtoile}
\end{eqnarray}
since $\tilde m_{i+1}^*=\tilde m_{i+1}\leq \tilde \tau_{i+1}(h_t)=\tilde \tau_{i+1}^*(h_t)$ on $\B_5$.
Now, notice that
the right hand side of \eqref{RHSTempsLocalDeXiEtoile}
is the supremum of the local times of $X_i^*-\tilde L_i^*$, up to its first hitting time of $\tilde \tau_{i+1}^*(h_t)-\tilde L_i^*$,
 over all locations in $[0, \tilde m_{i+1}^*-\tilde L_i^*]$.
Since $X_i^*-\tilde L_i^*$ is a diffusion in the potential
$\big(\wk(\tilde L_i^*+x)-\wk(\tilde L_i^*),\ x\in \R\big)$, which has on $[0,+\infty)$ the same law as $(\wk(x),\ x\geq 0)$
because $\tilde L_i^*$ is a stopping time for $\wk$,
the right hand side of \eqref{RHSTempsLocalDeXiEtoile} has the same law, under the annealed probability $\P$, as
$        \sup\nolimits_{x\in[0, m_1^*(h_t)]}
        \loX[H(\tau_1^*(h_t)),x]
$.
Consequently,
\begin{eqnarray}
&&
    \P\bigg(
        \bigcup_{i=1}^{n_t-1}
        \left\{
            \sup_{x\in\R}\left(\lx(H(\tilde m_{i+1}), x)-\lx(H(\tilde L_i), x)\right)
            >
            t   e^{[\k(1+3\d)-1]\phi(t)}
        \right\}
    \bigg)
\nonumber\\
& \leq &
    n_t
    \Big[
    \P\Big(
            \sup\nolimits_{x\in[0,m_1^*(h_t)]}\loX\big[H(\tau_1^*(h_t)),x\big]
        >
            t   e^{[\k(1+3\d)-1]\phi(t)}
    \Big)
    +\P\big(\overline{\B}_5\big)
    +\P\big(\overline{\B}_6\big)
    \Big]
\nonumber\\
& \leq &
    C_+ e^{-\k \delta \phi(t)}
\label{eqPourVallesAumoins1}
\end{eqnarray}
by Lemma \ref{LemmaMajorationMaxLocJusqueMontee1},
since
$
    \P\big(\overline{\B}_5\big)
\leq
    C_+ n_t h_t e^{-\k h_t}
$
by \eqref{6.2.1},
$
    \P\big(\overline{\B}_6\big)
\leq
    \P\left(\overline{\mathcal{B}_{2}}(n_t)\right)
\leq
    C_3 n_t e^{-\delta \kappa h_t}
$
by \eqref{HLLM}
and since $\phi(t)=o(\log t)$.
Notice that, as before, $\tilde m_1=m_1=m_1^*(h_t)$ on $\mV_t\cap \{M_0\leq 0\}$.
Finally,
\begin{align*}
    \P\Big(
        \sup\limits_{x\in[0, \tilde m_1]}\lx(H(\tilde m_1), x)
        >
         t   e^{[\k(1+3\d)-1]\phi(t)}
    \Big)
& \leq
    \frac{C_+}{e^{\k \delta \phi(t)}}
    +P\big(\overline{\mV_t}\big)
    +P(0<M_0<m_1)
\\
& \leq
    \frac{C_+}{e^{\k \delta \phi(t)}}
\end{align*}
also by Lemma \ref{LemmaMajorationMaxLocJusqueMontee1}, Lemma \ref{CVs}, and since
$P(0<M_0<m_1)\leq C_+  h_t e^{-\k h_t}$
due to \eqref{6.2.2}.
This and \eqref{eqPourVallesAumoins1} prove the lemma.
\hfill$\Box$


\subsubsection{Local time inside the valley $\big[\tilde L_j^-, \tilde L_j\big]$ but far from $\tilde m_j$}

\noindent \\ We introduce for $t>0$ and $j\geq 1$,
\begin{align}
    r_t
:=
    C_0 \phi(t),
\qquad
    {\Dt}_j
:=
    [\tm-r_t, \tm+r_t],
\label{Dj}
\end{align}
where  $C_0>0$ is a constant that can be chosen as large as needed.
We also define
$$
    \mathcal{B}_{4}(m)
:=
    \bigcap_{j=1}^{m}\Bigg\{ \sup_{  x \in \overline{\Dt_j}  \cap  [\tilde \tau_j^-(h_t^+),\Lt_j]  }
    \Big(\loX\big(H(\Lt_{j}),x\big)-\loX\big(H(\mt_{j}),x\big) \Big) <t e^{-2 \phi(t)} \Bigg\}
$$
for $m\geq 1$,
where $\overline{\Dt_j}$ is the complementary of ${\Dt}_j$.
Moreover, we recall that $\tilde L_j^-=\tilde \tau_j^-(h_t^+)$.

\bigskip

\begin{lemma}\label{LemmaTempsLocalBordVallees} \label{negloc2}
There exists $C_6>0$ such that if $C_0$ is large enough, for large $t$,
$$
    \P\big[\mathcal{B}_4(n_t)\big]
\geq
    1- C_6 n_t e^{-2\phi(t)}.
$$
\end{lemma}

\bigskip

\noindent{\bf Proof:} 
Let $j\in[1,n_t]$. Throughout the rest of the paper, for $y\in\R$,
we denote by $\P_y^{\wk}$ the law of $X$ starting from $y$ instead of $0$, conditionally on $\wk$.
As we are interested in the local time at $x$ after $X$ reaches $\mt_{j}$ we work under \bl{$\P_{\tilde m_j}^{\wk}$}.
So first, thanks to \eqref{1.2} and \eqref{1.3}, under \bl{$\P_{\tilde m_j}^{\wk}$},
there exists a Brownian motion $(B(s),\ s\geq 0)$, independent of $\tilde V^{(j)}$, such that
\[
                \loX\big[H(\tilde L_j),x\big]
=
    e^{-\tilde V^{(j)}(x)}\mathcal{L}_B[\tau^B( A^j(\tilde L_j)), A^j(x)],
\qquad
    x\in\R,
\]
where $A^j(x):=\int_{\tilde m_j}^x e^{\tilde V^{(j)}(s)} \dd s$.
{Let $\tilde B^j(.) := B^j((A^j(\tL))^2 .)/A^j(\tL)$. By scaling, and because $B$ is independent from $W_{\kappa}$,
we notice that conditionally to $W_{\kappa}$, $\tilde B^j$ is a standard Brownian motion. Therefore, even if  $\wk$ appears in the expression of $\tilde B^j$,
$\tilde B^j$  is (probabilistically) independent of $\wk$. We still denote it by $B$ in the sequel to simplify the notation. With this notation, we have}
\begin{equation}\label{eqFormuleTempsLocalJusqueTempsSortieValleeScaling}
                \loX\big[H(\tilde L_j),x\big]
=
    e^{-\tilde V^{(j)}(x)}  A^j(\tilde L_j) \mathcal{L}_{B}[\tau^{B}(1),A^j(x)/A^j(\tilde L_j)],
\qquad
    x\in\R.
\end{equation}

In order to bound the factors $\mathcal{L}_{B}\big[\tau^{B}(1),.\big]$ and $A^j(\tilde L_j)$
in \eqref{eqFormuleTempsLocalJusqueTempsSortieValleeScaling},
we first introduce
\begin{equation}\label{eqDefEnsemblesE1etEZ2pourTempsLocalVallees}
    \B_1
:=
    \big\{\sup\nolimits_{u\in\R} \mathcal{L}_{B}[\tau^{B}(1),u]\leq e^{2\phi(t)}\big\},
\qquad
    \B_2
:=
    \big\{A^j(\tilde L_j)\leq  2 e^{h_t+2\phi(t)/\k}\big\}.
\end{equation}
We have
$
    \P\big({\overline{\B}_1}\big)
\leq
    5e^{-2\phi(t)}
$
for large $t$ by Lemma \ref{LemmaContinuiteEn0}  eq. \eqref{Diel}  and \eqref{Diel2}.
Moreover on $\mV_t$, we have
by Remark \ref{RemEgaliteAvecouSansTilde} and Fact \ref{Fact_Williams}  {\bf (ii)} and {\bf (iii)},
\begin{eqnarray*}
    A^j\big(\tilde L_j\big)
& \leq &
    \big[ \bl{\tilde \tau_j(h_t)-\tilde m_j} \big]e^{h_t}+\bl{\int_{\tilde \tau_j(h_t) }^{\tilde L_j} e^{\tilde V  ^{(j)}(s)} \dd s}
\\
& = &
    \big[ {\tau_j(h_t)-m_j} \big]e^{h_t}+{\int_{\tau_j(h_t) }^{L_j} e^{V^{(j)}(s)} \dd s}
\\
& \egloi &
    e^{h_t}\tau^{\BP}(h_t)+G^+(h_t/2,h_t),
\end{eqnarray*}
where $\BP$ has law $BES(3,\k/2)$ and is independent of $G^+(h_t/2,h_t)$,
which is defined in \eqref{eqDefF+G+Bis}, and with
$\bl{\tilde L_j}=\inf\{s>\bl{ \tilde \tau_j(h_t), \tilde V^{(j)}(s) }=h_t/2\}$
as defined in \eqref{eqDefLiTilde},
and $L_j:=\inf\{s>\tau_j(h_t), V^{(j)}(s) =h_t/2\}$.
Consequently,
\begin{eqnarray*}
    P\big(\overline{\B}_2\big)
& \leq &
    P\big(\tau^{\BP}(h_t)>e^{2\phi(t)/\k}\big)
    +
    P\big(G^+(h_t/2,h_t)> e^{h_t+2\phi(t)/\k}\big)
    +
    P\big(\overline{\mV_t}\big)
\\
& \leq &
    C_+ e^{-2\phi(t)}
\end{eqnarray*}
for large $t$ by
Lemma \ref{Lemma72}
eq. \eqref{bessel4}, Lemma \ref{EstimR} eq. \eqref{G}
and Lemma \ref{CVs},
and since $\phi(t)=o(\log t)$ and $\log \log t=o(\phi(t))$.

Now, we would like to bound the factor $e^{-\tilde V^{(j)}(x)}$ that appears in \eqref{eqFormuleTempsLocalJusqueTempsSortieValleeScaling}.
To this aim, let
\begin{eqnarray*}
    \B_3
& := &
    \big\{ \bl{\tilde \tau_j}[\k \Cct\phi(t)/8] \leq \tilde m_j +\Cct \phi(t)\big\},
\\
    \B_4
& := &
    \Big\{
        \inf_{[ \tau_j[\k  \Cct \phi(t)/8], \tau_j(h_t)]}  V^{(j)}\geq \k \Cct\phi(t)/16
    \Big\},
\end{eqnarray*}
{with $\tilde \tau_j$ and $\tilde \tau_j^-$ defined in \eqref{eqDefTaui2} and \eqref{eqDefTauiMoins}},
and $\tau_j$ and $\tau_j^-$ in \eqref{eqDefTauj} and \eqref{eqDefTaujBis}.
{First, using \eqref{6.2.4},
$
    P\big(\overline{\B}_3\big)
\leq
    C_+ e^{-[\k^2 \Cct\phi(t)]/(16\sqrt{2})}
\leq
    e^{-2\phi(t)}
$
if we choose $\Cct$ large enough.
Moreover Fact \ref{Fact_Williams} together with \eqref{3.10}
(applied with $h=\Cct\phi(t)$, $\alpha=\kappa/8$, $\gamma=\kappa/16$
and $\omega=h_t/(\Cct\phi(t))$, see also the remark at the end of Lemma \ref{Lemma72})
give
$
    P\big(\overline{\B}_4\big)
\leq
    2e^{- \k^2\Cct\phi(t)/16}
\leq
    e^{-2\phi(t)}
$
for large $t$.}  \\
We notice that $\inf_{[\tilde m_j+\Cct\phi(t), \tilde \tau_j(h_t)]} \tilde V^{(j)} \geq \k \Cct\phi(t)/16$
on $\B_3\cap \B_4\cap \mV_t$, since $\tau_j=\tilde\tau_j$ and $V^{(j)}=\tilde V^{(j)}$ on $\mV_t$ thanks to Remark \ref{RemEgaliteAvecouSansTilde}.
We prove similarly that
$$
    P\big(\overline{\B}_5 \big)
\leq
    C_+ e^{- \k^2\Cct\phi(t)/(16\sqrt{2})}+P\big(\overline{\mV_t}\big)
\leq
    2e^{-2\phi(t)},
$$
where
\begin{eqnarray*}
    \B_5
& := &
    \bigg   \{
        \inf_{[\tilde \tau_j^-(h_t), \tilde m_j-\Cct\phi(t)]} \tilde V^{(j)}\geq \k \Cct\phi(t)/16
    \bigg\},
\\
    \B_6
& := &
    \bigg\{
        \inf_{[\tilde \tau_j^-(h_t^+),\tilde \tau_j^-(h_t)]} \tilde V^{(j)}\geq h_t/2
    \bigg\}.
\end{eqnarray*}
Also by \eqref{6.2.3},
$
    P\big(\overline{\B}_6\big)
\leq
    e^{-\k h_t/8}
$.
We also know that
$
    \tilde V^{(j)}(x)\geq h_t/2\geq \kappa C_0 \phi(t) /16
$
for all $\tilde \tau_j(h_t)\leq x \leq \tilde L_j$ by definition of $\tilde L_j$, uniformly for large $t$.
Consequently on $\cap_{i=3}^6 \B_i\cap \mV_t$, for all
$x\in \overline{\mathcal{D}}_j\cap[\tilde \tau_j^-(h_t^+), \tilde L_j]$,
we have
$
    e^{-\tilde V^{(j)}(x)}
\leq
     e^{-\k \Cct\phi(t)/16}
$.

Hence on $\cap_{i=1}^6 \B_i\cap \mV_t$, we have under $\P_{\tilde m_j}^{\wk}$,
by \eqref{eqFormuleTempsLocalJusqueTempsSortieValleeScaling} and \eqref{eqDefEnsemblesE1etEZ2pourTempsLocalVallees},
$$
            \sup_{x\in \overline{\mathcal{D}_j}\cap[\tilde \tau_j^-(h_t^+), \tilde L_j]}
                \loX\big[H(\tilde L_j),x\big]
\leq
    2 t
    e^{(1+2/\k)\phi(t)}
     e^{-\k \Cct\phi(t)/16}
<
    t
    e^{-2\phi(t)},
$$
if we choose $\Cct$ large enough.
So, conditioning by $\wk$ and applying the strong Markov property at time $H(\tilde m_j)$, we get
\begin{eqnarray*}
&&
    \P
        \bigg[
            \sup_{x\in \overline{\mathcal{D}_j}\cap[\tilde \tau_j^-(h_t^+), \tilde L_j]}
            \big(
                \loX\big[H(\tilde L_j),x\big]
                -
                \loX\big[H(\tilde m_j),x\big]
            \big)
            <
                t e^{-2\phi(t)}
        \bigg]
\\
& \geq &
    \E\big(\P_{\tilde m_j}^{\wk}\big(\cap_{i=1}^6 \B_i\cap \mV_t\big)\big)
\geq 1-C_+ e^{-2\phi(t)}
\end{eqnarray*}
uniformly for large $t$ due to the previous estimates and thanks to Lemma \ref{CVs}.
This proves the lemma.
\hfill$\Box$


\subsection{Approximation of the main contributions}
\noindent \\
\alex{
In this section we give an approximation of the exit time of each $h_t$-valley $[\tilde L_j^-, \tilde L_j]$
and of the local time at the bottom $\tilde m_j$ of this $h_t$-valley for every $1\leq j \leq n_t$. }
More precisely, we make a link between the family
$\big((U_j:=H(\tL)-H(\tm), \loX(H(\tL),\tm)),1\leq j \leq n_t\big)$,
and the  i.i.d. sequence
$\big((\H_j, \ell_j),1\leq j \leq n_t \big)$ described in the introduction.


In the following, $F_1^+(h_t)$, $G^{+}(h_t/2, h_t)$, $F_2^-(h_t/2)$ and $F_3^-(h_t/2)$
denote independent r.v. with law respectively
$F^+(h_t)$, $G^{+}(h_t/2, h_t)$, $F^-(h_t/2)$ and $F^-(h_t/2)$, defined in \eqref{eqDefF+G+} and \eqref{eqDefF+G+Bis}.

\medskip

\begin{proposition} \label{Pro3.3}  For $\delta>0$ small enough (recall that $\delta$ appears in the definitions of $n_t$ and $h_t^+$),
there exist $d_1=d_1(\delta,\kappa)>0$ and $D_1(d_1)>0$
such that for large $t$,
possibly on an enlarged probability space,
there exists a sequence
$((S_j,R_j,{\bf e}_j),\ 1\leq j \leq n_t)$ of i.i.d. random variables depending on $t$, with $S_j$, $R_j$ and ${\bf e}_j$  independent
for every $j$
and $S_j \egloi F^{+}_1(h_t)+G^{+}(h_t/2, h_t)$,
$ R_j \egloi F^{-}_2(h_t/2)+F^{-}_3(h_t/2)$
and ${\bf e}_j\egloi \mathcal{E}(1/2)$ (exponential variable with mean $2$) such that
\begin{equation}\label{EqProposition35}
    \P\left(\cap_{j=1}^{n_t}
        \Big\{\big|U_j-\H_j\big| \leq \varepsilon_t \H_j,\ \big|\loX(H(\tilde L_j),\tilde m_j)-\ell_j\big| \leq \varepsilon_t \ell_j \Big\}
        \right)
\geq
    1-e^{-D_1 h_t},
\end{equation}
where
$\ell_j:=S_j {\bf e}_j$,
$\H_j:=R_j \ell_j$ and
$\varepsilon_t:=e^{-d_1 h_t}$.

\end{proposition}

\noindent The proof of the above proposition, which is in the spirit of the proofs of Propositions 3.4 and 4.4 in \cite{AndDev},
 makes use of the following lemma:

\bigskip

\begin{lemma} \label{lemX}
For $\delta>0$ small enough,
there exist constants $d_->0$ and
$D_->0$,
possibly depending on $\kappa$ and $\delta$, such that
the two following statements are true for $t>0$ large enough.\\
{\bf (i)} There exists a sequence $({\bf e}_j,\, 1\leq j \leq n_t)$
of i.i.d. random variables with exponential law of mean 2 and independent of $W_{\kappa}$, such that
 \begin{equation} \label{3.29}
    \P\bigg(\bigcap_{j=1}^{n_t}\Big\{|U_j-\mathbb{\tilde H}_j| \leq e^{-(d_-)h_t} \mathbb{\tilde H}_j,\ \loX(H(\tilde L_j),\tm)=\mathbb L_j \Big\} \bigg)
\geq
    1- e^{-(D_-)h_t},
\end{equation}
where
$\mathbb L_j:={\bf e}_j \int_{\tm}^{\tL} {e^{\tV(x)}}\dd x$,
$ \tilde {R}_j:= \int_{\tilde \tau_j^-(h_t/2)}^{\tilde \tau_j(h_t/2)}e^{-\tV(x)} {\dd x}$
and
$ \mathbb{\tilde H}_j:=\mathbb L_j \tilde {R}_j $
 for all $1\leq j \leq n_t$.
Moreover the random variables $(\mathbb L_j,\mathbb{\tilde H}_j)$, $1\leq j \leq n_t$, are i.i.d.
\\
{\bf (ii)}  Possibly on an enlarged probability space,
there exist random variables $R_j$ and $S_j$, $1\leq j\leq n_t$, such that all the random variables
$R_j$,  $S_j$, ${\bf e}_j$, $1\leq j \leq n_t$ are independent,
with
$S_j \egloi F^{+}_1(h_t)+G^{+}(h_t/2, h_t)$,
and
$ R_j \egloi F^{-}_2(h_t/2)+F^{-}_3(h_t/2)$
for every $1\leq j\leq n_t$, such that
\begin{equation}
    P\Bigg(\bigcap_{j=1}^{n_t}
        \left\{
            \left| \int_{\tm}^{\tL} {e^{\tV(x)}}\dd x - {S}_j\right| \leq e^{-(d_-)h_t} {S}_j,\,
            \tilde R_j=R_j
        \right\}
    \Bigg)
\geq
    1- e^{-(D_-)h_t},
\label{3.30}
\end{equation}
\end{lemma}

\bigskip

\noindent \textbf{Proof of Lemma \ref{lemX}:}
\textit{We start with {\bf (i)}}.
Recall that $\tilde m_j<\tilde L_j<\tilde m_{j+1}$ for every $j\geq 1$, e.g. by \eqref{InegRvTilde1}.
By the strong Markov property applied under $\P^{\wk}$ at stopping times $H(\tilde m_j)$, the random variables
$\big(U_j,\loX[H(\tL),\tm]\big),\, 1\leq j \leq n_t$, are independent under $\P^{\wk}$. By the same Markov property
and formulas \eqref{1.2} and \eqref{1.3},
the sequence
$(U_j,\loX[H(\tL),\tm],\, 1\leq j \leq n_t)$ is equal to the sequence
$(H_j(\tL), $ $ \lo_j[H_j(\tL), \tm],\, 1\leq j \leq n_t)$, where
\begin{align}
&
    H_j(\tL)
:=
    \int_{-\infty}^{\tL} e^{-\tV(u)}\mathcal{L}_{B^j}\big[\tau^{B^j}(A^j(\tL)),A^j(u)\big]\dd u,
\nonumber \\
&
    \lo_j[H_j(\tL),\tm]
=
    \mathcal{L}_{B^j}\big[\tau^{B^j}(A^j(\tL)),0\big],
\qquad
    A^j(u)
:=
    \int_{\tm}^{u} e^{\tV(x)}\dd x,
\quad
    u\in\R,
\label{DefAj}
\end{align}
with $(B^j, 1\leq j \leq n_t)$ a sequence of independent standard Brownian motions independent of $W_{\kappa}$, such that $B^j$ starts at $A^j(\tm)=0$ and is killed when it first hits $A^j(\tL)$. Recall that $\mathcal{L}_{B^j}$ denotes the local time of $B^j$.
{Define
$
    \mathcal{A}_j
:=
    \big\{ \max_{u<\tL^-}\mathcal{L}_{B^j}\big[\tau^{B^j}(A^j(\tL)),A^j(u)\big]=0 \big\}
$,
$1\leq j\leq n_t$.
By \eqref{6.3.2},   there exists $c_->0$
(possibly depending on $\kappa$ and $\delta$) such that   $\P\big(\cap_{j=1}^{n_t}\mathcal{A}_j\big) \geq 1- e^{-(c_-) h_t} $}
for large $t$. So for large $t$,
\begin{equation} \label{eqProbaHjhj}
    \P\bigg(\bigcap_{j=1}^{n_t} \left\{H_j(\tL)= \tilde h_j \right\} \bigg)
\geq
    1-  e^{-(c_-)  h_t},
\end{equation}
where
$$
    \tilde  h_j
:=
    \int_{\tL^-}^{\tL} {e^{-\tV(u)}}\mathcal{L}_{B^j}\big[\tau^{B^j}(A^j(\tL)),A^j(u)\big]\dd u,
\qquad
    1\leq j\leq n_t.
$$
We also notice that for every $1\leq j \leq n_t$,
$\big(\tilde h_j,\lo_j[H_j(\tL),\tm]\big)$ is measurable with respect to the $\sigma$-field generated by
$(\tV(x+\tilde {L}^+_{j-1}),\ 0\leq x < \tilde {L}^+_j-\tilde {L}^+_{j-1} )$ and $B^j$,
where by \eqref{InegRvTilde1} and \eqref{InegRvTilde2},
$\tilde {L}^+_{j-1}< \tL^-<\tilde m_j<\tilde L_j<\tilde L_j^+$.
Hence, the random variables
$\big(\tilde h_j,\lo_j[H_j(\tL),\tm]\big)$, $1\leq j \leq n_t$ are i.i.d under $\P$ by the second fact of Lemma \ref{CVs}.
For the same reason, $\big(\tilde R_j, A^j(\tilde L_j)\big)$, $1\leq j \leq n_t$ are also i.i.d.
\\
For $1\leq j \leq n_t$, let {$\tilde B^j(.) := B^j\big((A^j(\tL))^2 .\big)/A^j(\tL)$.
Notice that
\begin{equation}\label{EqEgaliteLocalTimes}
    \mathcal{L}_{B^j}\big[\tau^{B^j}(A^j(\tL)),A^j(u)\big]
=
    A^j(\tL) \mathcal{L}_{\tilde B^j}\big[\tau^{\tilde B^j}(1),A^j(u)/A^j(\tL)\big],
\qquad \tL^- \leq u \leq  \tL.
\end{equation}
Moreover by scaling and because $B^j$ is independent from $\wk$,  $\tilde B^j$ is, conditionally to $W_{\kappa}$,
a standard Brownian motion starting from $0$ and killed when it first hits $1$. Furthermore, even if $\wk$ appears in the expression of $\tilde B^j$, $\tilde B^j$ is independent of $\wk$.
 }
Then, let
\begin{equation}\label{eqDefej}
    {\bf e}_j
:=
    \mathcal{L}_{\tilde B^j}\big[\tau^{\tilde B^j}(1),0\big]
=
    \mathcal{L}_{B^j}\big[\tau^{B^j}\big(A^j\big(\tL\big)\big),0\big]/A^j\big(\tL\big).
\end{equation}
Notice that by the first Ray-Knight theorem,  ${\bf e}_j$ is exponentially distributed with mean $2$.
Since $\tilde B^j$ is independent of $\wk$, ${\bf e}_j$ is also independent of $\wk$.
Also, the sequence ${\bf e}_j$, $1\leq j\leq n_t$ is i.i.d. because the $B^j$ are independent
and the $\big(\tilde R_j, A^j(\tilde L_j)\big)$ are i.i.d., so
$(\mathbb L_j,\mathbb{\tilde H}_j)$, $1\leq j \leq n_t$, are also i.i.d.
Moreover, \eqref{EqEgaliteLocalTimes} leads to
\begin{equation}\label{eqExpressionLjej}
    \lo_j\big[H_j\big(\tL\big),\tm\big]
=
    A^j\big(\tL\big) \mathcal{L}_{\tilde B^j}\big[\tau^{\tilde B^j}(1),0\big]
=
    A^j\big(\tL\big) {\bf e}_j
=
    \mathbb L_j.
\end{equation}
Now, for small $\varepsilon>0$, thanks to Lemma \ref{Lemma63}, we have for large $t$,
\begin{eqnarray*}
    \P\bigg( \bigcap_{j=1}^{n_t} \left\{\left | \tilde h_j-A^j(\tL) {{\bf e}_j} \tilde {R}_j \right|
        \leq
        2e^{-(1-3\varepsilon)h_t/6} A^j(\tL) {{\bf e}_j}  \tilde {R}_j\right\}
    \bigg)
& \geq &
    1-\frac{C_+ n_t}{ e^{(c_-) \varepsilon h_t}}
\\
& \geq &
    1-\frac{C_+}{ e^{(c_-/2) \varepsilon h_t}},
\end{eqnarray*}
since $n_t=e^{o(1)h_t}$.
Finally, this, together with \eqref{eqProbaHjhj}, \eqref{eqExpressionLjej} and the equality of sequences at the start of this proof
show \eqref{3.29} for some $D_->0$ and $d_->0$. So {\bf (i)} is proved.

{\it We now prove {\bf(ii)}.}
The r.v. $\tilde A_j(\tilde L_j)=\int_{\tilde m_j}^{\tilde L_j}e^{\tilde V^{(j)}(x)}\dd x$ and $\tilde R_j$ are not independent, so we
want to replace them by r.v. having better independence properties.
Applying Lemma \ref{LemmaConstructionS2} with subscript $2$ replaced by $j$ for $1\leq j \leq n_t$
gives the existence of $R_j$ and $S_j$, independent and independent of ${\bf e_j}$,
having the law claimed in {\bf (ii)} and satisfying \eqref{EqApproxIntbySCase2}
with $2$ replaced by $j$.
This gives \eqref{3.30} since $n_t=e^{o(1)h_t}$.
The fact that we can build these $R_j$ and $S_j$ with the claimed independence properties follows from the
fact that $\big({\bf e_j}, \tilde R_j, \tilde A^j(\tilde L_j)\big), 1\leq j \leq n_t$ are i.i.d.
\hfill $\Box$


\bigskip

\noindent \textbf{Proof of Proposition \ref{Pro3.3}:}
\noindent
The existence and the law of the ${\bf e_j}$ come from Lemma \ref{lemX} {\bf (i)}.
The existence and the law of the $R_j$ and ${S}_j$, and the independence of
$R_j$, ${S}_j$, ${\bf e_j}$, $1\leq j \leq n_t$ come from Lemma \ref{lemX} {\bf (ii)}.
Moreover, by Lemma \ref{lemX} {\bf (i)} and {\bf (ii)},
there exist $d_1>0$ and $D_1>0$ such that
for large $t$,
\begin{eqnarray*}
&&
    \P\left(\cap_{j=1}^{n_t}\Big\{\big|U_j-{\bf e}_j {S}_j {R}_j \big|
                                \leq
                                \varepsilon_t {\bf e}_j {S}_j {R}_j ,
                                \
                                \big|\loX\big(H\big(\tilde L_j\big),\tm\big)-{\bf e}_j {S}_j \big|
                                \leq  \varepsilon_t {\bf e}_j {S}_j
                            \Big\}
    \right)
\\
& \geq &
    1-e^{-D_1 h_t},
\end{eqnarray*}
which proves \eqref{EqProposition35}. So Proposition \ref{Pro3.3} is proved.
\hfill $\Box$

\section{Convergence toward the Lévy process $(\mathcal{Y}_1, \mathcal{Y}_2)$ and continuity}\label{Section4}


\subsection{Preliminaries}

We begin this section by the convergence of certain repartition functions.
These {key results} are in the same spirit as the second part of Lemma 5.1 in \cite{AndDev}.

\medskip

\begin{lemma} \label{lemproba} Recall from Proposition \ref{Pro3.3} that $\ell_1:={\bf e}_1 { S}_1$
and  $\mathcal{H}_1:={\bf e}_1 { S}_1 R_1$.
Then for any $\varepsilon\in(0, 1/3)$,
\begin{align}
    \underset{t \rightarrow +\infty}{\lim}
&
    \sup_{x \in \left [ e^{-(1-2 \varepsilon)\phi(t)}, + \infty \right [}
    \ \left | x^{\kappa} e^{\kappa \phi(t)} \mathbb{P} \big( \ell_1 /t > x \big) - \mathcal{C}_2 \right |  = 0,
\label{cvmesure7.1}
\\
    \underset{t \rightarrow +\infty}{\lim}
&
    \sup_{y \in \left [ e^{-(1-3\varepsilon)\phi(t)}, + \infty \right [}
    \ \left | y^{\kappa} e^{\kappa \phi(t)} \mathbb{P} \big(  \mathcal{H}_1/t > y \big) - \mathcal{C}_2 \mathbb{E}\big[  (\mathcal{R}_{\kappa})^{\kappa} \big] \right |  = 0, \label{cvmesure9.1}
\end{align}
with $\mathcal{C}_2$ a positive constant (see below \eqref{cvmesure7}). \\
Moreover, for any $\alpha>0$, $e^{\kappa \phi(t)} \P(\ell_1/t \geq x, \H_1/t \geq y )$ converges uniformly when $t$ goes to infinity on $[\alpha,+ \infty [ \times [\alpha,+ \infty [$ to $\nu \left ( [x, + \infty[ \times [y, + \infty[ \right )$, where $\nu$ is defined in \eqref{eqDefMesureNU}.
\end{lemma}

\medskip

\noindent{\bf Proof:} Let $\varepsilon\in(0,1/3)$.\\
\textbf{Proof of \eqref{cvmesure7.1}:}
We first prove that, as $t\to+\infty$,  $x^{\kappa} e^{\kappa \phi(t)} \mathbb{P} \left ( S_1 /t > x \right )$ converges uniformly in $x \in \left [ e^{-(1-\varepsilon)\phi(t)}, + \infty \right [$ to a \bl{constant $c$}, that is, we prove that
\begin{align}
\underset{t \rightarrow +\infty}{\lim} \ \sup_{x \in \left [ e^{-(1-\varepsilon)\phi(t)}, + \infty \right [} \ \left | x^{\kappa} e^{\kappa \phi(t)} \mathbb{P} \left ( S_1 /t > x \right ) - \bl{c} \right | = 0. \label{0cvmesure1}
\end{align}
For that, with the change of variables $y = e^{(1-\varepsilon)\phi(t)} x$, we just have to prove that
\begin{align}
\underset{t \rightarrow +\infty}{\lim} \ \sup_{y \in \left [ 1, + \infty \right [} \ \left | y^{\kappa} e^{\kappa \varepsilon \phi(t)} \mathbb{P} \left ( S_1 / e^{h_t + \varepsilon \phi(t)} > y \right ) - \bl{c} \right | = 0, \label{0cvmesure2}
\end{align}
but this is equivalent to prove that for any function $f\ :\ ]0, + \infty[ \rightarrow [1, + \infty[$,
\begin{align}
\underset{t \rightarrow +\infty}{\lim} \ (f(t))^{\kappa} e^{\kappa \varepsilon \phi(t)} \mathbb{P} \left ( S_1 / e^{h_t + \varepsilon \phi(t)} > f(t) \right ) = \bl{c}. \label{0cvmesure3}
\end{align}
\noindent First by definition (see Proposition \ref{Pro3.3}),
$S_1$ can be written as the sum of two independent random variables, that we denote by
$F^{+}_1(h_t)$ and $G^{+}(h_t/2, h_t)$ for simplicity. That is,
\begin{align}
    S_1 /t &
=
    \left ( F^{+}_1(h_t)+G^{+}(h_t/2, h_t) \right )/t
=
    e^{-\phi(t)} \left ( e^{-h_t} F^{+}_1(h_t) + e^{-h_t} G^{+}(h_t/2, h_t)  \right ).
\label{cvmesure1}
\end{align}
Since we know the asymptotic behavior of the Laplace transforms of
$F^+(h_t)/e^{h_t}$ and $G^+(h_t/2, h_t)/e^{h_t}$, the proof of $(\ref{0cvmesure3})$ is similar to the proof of a Tauberian theorem. First by $\eqref{Ff2}$ and $\eqref{Gpf2}$ we have, using the independence of $F_1^+(h_t)$ and $G^+(h_t/2, h_t)$,
\begin{eqnarray}
    \forall \gamma > 0,
\qquad
    \ \omega_{f, t} (\gamma)
& :=  &
    \frac{1}{\gamma}
    \left({1 - \mathbb{E} \left [ e^{- \gamma S_1 / (f(t) e^{h_t + \varepsilon \phi(t)}) } \right ]}\right)
\nonumber\\
& \underset{t \rightarrow +\infty}{\sim} &
    \bl{c'} \gamma^{\kappa - 1} (f(t))^{-\kappa} e^{-\kappa \varepsilon \phi(t)},
\label{0cvmesure4}
\end{eqnarray}
where \bl{$c'=\Gamma(1-\kappa)2^\k/\Gamma(1+\kappa)$}.
Note that by Fubini, $\omega_{f, t}$ is the Laplace transform of the measure
$\dd  U_{f, t} (z) := \mathds{1}_{\mathbb{R}_+}(z) \mathbb{P} \left ( S_1 / (f(t) e^{h_t + \varepsilon \phi(t)}) > z \right ) \dd z$,
that is, $\omega_{f, t}(\gamma)=\int_0^{\infty} e^{-\gamma z}\dd U_{f, t}( z)$.
From $(\ref{0cvmesure4})$, we have
\[
    \forall \gamma > 0, \qquad \frac{\omega_{f, t} (\gamma)}{\omega_{f, t} (1)} \underset{t \rightarrow +\infty}{\longrightarrow} \gamma^{\kappa - 1}.
\]
We can now follow the same line as in the proof of a classical Tauberian theorem, making the link between a Laplace transform and the repartition function,  (see for example \cite{Feller2} volume 2, section XIII.5, Theorem 1, page 442), we can deduce that
\[ \forall z > 0, \qquad \frac{U_{f, t} ([0, z])}{\omega_{f, t} (1)} \underset{t \rightarrow +\infty}{\longrightarrow} \frac{z^{1-\kappa}}{\Gamma (2 - \kappa)}. \]
%
%
Then, e.g. as in the proof of Theorem 4 of the same reference page 446,
or using inequalities similar to those at the end of the proof of Lemma 5.1 in \cite{AndDev},
we deduce from the monotony of the densities of measures $U_{f, t}$  that
\[
    \forall z > 0, \qquad \frac{\mathbb{P} \left ( S_1 / (f(t) e^{h_t + \varepsilon \phi(t)}) > z \right )}{\omega_{f, t} (1)} \underset{t \rightarrow +\infty}{\longrightarrow} z^{-\kappa} \frac{1 - \kappa}{\Gamma (2 - \kappa)}.
\]
\noindent Considering this convergence with $z = 1$ we get exactly $(\ref{0cvmesure3})$ for
\bl{$c = c'(1 - \kappa) / \Gamma (2 - \kappa)=2^\k/\Gamma(1+\kappa)$},
so \eqref{0cvmesure1} follows.

Now, let $a_t := e^{\varepsilon \phi(t)}$.
For any $x>0$,
\[ x^{\kappa} e^{\kappa \phi(t)} \mathbb{P} \left ( {\bf e}_1 S_1 /t > x, \ {\bf e}_1 < a_t \right ) = 2^{-1} \int_0^{a_t} (x/u)^{\kappa} e^{\kappa \phi(t)} \mathbb{P} \left ( S_1 /t > x/u \right ) u^{\kappa} e^{-u/2} \dd u, \]
because ${\bf e}_1$ has law $\mathcal{E}(1/2)$ and is independent of $S_1$. \\
Taking $x$ arbitrary in $[ e^{-(1-2 \varepsilon)\phi(t)}, + \infty [$,
we have $x/u \in [ e^{-(1-\varepsilon)\phi(t)}, + \infty [$ for every $u \in ]0, a_t]$, so thanks to $(\ref{0cvmesure1})$ we get
\begin{align}
    \underset{t \rightarrow +\infty}{\lim} \ \sup_{x \in \left [ e^{-(1-2 \varepsilon)\phi(t)}, + \infty \right [}
    \, \left |
        x^{\kappa} e^{\kappa \phi(t)} \mathbb{P} \left ( {\bf e}_1 S_1 /t > x, \ {\bf e}_1 < a_t \right )
        -
        \bl{\frac{c}{2}} \int_0^{+\infty} \frac{u^{\kappa}}{e^{u/2}} \dd u
    \right |
= 0.
\label{cvmesure6}
\end{align}
Now for $t$ large enough such that
$\forall y \geq 1$, $y^{\kappa} e^{\kappa \phi(t)} \mathbb{P} \left ( S_1 /t > y \right ) < 2 \bl{c}$
(see \eqref{0cvmesure1}), we have for any $x > 0$,
\[
    \left | x^{\kappa} e^{\kappa \phi(t)} \mathbb{P} \big( {\bf e}_1 S_1 /t > x, \ {\bf e}_1 < a_t \big)
    - x^{\kappa} e^{\kappa \phi(t)} \mathbb{P} \big( {\bf e}_1 S_1 /t > x \big)
    \right |
\]
\begin{align}
& =
    x^{\kappa} e^{\kappa \phi(t)} \mathbb{P} \big( {\bf e}_1 S_1 /t > x, \ {\bf e}_1 \geq a_t \big)
\nonumber\\
& = 2^{-1} \int_{a_t}^{+\infty} x^{\kappa} e^{\kappa \phi(t)} \mathbb{P} \left ( S_1 /t > x/u \right ) e^{-u/2} \dd u
\nonumber\\
& = 2^{-1} \int_{a_t}^{+\infty} u^{\kappa} (x/u)^{\kappa} e^{\kappa \phi(t)} \mathbb{P} \left ( S_1 /t > x/u \right ) \mathds{1}_{x \leq u} e^{-u/2} \dd u
\nonumber\\
& + 2^{-1} \int_{a_t}^{+\infty} u^{\kappa} (x/u)^{\kappa} e^{\kappa \phi(t)} \mathbb{P} \left ( S_1 /t > x/u\right ) \mathds{1}_{x > u} e^{-u/2} \dd u
\nonumber\\
& \leq 2^{-1} e^{\kappa \phi(t)} \int_{a_t}^{+\infty} u^{\kappa} e^{-u/2} \dd u +  \bl{c} \int_{a_t}^{+\infty} u^{\kappa} e^{-u/2} \dd u.
\label{DerniereIneg}
\end{align}
For the second term in the inequality we used the fact that
$$
    (x/u)^{\kappa} e^{\kappa \phi(t)} \mathbb{P} \left ( S_1 /t > x/u\right ) < 2\bl{c}
$$
when $x \geq u$. Since $a_t = e^{\varepsilon \phi(t)}$, the right hand side of \eqref{DerniereIneg}
converges to $0$ when $t$ goes to infinity. Combining this with $(\ref{cvmesure6})$, we get
\begin{align}
\underset{t \rightarrow +\infty}{\lim} \ \sup_{x \in \left [ e^{-(1-2 \varepsilon)\phi(t)}, + \infty \right [} \ \left | x^{\kappa} e^{\kappa \phi(t)} \mathbb{P} \left ( {\bf e}_1 S_1 /t > x \right ) - 2^{-1} \bl{c} \int_0^{+\infty} u^{\kappa} e^{-u/2} \dd u \right | = 0, \label{cvmesure7}
\end{align}
and this is exactly $(\ref{cvmesure7.1})$ with
$
    \mathcal{C}_2
: =
    2^{-1} \bl{c} \int_0^{+\infty} u^{\kappa} e^{-u/2} \dd u
=
    2^\k \Gamma(\kappa+1)c
=
    4^\kappa
$.
\bigbreak
\noindent \textbf{Proof of \eqref{cvmesure9.1}:}
Let $\mu_{R_1}$ be the distribution of $R_1$. For any $y>0$,  $a>0$ and $t > 0$,
we have by independence of ${\bf e}_1 S_1$ and $R_1$,
\[
    y^{\kappa} e^{\kappa \phi(t)} \mathbb{P} \left ( {\bf e}_1 S_1 R_1 /t > y, \ R_1 < a \right ) = \int_0^{a} (y/u)^{\kappa} e^{\kappa \phi(t)} \mathbb{P} \left ( {\bf e}_1 S_1 /t > y/u \right ) u^{\kappa} \mu_{R_1}(\dd u).
\]
Taking $a = a_t = e^{\varepsilon \phi(t)}$ and $y$ arbitrary in $[ e^{-(1-3\varepsilon)\phi(t)}, + \infty [$,
we have
$y/u \in [ e^{-(1-2 \varepsilon)\phi(t)}, + \infty [$
for all $u \in ]0, a_t]$,
so thanks to $(\ref{cvmesure7})$, we get
\begin{align*}
&
    \underset{t \rightarrow +\infty}{\lim}
    \ \sup_{y \in \left [ e^{-(1-3\varepsilon)\phi(t)}, + \infty \right [}
    \ \left |
        y^{\kappa} e^{\kappa \phi(t)} \mathbb{P} \big( {\bf e}_1 S_1 R_1 /t > y, \ R_1 < a_t \big)
        -
        \mathcal{C}_2 \int_0^{a_t} u^{\kappa} \mu_{R_1}(\dd u)
    \right |
\\
& = 0,
\end{align*}
where we used
$
    \int_0^\infty u^{\kappa} \mu_{R_1}(\dd u)
=
    \mathbb{E}\left [ (R_1)^{\kappa} \right ]
\leq
    \mathbb{E}\left [ (\mathcal{R}_{\kappa})^{\kappa} \right ]
<
    \infty
$,
as explained in the following lines.
By definition (see before Proposition \ref{Pro3.3} and \eqref{eqDefF+G+})
and with $\widetilde W_{\kappa}^{\uparrow}$ an independent copy of $W_{\kappa}^{\uparrow}$,
$R_1$ is equal in law to
$
    \int_0^{\tau^{W_{\kappa}^{\uparrow}}(h_t /2)} e^{- W_{\kappa}^{\uparrow}(x)} \dd x
+
    \int_0^{\tau^{\widetilde W_{\kappa}^{\uparrow}}(h_t /2)} e^{- \widetilde W_{\kappa}^{\uparrow}(x)} \dd x
$,
which itself converges almost surely to $\mathcal{R}_{\kappa}$ (defined in \eqref{Rk}) when $t$ goes to infinity.
This also shows that for each $t$, $R_1$ is stochastically inferior to $\mathcal{R}_{\kappa}$,
which admits finite moments of any positive order by Lemma \ref{LTR}.
In particular the family $(R_1)_{t > 0}$ is bounded in all $L^p$ spaces,
and more precisely,
$
    \mathbb{E}\left [ (R_1)^{p} \right ]
\leq
    \mathbb{E}\left [ (\mathcal{R}_{\kappa})^{p} \right ]
<
    \infty
$
for every $p\in\mathbb R_+$.
So by the dominated convergence theorem,
$\int_0^{+\infty} u^{\kappa} \mu_{R_1}(\dd u)$ converges
to $\mathbb{E}\left [ (\mathcal{R}_{\kappa})^{\kappa} \right ]$ when $t$ goes to
infinity. Hence,
\[
    \underset{t \rightarrow +\infty}{\lim} \ \sup_{y \in \left [ e^{-(1-3\varepsilon)\phi(t)}, + \infty \right [}
    \ \left | y^{\kappa} e^{\kappa \phi(t)} \mathbb{P} \big( {\bf e}_1 S_1 R_1 /t > y, \ R_1 < a_t \big)
    - \mathcal{C}_2 \mathbb{E}\big[ (\mathcal{R}_{\kappa})^{\kappa} \big] \right | = 0.
\]
Finally, as the family $(R_1)_{t > 0}$ is bounded in all $L^p$ spaces,
$e^{\kappa \phi(t)}\int_{a_t}^\infty u^\kappa \mu_{R_1}(\dd u)$ converges to $0$ as $t\to+\infty$.
So we can proceed as before (as in \eqref{DerniereIneg}, integrating with respect to $R_1$ instead of ${\bf e}_1$
and using \eqref{cvmesure7.1} instead of \eqref{0cvmesure1})
to remove the event $R_1< a_t$ and we thus get
\begin{align}
\underset{t \rightarrow +\infty}{\lim} \ \sup_{y \in \left [ e^{-(1-3\varepsilon)\phi(t)}, + \infty \right [}
\ \left | y^{\kappa} e^{\kappa \phi(t)} \mathbb{P} \big( {\bf e}_1 S_1 R_1 /t > y \big)
- \mathcal{C}_2 \mathbb{E}\big[ (\mathcal{R}_{\kappa})^{\kappa} \big] \right | = 0, \label{cvmesure9}
\end{align}
which is $(\ref{cvmesure9.1})$.

{
We now prove the last assertion.
For any $x>0$, $y>0$, $a>0$ and $t > 0$, we have
\begin{eqnarray*}
&&
    e^{\kappa \phi(t)} \mathbb{P} \left ( {\bf e}_1 S_1 /t > x, \ {\bf e}_1 S_1 R_1 /t > y, \ R_1 < a \right )
\\
& = &
    \int_0^{a} e^{\kappa \phi(t)} \mathbb{P} \left ( {\bf e}_1 S_1 /t > x, \ {\bf e}_1 S_1 /t > y/u \right ) \mu_{R_1}(\dd u)
\\
& = &
    \int_0^{a \wedge (y/x) } e^{\kappa \phi(t)} \mathbb{P} \left ( {\bf e}_1 S_1 /t > y/u \right ) \mu_{R_1}(\dd u)
\\
&&
    \qquad
    +
    \int_{a \wedge (y/x) }^{a} e^{\kappa \phi(t)} \mathbb{P} \left ( {\bf e}_1 S_1 /t > x \right ) \mu_{R_1}(\dd u),
\\
& = &
    \frac1{y^{\kappa}} \int_0^{a \wedge (y/x) } e^{\kappa \phi(t)} (y/u)^{\kappa} \mathbb{P} \left ( {\bf e}_1 S_1 /t > y/u \right ) u^{\kappa} \mu_{R_1}(\dd u)
\\
&&
    \qquad +
    \frac1{x^{\kappa}} \int_{a \wedge (y/x) }^{a} e^{\kappa \phi(t)} x^{\kappa} \mathbb{P} \left ( {\bf e}_1 S_1 /t > x \right ) \mu_{R_1}(\dd u).
\end{eqnarray*}
Taking $a = a_t = e^{\varepsilon \phi(t)}$ and $x$, $y$ arbitrary in $[ \alpha, + \infty [$
(for some $\alpha > 0$), we have $(y/u,x) \in [ e^{-(1-2\varepsilon)\phi(t)}, + \infty [^2$,
$\forall u \in ]0, a_t]$ whenever $t$ is large enough, so,
thanks to $(\ref{cvmesure7})$ we get that
$e^{\kappa \phi(t)} \mathbb{P} \left ( {\bf e}_1 S_1 /t > x, {\bf e}_1 S_1 R_1 /t > y, \ R_1 < a_t \right )$
converges uniformly in $(x, y) \in [ \alpha, + \infty [ \times [ \alpha, + \infty [$ toward
\[
    \mathcal{C}_2x^{-\kappa}\mathbb{P}\big(\mathcal{R}_{\kappa}>y/x\big)
    +
    \mathcal{C}_2y^{-\kappa} \mathbb{E}\big((\mathcal{R}_{\kappa})^{\kappa}\un_{\mathcal{R}_{\kappa}\leq y/x}\big)
=
    \nu\big([x,+\infty[\times[y,+\infty[\big).
\]
Then as before we can remove the event $ \{ R_1 < a_t\}$
since $e^{\kappa\phi(t)} \mathbb{P} (R_1\geq a_t)\to  0$ as $t\to+\infty$ because the family $(R_1)_{t > 0}$ is bounded in all $L^p$ spaces,
which gives the last assertion of Lemma \ref{lemproba}.}
\hfill$\Box$

\subsection{Proof of Proposition \ref{CVY1Y2}}

\noindent \\
We start with the finite dimensional convergence.
We recall that $(Y_1, Y_2)^t_{s}$ is defined just before Proposition \ref{CVY1Y2}, and
$(\mathcal{Y}_1,\mathcal{Y}_2)$ before \eqref{eqDefMesureNU}.
We sometimes use the notation
$
    (\mathcal{Y}_1,\mathcal{Y}_2)_s
=
    (\mathcal{Y}_1(s),\mathcal{Y}_2(s))
$
and
$
    (Y_1,Y_2)^t_s
=
    (Y_1^t(s),Y_2^t(s))
$.

\medskip

\begin{lemma} For any $k \in \N$ and $s_i>0, i \leq k$,
$((Y_1, Y_2)^t_{s_i}, i \leq k)$ converges  in law
as $t$ goes to infinity to $((\mathcal{Y}_1,\mathcal{Y}_2)_{s_i},i\leq k)$.
\end{lemma}

\medskip

\noindent{\bf Proof:}
The proof is basic here, however we give some details  as we deal with a two dimensional walk which increments  depend on $t$ itself. As $Y_1^t(s)$ and  $Y_2^t({s})$ are sums of i.i.d sequences we only have to prove the convergence in law for the couple $(Y_1, Y_2)^t_{s}$ for any $s>0$.
For $b\geq 0$, we define $\big(Y_1^{> b},Y_2^{>b}\big)$, obtained from $(Y_1,Y_2)^t$ by keeping only the increments larger than $b$, that is,
$
    Y_1^{> b}(s)
:=
    \frac{1}{t} \sum_{j = 1}^{ \lfloor s e^{\kappa \phi(t)} \rfloor} \ell_j \un_{\ell_j/t > b }
$
and
$
    Y_2^{> b}(s)
:=
    \frac{1}{t} \sum_{j = 1}^{ \lfloor s e^{\kappa \phi(t)} \rfloor} \mathcal{H}_j \un_{\mathcal{H}_j/t > b }
$  for every $s\geq 0$ and $t>0$.
Also let $Y_i^{ \leq b}(s):=Y_i^t(s)-Y_i^{> b}(s)$ for $i\in\{1,2\}$.
We first prove that for any $s>0$,
\begin{align}
    \lim_{\varepsilon \rightarrow 0}
    \limsup_{t \rightarrow + \infty}
    \P\big(\big|\big|(Y_1^{ \leq \varepsilon},Y_2^{\leq \varepsilon})_s^t\big|\big|>\varepsilon^{1- \kappa(2-\kappa)} \big)
=
    0,
\label{infeps}
\end{align}
\noindent where for any $a=(a_1,a_2)\in \R^2$, $||a||:= \max(|a_1|, |a_2|)$,
with $(Y_1^{ \leq \varepsilon},Y_2^{\leq \varepsilon})_s^t=\big(Y_1^{ \leq \varepsilon}(s), Y_2^{ \leq \varepsilon}(s)\big)$
and $1- \kappa(2-\kappa)>0$ since $\kappa<1$.

\noindent {Let $\varepsilon>0$ and $s>0$. We now give an upper bound for the first moments of $Y_1^{ \leq \varepsilon}(s)$ and  $Y_2^{ \leq \varepsilon}(s)$.
Let $\eta>0$ be such that $\kappa -(1 - 3 \eta) < 0$. Applying Fubini, we have for large $t$,
\begin{eqnarray}
&&
    e^{\kappa \phi(t)} \E\left(\frac{\ell_1}{t} \un_{\ell_1 /t\leq \varepsilon }\right)
\nonumber\\
& = & e^{\kappa \phi(t)} \mathbb{E} \left [ \frac{{\bf e}_1 S_1}{t} \ \mathds{1}_{ {\bf e}_1 S_1 /t \leq \varepsilon } \right ]
\nonumber\\
& \leq &
    \int_0^{\varepsilon} e^{\kappa \phi(t)} \mathbb{P} \left ( {\bf e}_1 S_1 /t > x\right ) \dd x
\nonumber\\
& = &
    \int_0^{e^{-(1 - 2 \eta)\phi(t)}} e^{\kappa \phi(t)} \mathbb{P} \left ( {\bf e}_1 S_1 /t > x\right ) \dd x + \int_{e^{-(1 - 2 \eta)\phi(t)}}^{\varepsilon} e^{\kappa \phi(t)} \mathbb{P} \left ( {\bf e}_1 S_1 /t > x\right ) \dd x
\nonumber\\
& \leq &
    e^{(\kappa -(1 - 2 \eta) )\phi(t)} + \int_{e^{-(1 - 2 \eta)\phi(t)}}^{\varepsilon} x^{- \kappa} x^{\kappa} e^{\kappa \phi(t)} \mathbb{P} \left ( {\bf e}_1 S_1 /t > x\right ) \dd x.
\label{InegDeuxTermes}
\end{eqnarray}
The first term in \eqref{InegDeuxTermes} converges to $0$ when $t$ goes to infinity because $\kappa -(1 - 2 \eta) <-\eta< 0$.
Moreover, according to \eqref{cvmesure7.1}, for $t$ large enough, we have
\[
    \forall x \geq e^{-(1 - 2 \eta)\phi(t)},
\qquad
    x^{\kappa} e^{\kappa \phi(t)} \mathbb{P} \left ( {\bf e}_1 S_1 /t > x\right )
\leq
    2 \mathcal{C}_2.
\]
For such $t$, the second term in \eqref{InegDeuxTermes}  is less than
\[
    2 \mathcal{C}_2 \int_0^{\varepsilon} x^{- \kappa} \dd x = 2 \mathcal{C}_2 \frac{\varepsilon^{1 - \kappa}}{1-\kappa}.
\]
So, we get for large $t$,
\begin{align}
    e^{\kappa \phi(t)} \E\left(\frac{\ell_1}{t} \un_{\ell_1 /t
\leq
    \varepsilon }\right) \leq e^{(\kappa -(1 - 2 \eta) )\phi(t)} + C_+ \varepsilon^{1 - \kappa}.
\label{cvsub1}
\end{align}
Using the same method and applying this time \eqref{cvmesure9.1}, we get for large $t$,
\begin{align}
    e^{\kappa \phi(t)} \E\left(\frac{\H_1}{t} \un_{\H_1 /t
\leq
    \varepsilon }\right) \leq e^{(\kappa -(1 - 3 \eta) )\phi(t)} + C_+ \varepsilon^{1 - \kappa}.
\label{cvsub2}
\end{align}
We thus obtain
\begin{align} \E\left( Y_1^{ \leq \varepsilon}(s)\right) \leq s e^{(\kappa -(1 - 2 \eta) )\phi(t)} + C_+ s \varepsilon^{1- \kappa}, \label{cvsub3} \\
\E\left( Y_2^{ \leq \varepsilon}(s)\right) \leq s e^{(\kappa -(1 - 3 \eta) )\phi(t)} + C_+ s \varepsilon^{1- \kappa}, \label{cvsub4}
\end{align}
then a Markov inequality  leads to   \eqref{infeps} since $\kappa -(1 - 3 \eta) < 0$.

 }

The next step is to prove that $(Y_1^{>  \varepsilon},Y_2^{>\varepsilon})^t_s$ can be written
as the integral of a point process which converges to the desired limit. We have
$$
    \big(Y_1^{>  \varepsilon},Y_2^{>\varepsilon}\big)^t_s
=
    \big(Y_1^{>  \varepsilon}(s),Y_2^{>\varepsilon}(s)\big)
=
    \left(\int_{x> \varepsilon} \int_0^s x \mathcal{P}_t^1(\dd x,\dd v), \int_{x> \varepsilon} \int_0^s x \mathcal{P}_t^2(\dd x,\dd v) \right)
$$
where the measures  $\mathcal{P}_t^1$ and $\mathcal{P}_t^2$ are defined by
$
    \mathcal{P}_t^1
:=
    \sum_{i=1}^{+ \infty} \delta_{(t^{-1} \ell_i,e^{-\kappa \phi(t)}i) }
$
and similarly
$
    \mathcal{P}_t^2
:=
    \sum_{i=1}^{+ \infty} \delta_{(t^{-1} \H_i,e^{-\kappa \phi(t)}i) }
$.
Recall that $\mathcal{P}_t^1$ and $\mathcal{P}_t^2$ are not independent.
We now prove that $(\mathcal{P}_t^1,\mathcal{P}_t^2)$ converges to a Poisson point measure.
For that just use Lemma \ref{lemproba} together with  Proposition 3.1 in \cite{Resnick}
after discretization, it implies that $(\mathcal{P}_t^1,\mathcal{P}_t^2)$  converges weakly to the Poisson random measure denoted by $(\mathcal{P}^1,\mathcal{P}^2)$ with intensity measure given by $\dd s \times \nu $.

\noindent Then using that for any $\varepsilon>0$, and $T< + \infty$, on $[0,T)\times (\varepsilon,+ \infty) \times(\varepsilon,+ \infty) $ $\dd s \times \nu$ is finite, we have that $(Y_1^{>  \varepsilon},Y_2^{>\varepsilon})^t_s$ converges weakly to
$$
(\mathcal{Y}_1^{>  \varepsilon},\mathcal{Y}_2^{>  \varepsilon})_s:= \left(\int_{x> \varepsilon} \int_0^s x \mathcal{P}^1(\dd x,\dd v), \int_{x> \varepsilon} \int_0^s x \mathcal{P}^2(\dd x,\dd v) \right).
$$

\noindent We are left to prove that $(\mathcal{Y}_1^{>  \varepsilon},\mathcal{Y}_2^{>  \varepsilon})$
converges to  $(\mathcal{Y}_1,\mathcal{Y}_2)$ when $\varepsilon \downarrow 0$. This is a straightforward computation, that we detail for completeness.
Let $\nu_1([x, + \infty[):= \nu \left ( [x, + \infty[ \times \R_+ \right ) =\mathcal{C}_2 / x^{\kappa}$, we have
\begin{align*}
\E \left(\int_{x \leq \varepsilon} \int_0^s x \mathcal{P}^1(\dd x,\dd v)\right) = s \int_{x \leq \varepsilon} x \nu_1(\dd x) = C \varepsilon^{1-\kappa},
\end{align*}
Then a Markov inequality proves that for any $s>0$, the process
$\int_{x \leq \varepsilon} \int_0^s x \mathcal{P}^1(\dd x,\dd v)$ converges to  zero
(when $\varepsilon$ goes to zero) in probability.
The same is true for $\int_{x \leq \varepsilon} \int_0^s x \mathcal{P}^2(\dd x,\dd v)$, so we obtain that $(\mathcal{Y}_1^{>  \varepsilon},\mathcal{Y}_2^{>  \varepsilon})_s$ converges in probability
to $(\mathcal{Y}_1,\mathcal{Y}_2)_s$ when $\varepsilon \rightarrow 0$.
\hfill$\Box$

\noindent \\
We now prove the tightness of $(\mathcal{D}(Y_1,Y_2)^t)_t$, the family of measures induced by processes $(Y_1,Y_2)^t$.

\medskip

\begin{lemma} The family of laws $(\mathcal{D}(Y_1,Y_2)^t)_t$ is tight on $(D([0,+ \infty),\R^2),J_1)$.
\end{lemma}

\medskip

\noindent{\bf Proof:}
We only have to prove that the family law of the restriction of the process to the interval $[0,T]$, $((Y_1,Y_2)^t |_{[0,T]})_t$ is tight. To prove this we use the following  restatement of Theorem 1.8 in \cite{Billingsley} using Aldous's tightness criterion (see Condition 1, and equation (16.22) page 176 in \cite{Billingsley}) also used in \cite{Bovier} page 100.
We have to check the two following statements: \\
1)  for any $\varepsilon>0$, there exists $a$ such that for any $t$ large enough,
$$
    \P\Big(\sup_{s \in [0,T]}||(Y_1,Y_2)^t_s || \geq a\Big)
\leq
    \varepsilon.
$$
2) for any $\varepsilon>0$, and $\eta>0$ there exists $ \delta$, $0< \delta<T$ and  $t_0>0$
such that for $t>t_0$,
\[
    \P\big[  \omega((Y_1,Y_2)^t, \delta,T) \geq \eta\big]
\leq
    \varepsilon,
\]
with $\omega((Y_1,Y_2)^t, \delta,T):=\sup_{0\leq r \leq T} \omega((Y_1,Y_2)^t, \delta,T,r)$, and
\begin{eqnarray*}
&&
    \omega((Y_1,Y_2)^t, \delta,T,r)
\\
& := &
    \sup_{0 \vee (r-\delta) \leq u_1 < u < u_2 \leq (r+ \delta) \wedge T}
    \Big( \min\Big\{\big|\big| (Y_1,Y_2)^t_{u_2}-(Y_1,Y_2)^t_{u}\big|\big|,
\\
&&
    \qquad\qquad\qquad\qquad\qquad\qquad\qquad\qquad\qquad
    \big|\big| (Y_1,Y_2)^t_{u}-(Y_1,Y_2)^t_{u_1}\big|\big|\Big\}\Big).
\end{eqnarray*}
Also
$$
    \P(  v((Y_1,Y_2)^t,0,\delta,T) \geq \eta) \leq \varepsilon, \textrm{ and }  \P(  v((Y_1,Y_2)^t,T,\delta,T) \geq \eta)\leq \varepsilon,
$$
where $v((Y_1,Y_2)^t,u, \delta,T):=\sup_{(u- \delta)\vee 0 \leq u_1 \leq u_2 \leq (u+ \delta)\wedge T}\{|| (Y_1,Y_2)^t_{u_1}-(Y_1,Y_2)^t_{u_2}||\}$. \\

\noindent \underline{We first check 1)} since the process is monotone increasing,
\begin{align}
\P(\sup_{s \in [0,T]} ||(Y_1,Y_2)^t_s|| \geq a)=\P(||(Y_1,Y_2)^t_T|| \geq a) \leq \P(Y_1(T) \geq a) + \P(Y_2(T) \geq a). \label{majotriviale1}
\end{align}
Recall that $Y_1^{> b}$ is obtained from $Y_1$  where we remove the increments $\ell_j/t$ smaller than $b$  and $Y_1^{\leq b}=Y_1-Y_1^{> b}$. Define $N^{>b}_u:=\sum_{i=1}^{\lfloor u e^{\kappa \phi(t)} \rfloor} \un_{\ell_j/t > b}$. Let $0<\delta_1<1$. A Markov inequality yields
\begin{align}
\P\big(Y_1^t(T) \geq a\big)
& \leq  \P \left ( Y_1^{\leq 1}(T) \geq \frac{a}{2} \right ) + \P \left ( Y_1^{ > 1}(T) \geq \frac{a}{2} \right )
\nonumber \\
& \leq \frac{2}{a} \mathbb{E} \left [ Y_1^{ \leq 1}(T) \right ] + \frac{1}{a^{\delta_1}}\E \left ( N^{>1}_T\right )+ \P \left ( Y_1^{ > 1}(T) \geq \frac{a}{2}, N^{>1}_T \leq a^{\delta_1} \right ).
\label{majotriviale2}
\end{align}
On $\{N^{>1}_T \leq a^{\delta_1} \}$ there is at most $a^{\delta_1}$ terms in the sum $Y_1^{ > 1}(T)$ so
\begin{eqnarray}
    \mathbb{P} \left ( Y_1^{ > 1}(T) >a/2, N^{>1}_T \leq a^{\delta_1} \right )
& \leq &
    \sum_{1 \leq i \leq a^{\delta_1}} \mathbb{P} \left ( \ell_i/t \geq (a^{1- \delta_1}/2) | \ell_i/t \geq 1 \right )
\nonumber \\
& \leq &
     a^{\delta_1} \mathbb{P} \left ( \ell_1/t \geq (a^{1- \delta_1}/2) | \ell_1/t \geq 1 \right )
\nonumber \\
& \leq &
    a^{\delta_1} 2 \frac{\mathcal{C}_2 e^{- \kappa \phi(t)}a^{-\kappa(1- \delta_1)} 2^{\kappa}}{\mathcal{C}_2e^{- \kappa \phi(t) }}
\nonumber\\
& = &
    2^{1 + \kappa} \ a^{\delta_1 - \kappa (1-\delta_1)},
\label{majotriviale3}
\end{eqnarray}
for all $t$ large enough thanks to $(\ref{cvmesure7.1})$ and $\delta_1$ such that $\delta_1- \kappa(1- \delta_1)<0$.

\noindent Also, notice that  for any $b>0$, $N^{>b}_T$ follows a binomial law with parameters
$\left ( \lfloor T e^{\kappa \Phi(t)} \rfloor, \P(\ell_1/t>b) \right )$.
So, using \eqref{cvmesure7.1} again and \eqref{cvsub3}, we obtain for  $t$ large enough,
\begin{align}
\E(N^{>b}_T) \leq 2 \mathcal{C}_2  T b^{-\kappa}, \qquad \mathbb{E} \left [ Y_1^{\leq b}(T) \right ] \leq 2 \mathcal{C}_2 T b^{1-\kappa}.  \label{majotriviale4}
\end{align}
Collecting $(\ref{majotriviale3})$, $(\ref{majotriviale4})$ and $(\ref{majotriviale2})$, we get the existence of $t_1>0$ such that
\begin{align}
\underset{a \rightarrow +\infty}{\lim} \ \sup_{t \geq t_1} \ \P(Y_1(T) \geq a) = 0. \label{majotriviale6}
\end{align}

\noindent The same arguments holds for $Y_2$ (using $(\ref{cvmesure9.1})$ instead of $(\ref{cvmesure7.1})$ and $(\ref{cvsub4})$ instead of $(\ref{cvsub3})$) so $(\ref{majotriviale6})$ also holds for $Y_2$ instead of $Y_1$. We conclude the proof of 1) by putting $(\ref{majotriviale6})$ and its analogous for $Y_2$ in $(\ref{majotriviale1})$. \\

\noindent \\ \underline{We now check 2)}
\noindent  We first write, as usual,
\begin{eqnarray*}
&&
    \big\{\omega\big((Y_1,Y_2)^t, \delta,T\big) \geq \eta \big\}
\\
& \subset &
    \big\{\omega\big((Y_1^{\leq b},Y_2^{ \leq b})^t, \delta,T\big) \geq \eta/2 \big\}
    \cup
    \big\{\omega\big(\big(Y_1^{>b},Y_2^{>b}\big)^t, \delta,T\big) \geq \eta/2 \big\}.
\end{eqnarray*}
\noindent \\ For $Y_.^{\leq b}$, we have
$$
    \P\big[\omega\big(\big(Y_1^{\leq b},Y_2^{ \leq b}\big)^t, \delta,T\big) \geq \eta/2 \big]
\leq
    \P\big[\omega\big(Y_1^{\leq b}, \delta,T\big) \geq \eta/2 \big]
    +
    \P\big[\omega\big(Y_2^{ \leq b}, \delta,T\big) \geq \eta/2 \big].
$$
Moreover, by positivity of the increments,
\begin{eqnarray}
&&
    \P \left( \omega\big(Y_1^{\leq b}, \delta,T\big) \geq \eta/2 \right )
\nonumber\\
& \leq &
    \P \left ( \cup_{k \leq \lfloor T/ 2 \delta \rfloor}
        \Big\{Y_1^{\leq b}\big((k+1) 2\delta\big)-Y_1^{\leq b}\big(k 2 \delta\big) \geq \eta/4Ò\Big\}
        \right )
\nonumber \\
& \leq &
    \sum_{k \leq \lfloor T/ \delta \rfloor}
    \P \left ( Y_1^{\leq b}\big((k+1) 2\delta\big)-Y_1^{\leq b}\big(k 2\delta\big) \geq \eta/4 \right )
    \label{modcont1}.
\end{eqnarray}
For any $k$, $Y_1^{\leq b}((k+1) 2\delta)-Y_1^{\leq b}(k 2\delta)$ is the sum of at most $\lfloor 2 \delta e^{\kappa \Phi(t)} \rfloor + 1$ {i.i.d.} random variables having the same law as $\ell_1 / t$. We get that for any integer $k$
\[ \P \left ( Y_1^{\leq b}((k+1)2 \delta)-Y_1^{\leq b}(k 2\delta) \geq \eta/4 \right ) \leq \P \left ( Y_1^{\leq b}(3\delta) \geq \eta/4 \right ) \leq 8 \mathcal{C}_2 \delta b^{1-\kappa} / \eta, \]
where the first inequality holds for $t$ large enough so that $2\delta e^{\kappa \Phi(t)} \geq 1$ and the second from the second expression in \eqref{majotriviale4} (replacing $T$ by $2\delta$). Combining with $(\ref{modcont1})$ we get for large $t$
\begin{align}
\P \left ( \omega(Y_1^{\leq b}, \delta,T) \geq \eta/2 \right ) \leq 24 \mathcal{C}_2 T (1 + 2\delta) b^{1-\kappa} / \eta,  \label{modcont2}
\end{align}
[note that $\delta$ will be chosen later (and will be less than $1$)]. $T$ and $\eta$ are fixed so we choose $b$ small enough so that the right hand side of $(\ref{modcont1})$ is less than $\varepsilon/4$.
A similar estimate can be proved for $\P( \omega(Y_2^{\leq b}, \delta,T) \geq \eta/2) $.
\noindent \\ For $Y_.^{> b}$, we have again
$$
\P(\omega((Y_1^{> b},Y_2^{ > b})^t, \delta,T) \geq \eta/2 ) \leq \P(\omega(Y_1^{> b}, \delta,T) \geq \eta/2 ) +\P(\omega(Y_2^{ > b}, \delta,T) \geq \eta/2 ).
$$
Since $Y_1^{>b}$ is piecewise constant with jumps larger than $b$, $\{\omega(Y_1^{>b}, \delta,T) > \eta/2\} $ implies that two jumps larger than $b$ for $Y_1^t$ occur in an interval smaller than $2 \delta$. That is $\{\omega(Y_1^{>b}, \delta,T) > \eta/2\} \subset \cup_{j=1}^{\lfloor T e^{\kappa \phi(t)}\rfloor}\cup_{i>j, (i-j)/e^{\kappa \phi(t)}\leq 2 \delta}^{\lfloor T e^{\kappa \phi(t)} \rfloor} \{ \ell_j \wedge \ell_{i}/t  >b \}$. Applying $(\ref{cvmesure7.1})$ for $t$ large enough,
\[
    \P \left(\cup_{j=1}^{\lfloor T e^{\kappa \phi(t)}\rfloor}
             \cup_{i>j, (i-j)/e^{\kappa \phi(t)}\leq 2 \delta}^{\lfloor T e^{\kappa \phi(t)} \rfloor}
             \Big\{ \ell_j \wedge \ell_{i}/t  >b \Big\} \right)
\leq
    8 \mathcal{C}_2^2 \delta T b^{-2 \kappa},
\]
which can be small choosing this time $\delta = \delta(b)$ properly.
Again the same argument can be used for $\omega(Y_2^{>b}, \delta,T)$.
To finish the proof, we have to deal with $v()$, as again our processes are increasing,
\begin{align*}
\P(v((Y_1,Y_2)^t,0,\delta,T) \geq \eta) \leq \P(||(Y_1,Y_2)^t_{\delta}|| \geq  \eta)
\end{align*}
we can then {proceed as for 1)}  decreasing the value of $\delta$ if needed, this also applies to $\P(v((Y_1,Y_2)^t,T,\delta,T) \geq \eta)$. \hfill$\Box$ \\

\noindent Putting together the two preceding lemmata we obtain Proposition \ref{CVY1Y2}.

\subsection{Continuity of some functionals of $(\mathcal{Y}_1,\mathcal{Y}_2)$ in $J_1$ topology\label{secCont}}
In this section, we study the continuity of some functionals which will be applied later to
$(Y_1, Y_2)^t$ and to the Lévy processes $(\mathcal{Y}_1,\mathcal{Y}_2)$.
\noindent \\
For our purpose, we are interested in the following mappings.
We have already mentioned the first two in the introduction:
\begin{align*}
& \begin{array}{llll}
J : &D \left ( \mathbb{R}_+, \mathbb{R} \right ) & \longrightarrow D \left ( \mathbb{R}_+, \mathbb{R} \right ) \\
   &f & \longmapsto f^{\natural}
\end{array}
& \begin{array}{llll}
I : & \left (  D \left ( \mathbb{R}_+, \mathbb{R} \right ), J_1 \right ) & \longrightarrow \left (  D \left ( \mathbb{R}_+, \mathbb{R} \right ), U \right ) \\
   &f & \longmapsto f^{-1}
\end{array}
\end{align*}
where $U$ denotes uniform convergence on every compact subset of $\R_+$.
Then we also need the compositions of these two:  \bl{for any positive $a$, let
\begin{align}
\begin{array}{llll}
    J_{I, a} :
&
    D \left ( \mathbb{R}_+, \mathbb{R}^2 \right ) & \longrightarrow \mathbb{R}
\\
&
    f = (f_1, f_2) & \longmapsto f_1^{\natural} \left ( f_2^{-1} (a) \right ),
\end{array}
\begin{array}{llll}
    J_{I, a}^- :
&
    D \left ( \mathbb{R}_+, \mathbb{R}^2 \right ) & \longrightarrow \mathbb{R}
\\
&
    f = (f_1, f_2) & \longmapsto f_1^{\natural} \left ( f_2^{-1} (a)^- \right ), \nonumber
\end{array}
\end{align}
$J_{I, a}$ (respectively $J_{I, a}^-$) produces the largest jump of $f_1$,
between $0$ and the time just  after  (respectively before) $f_2$ first reaches $(a,+\infty)$.
We also define $K_{I, a}$, $K_{I, a}^-$, $\tilde K_{I, a}$ and $\tilde K_{I, a}^-$ as follows.
}
$$
\begin{array}{llll}
    K_{I, a} :
&
    D \left ( \mathbb{R}_+, \mathbb{R}^2 \right ) & \longrightarrow \mathbb{R}
\\
&
    f = (f_1, f_2) & \longmapsto f_1 \left ( f_2^{-1} (a) \right ) ,
\end{array}
$$
\begin{equation}
\begin{array}{llll}
    K_{I, a}^- :
&
    D \left ( \mathbb{R}_+, \mathbb{R}^2 \right ) & \longrightarrow \mathbb{R}
\\
&
    f = (f_1, f_2) & \longmapsto f_1 \left ( f_2^{-1} (a)^- \right ),
\end{array}
\label{DefK}
\end{equation}
$$
\begin{array}{llll}
\tilde K_{I, a} : &D \left ( \mathbb{R}_+, \mathbb{R}^2 \right ) & \longrightarrow \mathbb{R} \\
   &f = (f_1, f_2) & \longmapsto f_2 \left ( f_2^{-1} (a) \right ),
\end{array}
$$
$$
\begin{array}{llll}
\tilde K_{I, a}^- : &D \left ( \mathbb{R}_+, \mathbb{R}^2 \right ) & \longrightarrow \mathbb{R} \\
   &f = (f_1, f_2) & \longmapsto f_2 \left ( f_2^{-1} (a)^- \right ).
\end{array}
$$
Finally, with $\Delta f_1(s):=f_1(s)-f_1(s^-)$, define  $F^*$ by
$$\begin{array}{llll}
    F^* :
&
    D \left ( \mathbb{R}_+, \mathbb{R}^2 \right ) & \longrightarrow \mathbb{R}
\\
&
    f = (f_1, f_2) & \longmapsto \inf \left \{ s \in \big(0, f_2^{-1} (1)\big),\
    \Delta f_1(s) = f_1^{\natural} \left ( f_2^{-1} (1)^- \right ) \right \}.
\end{array}
$$
We need this functional $F^*$ for the characterization of the favorite site.

\medskip

\begin{lemma} \label{continuityofJ}
 $J$ is continuous in the $J_1$ topology.
\end{lemma}

\medskip

\noindent{\bf Proof:}
This fact is basic. However, we have not found a proof in the literature, so we give  some details.
To prove the continuity on $D \left ( \mathbb{R}_+, \mathbb{R} \right )$, we only have to prove it for every compact subset of  $\mathbb{R}_+$, (see  \cite{Whitt} Theorem 12.9.1).   So let $f \in D \left ( \mathbb{R}_+, \mathbb{R} \right )$ and $T > 0$ at which $f$ is continuous,  let us prove that $J_T$ defined by $$\begin{array}{llll}
J_T : &D \left ( [0, T], \mathbb{R} \right ) & \longrightarrow D \left ( [0, T], \mathbb{R} \right ) \\
   &g & \longmapsto g^{\natural}
\end{array}$$
is continuous at the restriction $f_{|[0, T]}$.
Let $\varepsilon > 0$ and $g \in D \left ( [0, T], \mathbb{R} \right )$ such that $d_T (f_{|[0, T]}, g) \leq \frac{\varepsilon}{2}$. $d_T$ is the usual metric $d$ of the $J_1$-topology restricted to the interval $[0,T]$.
By definition  of $d_T$ there exists a strictly increasing continuous mapping of $[0,T]$ onto itself, $e : [0, T] \longrightarrow [0, T]$ such that
\begin{align*}
\sup_{s \in [0, T]} \left | e(s) - s \right | \leq \frac{\varepsilon}{2} \textrm{ and } \sup_{s \in [0, T]} \left | g \left ( e(s) \right ) - f_{|[0, T]}(s) \right | \leq \frac{\varepsilon}{2}.
\end{align*}
So for every $s \in [0, T]$ we have
\begin{align*}
\left | \Delta g \left ( e(s) \right ) -\Delta f_{|[0, T]} (s) \right | & = \left | \left ( g \left ( e(s) \right ) - g \left ( e(s)- \right ) \right ) - \left ( f_{|[0, T]} (s) - f_{|[0, T]} (s-) \right ) \right | \\
& \leq \left | g \left ( e(s) \right ) - f_{|[0, T]} (s) \right | + \left | g \left ( e(s)- \right ) - f_{|[0, T]} (s-) \right | \\
& \leq 2 \frac{\varepsilon}{2} = \varepsilon,
\end{align*}
where $\Delta h(s)=h(s)-h(s^-)$. This implies $d_T \left (J_T \left ( f_{|[0, T]} \right ), J_T(g) \right ) \leq \varepsilon$.
\hfill$\Box$

\medskip

\begin{lemma}  \label{lemContJI} \bl{Fix $a > 0$. The mappings $J_{I, a}^-$, $J_{I, a}$, $K_{I, a}^-$, $K_{I, a}$, $\tilde K_{I, a}^-$ and $\tilde K_{I, a}$} are continuous for $J_1$-topology at every couple $(f^1,f^2) \in D(\R_+,\R^2)$ such that
\begin{enumerate}
\item For any $\varepsilon>0$, $f^1$ and $f^2$ have a finite number of jumps greater than $\varepsilon$ on every compact subset of $\R_+^*$, \\
\item $f^2$ is strictly increasing, with a limit equal to $+ \infty$, \\
\item $f^2(0)=0$, \\
\item $f^2$ has a jump at $I(f^2)(a)$ and $\bl{f^2(I(f^2)(a)-)<a<f^2(I(f^2)(a))}$. \\
\end{enumerate}
\end{lemma}

\noindent{\bf Proof:}
This fact may also be known as we are looking at randomly stopped process, but once again we did not find what we need in the literature
(\cite{Silverstrov}, \cite{Whitt}).
\\
Let $(f^1_n,f^2_n)_n$ be a sequence of elements of $D(\R_+,\R)$ which converges to $(f^1,f^2)$ for the $J_1$ topology. To prove continuity, we prove that the sequence \gr{$(J_{I, a}^-(f^1_n,f^2_n))_n$} converges to\gr{ $J_{I, a}^-(f^1,f^2)$}, and the equivalent for $J_{I, a}$. \\
The first hypothesis {guaranties} that there exist neighborhoods of $I(f^2)\gr{(a)}$ for which {$f^1$ makes no jump greater than $1/4$ times its higher previous jump, that is to say there exists  $\delta \in ]0, I(f^2)\gr{(a)}[$  (notice that {$I(f^2)\gr{(a)}$ exists tanks to (2) and is positive thanks to $(3)$}) such that $f^1$ makes no jump greater than $J(f^1)(I(f^2)\gr{(a)} - \delta)/4$ on $[I(f^2)\gr{(a)} - \delta, I(f_2)\gr{(a)}[$ and on $]I(f^2)\gr{(a)}, I(f^2)\gr{(a)} + \delta]$. Note also that $J(f^1)$ is constant on $[I(f^2)\gr{(a)} - \delta, I(f_2)\gr{(a)}[$ and on $]I(f^2)\gr{(a)}, I(f^2)\gr{(a)} + \delta]$.} \\
Also $\delta$ can be made smaller (if needed) in such a way that $I(f^2)\gr{(a)} + \delta$ is a point of continuity of $(f^1, f^2)$ and $(f^1_n,f^2_n)_n$ for every $n\in \N$.
By hypothesis $d \left ( (f^1_n,f^2_n), (f^1, f^2) \right ) \longrightarrow_{n \rightarrow + \infty} 0$ so \[ d_n := d_{[0, I(f^2)\gr{(a)} + \delta]} \left ( (f^1_n,f^2_n)_{|[0, I(f^2)\gr{(a)} + \delta]}, (f^1, f^2)_{|[0, I(f^2)\gr{(a)} + \delta]} \right ) \longrightarrow_{n \rightarrow + \infty} 0,\]
where ${|[0, I(f^2)\gr{(a)} + \delta]}$ in index means restriction to $[0, I(f^2)\gr{(a)} + \delta]$.
Also by continuity of $J$ (see Lemma \ref{continuityofJ}) we also have $d \left ( J (f^1_n), J(f^1) \right ) \longrightarrow_{n \rightarrow + \infty} 0$ and therefore $$d_n' := d_{[0, I(f^2)\gr{(a)} + \delta]} \left ( \left ( J (f^1_n)\right )_{|[0, I(f^2)\gr{(a)} + \delta]}, \left ( J (f^1) \right )_{|[0, I(f^2)\gr{(a)} + \delta]} \right ) \longrightarrow_{n \rightarrow + \infty} 0.$$
Let $h^{-}$ (respectively $h^+$) be the largest jump of $f^1$ just before (resp. just after) $I(f^{2})\gr{(a)}$. By definition of $\delta$ we have $$h^{-} = J(f^{1}) \left ( I(f^{2})\gr{(a)} - \delta \right ), h^{+} = J(f^{1}) \left ( I(f^{2})\gr{(a)} + \delta \right ).$$
We have two cases, either $J(f^{1})$ is continuous at $I(f^{2})\gr{(a)}$ or it  makes a jump.

\noindent \textit{Case $J(f^{1})$ makes a jump}, in this case the size of the jump is $h^{+} - h^{-} > 0$.

\noindent Let
$\alpha = 8^{-1} \min \left ( h^{-}, \delta, 1 - f^{2} \left ( I(f^{2})\gr{(a)}^- \right ), f^{2} \left ( I(f^{2})\gr{(a)} \right ) - 1 \right ),$
and $n_0\in\N$ be such that  for any $n \geq n_0$,  $d_n < \alpha$ and $d_n' < \alpha$. $T  {= I(f^{2})\gr{(a)} + \delta}$, there exist two homeomorphisms $e_n, e_n' : [0, T] \longrightarrow [0, T]$ such that:

\begin{itemize}

\item $\sup_{s \in [0, T]} \left | e_n (s) - s \right | \leq d_n,$

\item
$
    \sup_{s \in [0, T]}
    \left |
        \left|
            \big(f^1_n\left ( e_n (s) \right ),f^2_n\left ( e_n (s) \right )\big)_{|[0, I(f^2)\gr{(a)} + \delta]}
        \right.
    \right.
$
\\
\hphantom{azertyazertyazerty}
$
    \left.
        \left.
            -
            \big(f^1(s), f^2(s)\big)_{|[0, I(f^2)\gr{(a)} + \delta]}
        \right |
    \right |_{\infty}
 \leq
    d_n.
$

\item $\sup_{s \in [0, T]} \left | e_n' (s) - s \right | \leq d_n',$

\item $\sup_{s \in [0, T]} \left| \left ( J (f^1_{ n}) \right )_{|[0, I(f^2)\gr{(a)} + \delta]} \left ( e_n' (s) \right ) - \left ( J (f^1) \right )_{|[0, I(f^2)\gr{(a)} + \delta]}(s) \right | \leq d_n'.$
\end{itemize}
\noindent The second inequality implies that for any $n \geq n_0$, $$f^2_n \left ( e_n \left ( I(f^2)\gr{(a)}^- \right ) \right ) < a < f^2_n \left ( e_n \left ( I(f^2)\gr{(a)} \right ) \right ),$$ so as we also have $f^2_n \left (  I(f^2_n)\gr{(a)}^- \right ) \leq a \leq f^2_n \left (   I(f^2_n)\gr{(a)}  \right )$ we get
\begin{align}
I(f^2_{ n})\gr{(a)} = e_n \left ( I(f^2)\gr{(a)} \right ). \label{b1}
\end{align}
The fourth point implies that for any  $n \geq n_0$,
\begin{align} J(f^1_{n}) \left ( e_n' \left ( I(f^2)\gr{(a)} - \frac{1}{2} \delta \right ) \right ) \geq J(f^1) \left (  I(f^2)\gr{(a)} - \frac{1}{2} \delta \right )-\alpha =h^--\alpha > \frac{1}{2} h^{-}. \label{b1.1}
\end{align}
The second point and the argument of the previous proof imply that for any $n \geq n_0$,
each jump of $f^1_{n}$ on $[e_n \left ( I(f^2)\gr{(a)} - \delta \right ), e_n \left ( I(f^2)\gr{(a)} \right )[$
is $2 \alpha$-close to a jump of $f^1$
on $[I(f^2)\gr{(a)} - \delta, I(f^2)\gr{(a)} [$,
but such jumps are less than $h^-/4$ because of the definition of $\delta$.
Thus, $f^1_{n}$ makes no jump larger than $h^{-}/2$
on the interval $[e_n \left ( I(f^2)\gr{(a)} - \delta \right ), e_n \left ( I(f^2)\gr{(a)} \right ))$.
Moreover, the increases of $e_n'$ and the first and third points imply that
\[
    e_n \left ( I(f^2)\gr{(a)} - \delta \right ) \leq e_n' \left ( I(f^2)\gr{(a)} - \delta/2 \right ) \leq e_n \left ( I(f^2)\gr{(a)}  \right ).
\]
So, combining this with \eqref{b1.1}, we get that $J(f^1_{n})$ is constant on the interval
$[e_n' \left ( I(f^2)\gr{(a)} - \delta/2 \right ), e_n \left ( I(f^2)\gr{(a)} \right ))$.

\noindent Now by definition of \gr{$J_{I, a}^-$}, with $(\ref{b1})$ and then collecting what have just done above yields
\begin{align} \forall n \geq n_0,\qquad  \gr{J_{I, a}^- \left ( (f^1_n,f^2_n) \right)} & = J(f^1_n) \left ( I(f^2_n)\gr{(a)}^-\right ) = J(f^1_n) \left ( e_n \left ( I(f^2)\gr{(a)} \right )^- \right ) \nonumber \\
& = J(f^1_n) \left ( e_n' \left ( I(f^2)\gr{(a)} - \delta/2 \right ) \right ). \label{JI1}
\end{align}

\noindent From definition of \gr{$J_{I, a}^-$} and the constantness of $J(f^1)$ on $[I(f^2)\gr{(a)} - \delta, I(f_2)\gr{(a)}[$ we also have
\begin{align}
\gr{J_{I, a}^- (f^1, f^2)} := J(f^{1}) \left ( I(f^{2})\gr{(a)}- \right ) = J(f^{1}) \left ( I(f^2)\gr{(a)} - \delta /2 \right ). \label{JI2}
\end{align}
Combining \eqref{JI1}, \eqref{JI2} and the fourth point
gives that, as $n$ goes to infinity,
\gr{$J_{I, a}^- \left ( (f^1_n,f^2_n) \right)$} converges to \gr{$J_{I, a}^- \left ( (f^1,f^2) \right)$}.

\noindent For \gr{$J_{I, a}$}, we prove in a similar way as above that $J(f^1_{n})$ is constant on $[e_n \left ( I(f^2)\gr{(a)} \right ), $ \\
$e_n' \left ( I(f^2)\gr{(a)} + \delta/2 \right ) ]$ so, as in \eqref{JI1} we have for $n$ large enough
\[ \gr{J_{I, a} \left ( (f^1_n,f^2_n) \right)} = J(f^1_n) \left ( e_n' \left ( I(f^2)\gr{(a)} + \delta/2 \right ) \right ), \]
which, combined with the analogous of \eqref{JI2}
\[ \gr{J_{I, a} (f^1, f^2)} = J(f^{1}) \left ( I(f^2)\gr{(a)} + \delta /2 \right ) \]
allows us to conclude, using the fourth point, that \gr{$J_{I, a} \left ( (f^1_n,f^2_n) \right)$} converges to \gr{$J_{I, a} \left ( (f^1,f^2) \right)$} as $n$ goes to infinity. Therefore, both \gr{$J_{I, a}^-$} and \gr{$J_{I, a}$} are continue at $(f^1,f^2)$.
{The continuity of the other functionals are proved similarly. }
\hfill$\Box$

\medskip

\begin{lemma} \label{lem4.7}
For any $(f^{1},f^{2})$ in $D(\R_+, \R^2)$  that satisfy the hypothesis of lemma \ref{lemContJI} and such that the sizes of the jumps of $f^{1}$ are all distinct,
$F^*$ is continuous at $(f^{1},f^{2})$ in the $J_1$ topology.
\end{lemma}

\medskip

\noindent{\bf Proof:}
The proof follows mainly the steps of Lemma $\ref{lemContJI}$, we keep the same notation. The jump which takes place at the instant $F^*(f^{1}, f^{2})$ has value $h^{-}$. With the additional hypothesis that the values of the jumps for  $f^{1}$ are all different we have unicity for the value $h^{-}$. Let us define $h'$, the second highest jump $f^{1}$ before instant $I(f^{2})(1)$.
With the additional condition that  $\alpha < \frac{1}{8} (h^{-} - h')$ we have with the same arguments as in the proof of the continuity of $J$ that for any $n \geq n_0$,  $f^{1}_n$ effectuates at  $e_n \left ( F^*(f^{1}, f^{2}) \right )$ a jump larger than $ h^{-} - 2 \alpha$, and larger than all the other jumps of  $f_{n}^1$ before $e_n (I(f^{2})(1)-) = I(f^{2}_n)(1)$ which are smaller than $ h' + 2 \alpha$.
So for $n \geq n_0$, the largest jump of $f^{1}$ before $I(f^2_{n})(1)$ is obtained for $e_n \left ( F^*(f^{1}, f^{2}) \right )$, that is to say for any  $n \geq n_0$, $$  F^* \left ( (f^{1}_n, f^{2}_n) \right ) = e_n \left ( F^*(f^{1}, f^{2}) \right ), $$
this implies $F^* \left ( (f^{1}_n, f^{2}_n) \right ) \longrightarrow_{n \rightarrow \infty} F^*(f^{1}, f^{2})$.
\hfill$\Box$

\section{Supremum of the Local time - and other functionals}\label{Section5}

\subsection{Supremum of the local time (proof of Theorem \ref{ADV})}

{
\noindent \\
First, notice that since the diffusion $X$ is almost surely transient to the right,
the random variable
$\sup_{x<0} \lo(+\infty,x)$ is $\mathbb P$-almost surely  finite. So almost surely,
$$
    \lim_{t \rightarrow +\infty } \sup_{x<0} \lo(t,x)/t=0.
$$
As a consequence, we only have to study the asymptotic behavior of $\sup_{x\geq 0} \lo(t,x)/t$} as $t\to+\infty$. \\
 We start with the proof of the following proposition, which makes a link between the supremum of the local time and the process $(Y_1,Y_2)^t$.

\medskip

\begin{proposition} \label{ProSupL}
Let $\alpha>0$.
For any $\varepsilon>0$ and large $t$,
\begin{align*}
    \mathcal{P}_1^{-} -v(\varepsilon,t)
\leq
    \P\left( \sup_{x  \geq 0} \loX(t,x)/t \leq \alpha  \right)
\leq
    \mathcal{P}_1^{+}+v(\varepsilon,t),
\end{align*}
where
$$
    \mathcal{P}_1^{\pm}
:=
    \P\left[
        \big(1- \bar{\mathcal{H}}_{ \mathcal{N}^{2\varepsilon}_t-1}\big)
        \frac{\bar{\ell}_{ \mathcal{N}^{2 \varepsilon}_t}-\bar{\ell}_{ \mathcal{N}^{2\varepsilon}_t-1} } {(\bar{\mathcal{H}}_{ \mathcal{N}^{2\varepsilon}_t}-\bar{\mathcal{H}}_{ \mathcal{N}^{2\varepsilon}_t-1})} \leq \alpha_t^{\pm},
        \ {\max_{ 1\leq j \leq  \mathcal{N}^{2\varepsilon}_t-1} \frac{\ell_j}{t} \leq \alpha_t^{\pm} }
    \right],
$$
and with
$\bar{\mathcal{H}}_k
    := Y_2^t\big(ke^{-\kappa \phi(t)}\big)= \frac{1}{t}\sum_{i=1}^k\mathcal{H}_i
$,
\
$
    \bar{\ell}_k
:=
    Y_1^t\big(ke^{-\kappa \phi(t)}\big)
=
    \frac{1}{t}\sum_{i=1}^k\ell_i
$ for any $k\in\N$,
\
$
    \mathcal{N}_t^{2\varepsilon}
:=
    \inf\big\{m \geq 1, \bar{\mathcal{H}}_{m} >1- 2\varepsilon \big\}
$,
\
$
    \alpha^{\pm}_t:=\alpha \big(1\pm (\log \log t)^{-1/2}\big)
$,
and $v$ is a positive function such that $\lim_{t \rightarrow + \infty} v(\varepsilon,t) \leq \textrm{const} \times \varepsilon^{ \kappa  \wedge (1- \kappa)}$.
\end{proposition}

\medskip

The proof of this proposition relies on the three following lemmata.
The first one deals with the local time at the $h_t$-minima for which the diffusion $X$  already escaped before time $t$.
The second deals with the local time at the last $h_t$-minimum $m_{N_t}$  in the remaining time before time $t$.
Finally the last one is a technical point.

\medskip

\begin{lemma} \label{Lem6.2} For any large $t>0$, $2 \leq k \leq n_t $,
any $x>0$
and $\gamma>0$ possibly depending on $t$, define the repartition function
$$
    F_{\gamma}(x)
:=
    \P\left(\max_{ 1\leq j \leq k-1} \loX\big(H\big(\Lt_j\big),\mt_j\big)  \leq \gamma t,\ \sum_{i=1}^{k-1}U_i \leq \bl{x t} \right).
$$
Then for large $t$, for all $2 \leq k \leq n_t $, $x>0$ and $\gamma>0$,
 \begin{align*}
F^-_{\gamma}(x) -e^{-D_1 h_t} \leq F_{\gamma}(x)  \leq F^+_{\gamma}(x) +e^{-D_1 h_t},
\end{align*}
where $F^{\pm}_{\gamma}(x):= \P\left(\max_{ 1\leq j \leq k-1} \ell_j \leq \gamma_t^{\pm} t, \sum_{i=1}^{k-1}\mathcal{H}_i \leq x_t^{\pm} t  \right)$ with $\gamma_t^{\pm}:= \gamma(1 \pm 2 \varepsilon_t)$, $x_t^{\pm}:=x(1 \pm  2 \varepsilon_t)$, $\varepsilon_t$ and $D_1$ are given in Proposition \ref{Pro3.3}.
\end{lemma}

\medskip

\begin{lemma} \label{lem6.3b} For any $t>0$, define for every $\gamma>0$ and $ 0<x<1 $ possibly depending on $t$,
$$
    f_{\gamma}(x)
:=
    E\left[
        \P_{\mt_1}^{W_{\kappa}}\Big(  \lo_{X'}(t(1-x),\mt_1) \leq \gamma t, H'(\Lt_1)>t(1-x), H'(\Lt_{1})< H'(\Lt_{1}^-) \Big)
   \right].
$$
For such  $t$, $\gamma$ and $x$, we also introduce
\begin{align*}
&
    \tilde f_{\gamma}(x)
\\
& :=
     E\left(\Pw_{\mt_{1}} \left(\sup_{y \in \Dt_{1} }  \lo_{X'}[t(1-x),y] \leq  \gamma t, H'(\Lt_{1})>t(1-x), H'({\Lt_1})<H'({\Lt^-_1}) \right) \right),
\end{align*}
with $\Dt_{1}$ defined in \eqref{Dj}.
Here $X'$ is an independent copy of $X$ starting at ${\mt_1}$, and the definition of $H'$ for $X'$ is the same
as the definition of $H$ for $X$.
Let $\varepsilon\in(0,1/2)$.
There exists $c_2>0$ such that for large $t$, for every $x\in[\varepsilon, 1-\varepsilon]$,
\begin{align}
    f^{-}_{\gamma}(x)-o(n_t^{-1})
\leq
    \tilde f_{\gamma}(x)
\leq
    f_{\gamma}(x) \leq   f^{+}_{\gamma}(x)+o(n_t^{-1}),
\label{tft}
\end{align}
with $f_{\gamma}^{\pm}(x):=  \P\left(\frac{1}{R_1} \leq \frac{\gamma}{1-x}(1\pm \varepsilon_t'), \mathcal{H}_1 > t(1-x)(1 \mp \varepsilon_t') \right)$ and $\varepsilon_t'= e^{-c_2 h_t}$.
\end{lemma}

\medskip

\begin{lemma}   \label{lem5.4}
For any $ 0  <a <1/4 $, we have for any $t>0$,
\begin{align}
    \sum_{1\leq k \leq n_t}\P\left[ \bar{\mathcal{H}}_{k} >   1-a/2,\ 1-2a < \bar{\mathcal{H}}_{k-1}\leq 1-3a/4 \right]
\leq
    s(a,t),  \label{rec1}
\end{align}
with $s(a,t)$ such that $\lim_{t \rightarrow + \infty} s(a,t)=\textrm{const} \times a^{1-\kappa}$.
For any $\varepsilon\in(0,1/2)$,
\begin{align}
    \forall t>0,
\qquad
    \P\left[ \varepsilon t \leq  H(m_{N_t}) \leq (1-\varepsilon) t  \right]
\geq
    1- \tilde s(\varepsilon,t) \label{NNeps},
\end{align}
with $\tilde s(\varepsilon,t)$ such that $\lim_{t \rightarrow + \infty} \tilde s(\varepsilon,t)=\textrm{const} \times \varepsilon^{(1- \kappa) \wedge \kappa}$.
\end{lemma}

\medskip

We postpone the proof of these lemmata after the proof of Proposition \ref{ProSupL}.

\bigskip

\noindent {\bf Proof of Proposition \ref{ProSupL}:}
\noindent Recall from \eqref{eqDefNt} that $N_t$ is the largest index such that $\sup_{s \leq t} X(s) \geq m_{N_t}$.
In particular, $H\big(\tilde L_j\big)\leq H(\tilde m_{j+1})=H(m_{j+1}) \leq H(m_{N_t})\leq t$ for every $1\leq j\leq N_t-1$
on $\mV_t\cap\{N_t\leq n_t\}$.
The main idea is to use the fact that the supremum of the local time at time $t$ is achieved in the neighborhood of the $h_t$-minima $m_i$,
$1\leq i\leq n_t$. \\

\noindent \underline{\textbf{We start with the upper bound.}}
Let $\alpha>0$ and $0<\varepsilon<1/2$.
Notice that $\P(N_t=0, \mV_t)\leq \P[H(\tilde m_1)>t]\leq C_2 v_t$ by \eqref{TpsNeg}.
Using \eqref{TpsNeg}, \eqref{HLLM}, \eqref{Bas1}, \eqref{NNeps}, Lemma \ref{CVs} and Remark \ref{RemEgaliteAvecouSansTilde},
we have for $t$ large enough,
\begin{align}
&
    \P\left(
        \sup_{x \in \R} \loX(t,x) \leq \alpha t
      \right)
\leq
    E \left[\P^{\wk}\left( \max_{ 1\leq j \leq N_t} \loX(t,m_{j}) \leq \alpha t \right)\right]
    \label{proba1}
\\
\leq &
    E \left[
        \P^{\wk} \left(
                    \max_{ 1\leq j \leq N_t-1} \loX\big[H\big(\tilde L_j),\tilde m_{j}\big]
                    \leq \alpha t,
                    \loX(t, \tilde m_{N_t})\leq \alpha t,
                    \mathcal{Q}, \mathcal{B}_1(n_t),\mathcal{B}_{2}(n_t), { \mathcal V}_t
                 \right)
      \right]
\nonumber\\
 &+  \bar s(\varepsilon,t). \nonumber
\end{align}
with
$
    \mathcal{Q}
:=
    \{\varepsilon t \leq  H(m_{N_t}) \leq (1-\varepsilon) t ,\ 1\leq N_t \leq n_t\}
$
and $\bar s$ satisfying  $\lim_{t \rightarrow + \infty  } \bar{s}(\varepsilon,t) \leq C_+ \varepsilon^{(1-\kappa) \wedge \kappa}$.
We will introduce in what follows different measures denoted by the letter $\nu$; they depend on $k$
but we do not write $k$ as a subscript to simplify the notation.
First, define two measures $\nu_1^{\wk}$ and $\nu_2^{\wk}$ on $(0,1)$
by, for every $ 0 < y  < 1$,
\begin{align*}
    \nu_1^{\wk}(y)
:=  &
    \nu_1^{\wk}( [0,y])
\\
 := &
    \Pw\Big(\max_{ 1\leq j \leq k-1} \loX\big[H\big(\Lt_j\big),\mt_j\big] \leq \alpha t,
    H(\tilde m_{k})-\sum_{i=1}^{k-1}U_i< \tilde v_t,
    H(\mt_{k}) \leq y t \Big),
\\
    \nu_{2}^{\wk}(y)
 := &
    \nu_2^{\wk}( [0,y])
\\
 := &
    \Pw_{\mt_{k}}\Big(  \lo_{X'}[t(1-y),\mt_{k}] \leq \alpha t,
    \ H'(\mt_{k+1})>t(1-y),
\\
&  \qquad\qquad\qquad
    H'\big(\mt_{k+1}\big)< H'\big(\Lt_{k}^-\big),\
    H'\big(\mt_{k+1}\big)-H'\big(\Lt_k\big) \leq \tilde v_t \Big),
\end{align*}
with $X'$ a diffusion starting from $\mt_k$ independent of $X$  (conditionally on $W_{\kappa}$),
and $H'$ has the same definition as $H$ (see \eqref{DefH}) but for $X'$.
Partitioning on the values of $N_t$, and $H(\tilde m_{k})$,
we obtain  by the strong markov property (applied at time $H(\tilde m_k)$ under $\Pw$),
that the probability $E\big[\P^{\wk}(.)\big]$ in the line below \eqref{proba1} is smaller than
\begin{align}
&
    \sum_{1\leq k \leq n_t} \int_{\varepsilon}^{1-\varepsilon} E \left(\nu_2^{\wk}(x) \dd\nu_1^{\wk}(x)\right)  =
 \sum_{1\leq k \leq n_t} E\left [  \int_{\varepsilon}^{1-\varepsilon} \nu_2^{\wk}(x) \dd\nu_1^{\wk}(x)  \right]. \label{lastup}
\end{align}
 The next step is to prove that the previous expectation can be approximated by a product of expectations. First notice that both $ y \rightarrow \nu_1^{\wk}(y)$ and $ y \rightarrow \nu_2^{\wk}(y)$ are  positive increasing. 
So integrating by parts
 \begin{align}
    \int_{\varepsilon}^{1-\varepsilon} \nu_2^{\wk}(x) \dd\nu_1^{\wk}(x)
 &  =
    \left[\nu_2^{\wk}(x) \nu_1^{\wk}(x)  \right]_{\varepsilon}^{1- \varepsilon}-\int _{\varepsilon}^{1- \varepsilon}\nu_1^{\wk}(x) \dd\nu_2^{\wk}(x)
 \nonumber \\
 & \leq
    \left[\nu_2^{\wk}(x) \nu_1^{\wk}(x) \right]_{\varepsilon}^{1- \varepsilon}
    - \int_{\varepsilon}^{1- \varepsilon} \widetilde{\nu}^{\wk}_1(x) \dd\nu_2^{\wk}(x)
 \nonumber \\
 & =
    \left[\nu_2^{\wk}(x) \Big(\nu_1^{\wk}(x)-\widetilde{\nu}_1^{\wk}(x)\Big) \right]_{\varepsilon}^{1- \varepsilon} + \mathcal{I}_1,
 \label{term1}
 \end{align}
with
$
    \widetilde{\nu}^{\wk}_1(x)
:=
    \Pw \left(\mathcal{G}_1, H(\tilde m_{k})-\sum_{i=1}^{k-1}U_i< \tilde v_t,
        \sum_{i=1}^{k-1}U_i+ \tilde v_t  \leq x t\right)
\leq
    \nu_1^{\wk}(x)
$
and
$\mathcal{G}_1:=\{\max_{ 1\leq j \leq k-1} $
$ \loX(H(\Lt_j),\mt_j) \leq \alpha t \}$ and
 \begin{align*}
    \mathcal{I}_1
 &:=
    \int_{\varepsilon}^{1- \varepsilon} \nu_2^{\wk}(x)\dd \tilde{\nu}_1(x)
 \leq
    \int_{\varepsilon}^{1- \varepsilon} \nu_2^{\wk}(x)\dd \nu_3^{\wk}(x)
 =:
    \mathcal{I}_1',
 \\
    \nu_3^{\wk}(x)
 & :=
    \Pw \left(\mathcal{G}_1,\sum_{i=1}^{k-1}U_i+ \tilde v_t  \leq xt  \right).
 \end{align*}
{\bf First,} we deal with what is going to be a negligible part, that is to say the first term in \eqref{term1}.
As $\nu_1^{\wk}(x) \leq \Pw \left(\mathcal{G}_1, H(\tilde m_{k})-\sum_{i=1}^{k-1}U_i< \tilde v_t,\sum_{i=1}^{k-1}U_i \leq x t\right)$
because by definition $\sum_{i=1}^{k-1}U_i<H(\tilde m_{k})$, we have, for $\varepsilon<x<1-\varepsilon$,
\[
    \Big|\nu_1^{\wk}(x)-\widetilde{\nu}_1^{\wk}(x)\Big|
\leq
    \Pw \left( xt-\widetilde v_t<\sum_{i=1}^{k-1}U_i \leq xt \right)
=:
    h_k(x).
\]
so
$
    \big[\nu_2^{\wk}(x) \big(\nu_1^{\wk}(x)-\widetilde{\nu}_1^{\wk}(x)\big)\big]_{\varepsilon}^{1-\varepsilon}
\leq
    \nu_2^{\wk}(1-\varepsilon)h_k(1- \varepsilon)+\nu_2^{\wk}(\varepsilon) h_k(\varepsilon)
$.
Notice that
$\sum_{i=1}^{k-1}U_i $ is measurable with respect to
$\sigma\big( X(s), 0\leq s \leq H\big(\tilde L_{k-1}\big); W_{\kappa}(x), x \leq \tilde L_{k-1}^+\big)$,
since $\tilde L_{k-1}\leq \tilde L_{k-1}^+$,
whereas the event in the definition of $\nu_2^{\wk}$  belongs to
$$
    \sigma\big( X'(s), 0 \leq s \leq \min \big(H'\big(\tilde L_{k}^-\big), H'\big(\tilde m_{k+1})\big);
        W_{\kappa}(x)-W_{\kappa}(\mt_k),\tilde L_{k-1}^+ \leq x \leq \tilde L_{k+1}^+
        \big),
$$
with $X'$ an independent copy of $X$ starting at $\mt_k$. \\
So  independence of $X$ and $X'$, and independence of the two portions of the environment involved (see Lemma \ref{CVs}) imply independence between $\nu_2^{\wk}$ and $h_k$. Hence,
\begin{eqnarray}
&&
    E\Big(\Big[\nu_2^{\wk}(x) \Big(\nu_1^{\wk}(x)-\widetilde{\nu}_1^{\wk}(x)\Big)\Big]_{\varepsilon}^{1-\varepsilon}\Big)
\nonumber\\
& \leq &
    E\Big[ \nu_2^{\wk}(1-\varepsilon)\Big] E\big[h_k(1- \varepsilon)\big]
    +E\Big[\nu_2^{\wk}(\varepsilon)\Big] E\big[ h_k(\varepsilon)\big]
\nonumber \\
& = &
    E\Big[\widetilde \nu_2^{\wk}(1-\varepsilon)\Big] E\big[h_k(1- \varepsilon)\big]
    +E\Big[\widetilde \nu_2^{\wk}(\varepsilon)\Big] E\big[ h_k(\varepsilon)\big]
\label{2term}.
\end{eqnarray}
with for any $x$, \\
\begin{eqnarray*}
    \widetilde \nu_2^{\wk}(x)
& := &
    \Pw_{\mt_1} \Big(
            \lo_{X}[t(1-x),\mt_{1}] \leq \alpha t,
            \ H(\mt_{2})>t(1-x),
\\
&&
            \qquad\quad\qquad\qquad   H(\mt_{2})< H\big(\Lt_{1}^-\big),
            \ H\big(\mt_{2}\big)-H\big(\Lt_1\big) \leq \tilde v_t \Big).
\end{eqnarray*}
As $E\big(\widetilde \nu_2^{\wk}(x)\big) \leq \P [U_1>t(1-x)-\tilde v_t ] $
and for every small $\varepsilon>0$ and t large enough $h_k(x) \leq \Pw \left( (x-\varepsilon)t<\sum_{i=1}^{k-1}U_i \leq xt  \right) $ we can apply Proposition \ref{Pro3.3}, we get
\begin{eqnarray*}
&&
     E\big[h_k(1- \varepsilon)\big] E\Big[\tilde \nu_2^{\wk}(1-\varepsilon)\Big]
\\
& \leq &
    \P \left( \frac{1-2\varepsilon}{1+\varepsilon_t}<\sum_{i=1}^{k-1}\frac{\mathcal{H}_i}{t}
    \leq \frac{1-\varepsilon}{1- \varepsilon_t} \right) \P\bigg(\mathcal{H}_1>\frac{t\varepsilon-\tilde v_t}{1+\varepsilon_t} \bigg)
+3e^{-D_1h_t}.
\end{eqnarray*}
By \eqref{cvmesure9.1} and the first part of Lemma \ref{5.4}, for any $0 <a<1$ and $b>0$,
\begin{eqnarray}
    \lim_{t \rightarrow +\infty} \sum_{1\leq k \leq n_t} \P \left( 1-a <\sum_{i=1}^{k-1}\frac{\mathcal{H}_i}{t} \leq 1  \right) \P(\mathcal{H}_1>t b )
& = &
    \frac{\textrm{const}}{b^{\kappa}}
     \int_{1-a}^{1} y^{\kappa-1}\dd y
\nonumber\\
& \leq &
    \frac{\textrm{const}}{b^{\kappa}}(1-(1-a)^{\kappa}).
    \hphantom{blab}
    \label{eqPourDeuxiemeProba}
\end{eqnarray}
Therefore, we obtain
\begin{align*}
    \sum_{1\leq k \leq n_t} E\Big[\tilde \nu_2^{\wk}(1-\varepsilon)\Big] E\big[h_k(1- \varepsilon)\big]
& \leq
    C_+ \cdot u(t,\varepsilon)
\end{align*}
with $u$ a positive function such that
$\lim_{t \rightarrow + \infty} u(t,\varepsilon) = \max(\varepsilon^{1- \kappa}, \varepsilon^{\kappa})$.
A similar argument also works for the second term in \eqref{2term}, which yields
\begin{align}
    \sum_{1\leq k\leq n_t}   E \left[\left[\nu_2^{\wk}(x) \Big(\nu_1^{\wk}(x)-\tilde{\nu}_1^{\wk}(x)\Big) \right]_{\varepsilon}^{1- \varepsilon} \right]
\leq
    2C_+ \cdot u(t,\varepsilon). \label{last0}
\end{align}

\noindent {\bf We now deal with $\mathcal{I}_1'$}. By independence between $X$ and $X'$, and the independent parts of the potential $\wk$ involved in $\nu_2^{\wk}(x)$ and $\nu_3^{\wk}(x)$,
\begin{align}
E(\mathcal{I}_1')=  \int_{\varepsilon}^{1- \varepsilon} \nu_2(x)\dd \nu_3(x),
\label{last11}
\end{align}
with $\nu_2(x):=E\big(\nu_2^{\wk}(x)\big)=E\big(\widetilde \nu_2^{\wk}(x)\big)$ and $\nu_3(x):=E\big(\nu_3^{\wk}(x)\big)$.

{ \noindent
By the lower bound in Lemma \ref{Lem6.2}, we have
$
    \nu_3(x)
=
    F_{\alpha}(x- \tilde v_t/t)
\geq
    F^-_{\alpha}(x- \tilde v_t/t) -e^{-D_1 h_t}
$
for every $x>\varepsilon$ for large $t$.
So,
again since $ y \rightarrow \nu_2(y)$ is positive increasing and $\nu_3$ is a repartition function,
integrating by parts twice as in \eqref{term1}   gives with the change of variables $u=x-\tilde v_t/t$,
\begin{eqnarray}
    \int_{\varepsilon}^{1- \varepsilon}  \nu_2(x)\dd \nu_3(x)
&\leq &
    \int_{\varepsilon-\tilde v_t/t}^{1- \varepsilon-\tilde v_t/t} \nu_2(x+\tilde v_t/t)\dd F_{\alpha}^-(x)+e^{-D_1 h_t}
\nonumber\\
&&
    +\Big[\Big(F_{\alpha}(x)-F_{\alpha}^-(x)\Big)\nu_2(x+\tilde v_t/t)\Big]^{1- \varepsilon-\tilde v_t/t}_{\varepsilon-\tilde v_t/t}.
\label{5.65b}\hphantom{azeert}
\end{eqnarray}
Recall (see before Lemma \ref{lemtps}) that $\tilde v_t/t=2 / \log (h_t)=o(1)$ as $t\to+\infty$.
Then we can prove in a similar way we have obtained \eqref{last0} that:
\begin{align}
    \sum_{1\leq k \leq n_t}
    \Big[\Big(F_{\alpha}(x)-F_{\alpha}^-(x)\Big)\nu_2(x+\tilde v_t/t)\Big]^{1- \varepsilon-\tilde v_t/t}_{\varepsilon-\tilde v_t/t}
& \leq
    C_+ \cdot u(t,\varepsilon),
\label{eqQuasiSimilaire}
\end{align}
with as usual a possibly enlarged $C_+$.
Indeed by Lemma \ref{Lem6.2},
$
    -\big(F_{\alpha}(\varepsilon-\tilde v_t/t)-F_{\alpha}^-(\varepsilon-\tilde v_t/t)\big)\nu_2(\varepsilon)
\leq
    e^{-D_1 t}
=o(n_t^{-1})
$ for every $1\leq k\leq n_t$ for large $t$,
and
$
    \big(F_{\alpha}-F_{\alpha}^-\big)(1-\varepsilon-\tilde v_t/t)
\leq
    \big(F_{\alpha}^+-F_{\alpha}^-\big)(1-\varepsilon-\tilde v_t/t)
    +e^{-D_1 t}
\leq
    \P\big(\max_{ 1\leq j \leq k-1} \ell_j \in[\gamma_t^{-} t, \gamma_t^{+} t]\big)
    + \P\big(\sum_{i=1}^{k-1}\mathcal{H}_i \in[x_t^{-} t, x_t^{+} t ] \big)
    +e^{-D_1 t}
$
for every $k\leq n_t$ for large $t$,
with $\gamma=\alpha$ and $x=1-\varepsilon-\tilde v_t/t$.
The first probability is less than
$
    n_t\P(S_1 {\bf e}_1\in[\gamma_t^{-} t, \gamma_t^{+} t])
=
    n_t\mathbb E\big(\int_{\gamma t (1-2\varepsilon_t)/S_1}^{\gamma t (1+2\varepsilon_t)/S_1}e^{-u/2}\dd u/2\big)
\leq
    8 n_t\varepsilon_t \sup_{v\geq 0} (ve^{-v})
=
    o(1/n_t)
$,
whereas the second one is treated as \eqref{eqPourDeuxiemeProba}, which leads to \eqref{eqQuasiSimilaire}.

So the important term in
the right hand side of inequality  \eqref{5.65b} comes from the integral.
We now work on $\nu_2(x)$. We have,
\begin{align*}
&
    \nu_2(x)
\\
& \leq
    E\big( \Pw_{\mt_1} \big[  \lo_{X}(t(1-x),\mt_{1}) \leq \alpha t,
            \ H\big(\Lt_1\big)>t(1-x)-\tilde v_t,
            \ H\big(\Lt_1\big)< H\big(\Lt_{1}^-\big)\big]\big)
\\
& \leq
    E\big( \Pw_{\mt_1} \big[  \lo_{X}(t(1-x)-\tilde v_t,
    \ \mt_{1}) \leq \alpha t,
\\
&
    \qquad\qquad\qquad\qquad
    \ H\big(\Lt_1\big)>t(1-x)-\tilde v_t,
    \ H\big(\Lt_1\big)< H\big(\Lt_{1}^-\big)\big]\big)
    =f_\alpha(x+\tilde v_t/t),
\end{align*}
as defined in Lemma \ref{lem6.3b}.
Then, as $F_{\alpha}^-(x)$ is positive and increasing in $x$, using Lemma \ref{lem6.3b} with $\gamma= \alpha$,
we obtain
\begin{align}
     \int_{\varepsilon-\tilde v_t/t}^{1- \varepsilon-\tilde v_t/t} \nu_2(x+\tilde v_t/t)\dd F_{\alpha}^-(x)
\leq
    \int_{\varepsilon-\tilde v_t/t}^{1- \varepsilon-\tilde v_t/t} f_{\alpha}^+(x+2\tilde v_t/t)\dd F_{\alpha}^-(x)    +o(n_t^{-1}). \label{Fpfp}
\end{align}
Now, as $ f_{\alpha}^{+}(x+2\tilde v_t/t)$ can be written (since $\mathcal{H}_k =\ell_k R_k$, see Proposition \ref{Pro3.3}),
\[
    f_{\alpha}^{+}(x+2\tilde v_t/t)
=
    \P\left((1-x-2\tilde v_t/t)  \frac{\ell_k}{\mathcal{H}_k} \leq \alpha(1+\varepsilon_t'),
    \mathcal{H}_k > t(1-x-2\tilde v_t/t)(1 -\varepsilon_t')  \right),
\]
we get by independence  of the random variables $((\ell_j,\mathcal{H}_j), j \leq n_t)$,
}
\begin{eqnarray}
&&
    \int_{\varepsilon-\tilde v_t/t}^{1- \varepsilon-\tilde v_t/t} f_{\alpha}^+(x+2\tilde v_t/t) \dd F_{\alpha}^-(x)
\nonumber\\
& \leq &
    \P\left[
            \big(1- \bar{\mathcal{H}}_{k-1}\big)
            \frac{\bar{\ell}_{k}-\bar{\ell}_{k-1} } {\bar{\mathcal{H}}_{k}-\bar{\mathcal{H}}_{k-1}}
        \leq
            \alpha+\tilde \varepsilon_t(k),\
        \bar{\mathcal{H}}_{k} \geq 1-\delta_t',
      \right.
\nonumber\\
&&
\qquad\qquad\qquad\qquad\qquad
      \left.
        {\max_{ 1\leq j \leq k-1} \frac{\ell_j}{t} \leq
        \alpha,\
        \bar{\mathcal{H}}_{k-1}\leq 1-\varepsilon+ \delta_t' }
      \right],
\label{eqComparaisonIntegraleProba}
\end{eqnarray}
with $\delta_t':=3\tilde v_t/t$,
$\tilde \varepsilon_t(k):=\left(\alpha+\ell_k/\mathcal{H}_{k}\right)\delta_t'$.\\
\noindent The idea now is to make appear the event $\big\{\mathcal{N}_{t}^{2 \varepsilon}=k\big\}$ in the above probability  (recall the definition of
$\mathcal{N}_{t}^{2 \varepsilon}$ given in Proposition \ref{ProSupL}) and then sum over $k$.
\\
We first prove that the sum over $k \leq n_t$,  of the above probability  is small if we intersect its event with the event
$
    \big\{\mathcal{N}_{t}^{2\varepsilon} \neq k\big\}
$.
In other words, let us prove that
\begin{align}
    \text{$\textstyle \sum_1$}
:=
    \sum_{1\leq k\leq n_t}  \P\left[\bar{\mathcal{H}}_{k} \geq 1-\delta_t',\
        { \bar{\mathcal{H}}_{k-1}\leq 1-\varepsilon+\delta_t' },\ \mathcal{N}_{t}^{2\varepsilon}  \neq k \right] \label{sum1}
\end{align}
is small. As
$
    \big\{\mathcal{N}_{t}^{2\varepsilon}  \neq k\big\}
=
    \big\{ \bar{\mathcal{H}}_{k} \leq 1-2\varepsilon \big\}\cup \big\{ \bar{\mathcal{H}}_{k-1}  >  1-2\varepsilon \big\} $,
and since for $t$ large enough,
$
    \big\{\bar{\mathcal{H}}_{k} \geq 1-\delta_t' \big\} \cap \big\{\bar{\mathcal{H}}_{k} \leq 1-2\varepsilon \big\}= \emptyset
$,
we have
$$
    \text{$\textstyle \sum_1$}
\leq
    \sum_{1\leq k\leq n_t}  \P\left[\bar{\mathcal{H}}_{k} \geq 1-\delta_t',\
    1- 2\varepsilon < { \bar{\mathcal{H}}_{k-1}\leq 1-\varepsilon+\delta_t' } \right].
$$
Therefore, for $t$ large enough, with $s(\varepsilon,t)$ defined in Lemma \ref{lem5.4},
\begin{equation}\label{eqQuasiFinaleUpperBound}
        \text{$\textstyle \sum_1$}
\leq
    \sum_{1\leq k\leq n_t}  \P\left[\bar{\mathcal{H}}_{k} > 1-\varepsilon/2,\
    1- 2\varepsilon < { \bar{\mathcal{H}}_{k-1}\leq 1-3\varepsilon/4 } \right] \leq s(\varepsilon,t)
\end{equation}
by \eqref{rec1}.
Finally, combining equations from \eqref{last11} to \eqref{eqQuasiFinaleUpperBound} leads to
\begin{eqnarray}
&&
    \sum_{1\leq k \leq n_t} E\big(\mathcal{I}_1'\big)
\nonumber \\
& \leq &
    \P\left[
        (1- \bar{\mathcal{H}}_{ \mathcal{N}_{t}^{2\varepsilon}-1})
        \frac{\bar{\ell}_{ \mathcal{N}_{t}^{2\varepsilon}}-\bar{\ell}_{ \mathcal{N}_{t}^{2\varepsilon}-1} }
            {\bar{\mathcal{H}}_{ \mathcal{N}_{t}^{2\varepsilon}}-\bar{\mathcal{H}}_{ \mathcal{N}_{t}^{2\varepsilon}-1}}
            \leq \alpha+\tilde \varepsilon_t(\mathcal{N}_{t}^{2\varepsilon}),
             \ {\max_{ 1\leq j \leq  \mathcal{N}_{t}^{2\varepsilon}-1} \frac{\ell_j}{t} \leq
             \alpha
             }
             \right]
\nonumber\\
&&
    +s(\varepsilon,t)+C_+u(t,\varepsilon)+o(1) \label{last1}.
\end{eqnarray}
To finish we have to deal with $\tilde \varepsilon_t( \mathcal{N}_{t}^{2\varepsilon} )$, a basic computation partitioning on the values of  $\mathcal{N}_{t}^{2\varepsilon}$, shows that
$
    \P\big[\tilde \varepsilon_t (\mathcal{N}_{t}^{2\varepsilon})) \geq \alpha \sqrt{\delta_t'}/6\big]
\leq
    C_+ \P\left( {R_{1}} \leq \sqrt \delta_t' \right)
=
    o(1)
$
as \bl{$R_1$ converges in distribution to $\mathcal{R}_{\kappa}$ which is almost surely positive.}
Collecting this last fact, \eqref{proba1}, \eqref{lastup}, \eqref{term1}, \eqref{last0}  and \eqref{last1} finish the proof of the upper bound.

\medskip

\noindent \underline{\textbf{Proof of the lower bound:}} \\
The proof here follows the same line as the upper bound.
The main difference comes from the fact that we can no longer use the inequality
$\sup_{x \in \R} \loX(t,x) \geq \sup_{ 1\leq j \leq N_t}\loX(t,m_j)$.
So for this part of the proof we stress on what is different from the upper bound,
 and refer to the previous computations when very few changes occur.  \\
Assume for the moment that
\begin{equation}
    \P\Big( \Big\{\sup_{x \in \R} \loX(t,x)  \geq 2 \tilde w_t \Big\}=: \mathcal{E}_2 \Big)
\geq
    1-o(1),
\label{EqProuveeDebutPreuveTheorem}
\end{equation}
with
\bl{$
    \tilde w_t
:=
    t  e^{(\kappa(1+ 3\delta)-1) \phi(t)}
$},
and recall that $\delta$ is chosen small enough such that $\kappa(1+ 3\delta)<1$ (see Lemma \ref{LemmaProbaMaxLocHorsdeTOUTESVallees}).
This fact \eqref{EqProuveeDebutPreuveTheorem} is a direct consequence of the upper-bound of $\P( \sup_{x \in \R} \loX(t,x) \leq \alpha t)$
(see at the beginning of the proof of Theorem \ref{ADV} page \pageref{PageReferencePourPreuveEq518} for a
 proof of \eqref{EqProuveeDebutPreuveTheorem}).
Recall \eqref{Dj}, and define for any $\ell \geq 1$,
\begin{align*}
&
    \mathcal{E}_3(\ell)  := \mathcal{E}^1_3(\ell) \cap \mathcal{E}^2_3(\ell), \quad \textrm{with}
\\
&
        \mathcal{E}^1_3(\ell)
    :=
        \bigcap_{j=1}^{\ell-1}
        \left\{ \sup_{x \in \Dt_j } \Big[\loX\big(H\big(\tL\big),x\big)-\loX\big(H\big(\tm\big),x\big)\Big]\leq  t  \tilde \alpha_t    \right\},
\\
&
        \mathcal{E}^2_3(\ell)
    :=
        \left\{ \sup_{x \in \Dt_{\ell} }
            \Big[ \loX(t,x)- \loX\big(H\big(\tilde m_{\ell}\big),x\big)\Big] \leq t \tilde \alpha_t
            \right \} ,
\end{align*}
with $\tilde \alpha_t:= (\alpha  t-2\tilde w_t)/t  $.
Recall the definitions of the events $\mathcal{B}_i$, $1\leq i \leq 4$
in Sections \ref{sectionEventB} and \ref{sectionEventC}.
We have for large $t$,
\begin{eqnarray*}
&&
    \big\{ \sup\nolimits_{x \in \R_+} \loX(t,x) \leq \alpha t \big\}
\cap
    \mathcal{V}_t \cap \mathcal{E}_2 \cap\{N_t\leq n_t \} \cap \cap_{i=1}^4\mathcal{B}_i(n_t)
\\
& \supset &
    \mathcal{E}_3(N_t) \cap  \mathcal{V}_t \cap \mathcal{E}_2 \cap\{ N_t\leq n_t\} \cap \cap_{i=1}^4\mathcal{B}_i(n_t).
\end{eqnarray*}
Indeed, $\loX(t,x)\leq \tilde w_t$ for every $x\in \big(\R_+-\cup_{j=1}^{n_t}\big[\tilde L_j^-, \tilde L_j\big]\big)$
on $\mathcal{B}_2(n_t)\cap \mathcal{B}_3(n_t)\cap \mathcal{V}_t\cap\{N_t\leq n_t\}$,
and on the same event intersected with $\mathcal{B}_4(n_t)$,
$\loX(t,x)\leq \tilde w_t+t e^{-2\phi(t)}<2\tilde w_t$
for every $x\in \cup_{j=1}^{n_t}\big(\big[\tilde L_j^-, \tilde L_j\big]\cap \overline{\Dt_{j}}\big)$,
whereas for $x\in \Dt_{j}$,
$\loX(H(\tilde m_j),x)\leq \tilde w_t$ if $j\leq n_t$ and
$\loX(t,x)-\loX\big(H\big(\tilde L_j\big),x\big)\leq \tilde w_t$ if $j<N_t$.
Notice that by Lemmata \ref{CVs}, \ref{lemtps}, \ref{negloc1}, \ref{negloc2} and the above assumption \eqref{EqProuveeDebutPreuveTheorem},
$$
    \P\big(\mathcal{V}_t \cap \mathcal{E}_2\cap\{N_t \leq n_t\}\cap\cap_{i=1}^4\mathcal{B}_i(n_t)\big)
\geq
    1-o(1).
$$
We now deal with
$\P(\mathcal{E}_3(N_t) \cap \mathcal{B}_1(N_t) \cap \mathcal{B}_2(n_t) \cap \mathcal{V}_t\cap\{N_t\leq n_t\} )$.
Using Lemma \ref{CVs}, the fact that
$ H(\Lt_{k}) \leq H(\mt_{k+1}) $ and the strong Markov property with respect to $\P^{\wk}$, we obtain
\begin{align*}
&
    \P(\mathcal{E}_3(N_t) \cap \mathcal{B}_1(N_t) \cap \mathcal{B}_2(n_t) \cap \mathcal{V}_t\cap \mathcal{Q})
\\
& \geq
    \sum_{k=1}^{n_t} E\left( \int_{\varepsilon}^{1-\varepsilon}
     \nu_4^{\wk}(y) \Pw \left( \mathcal{E}_3^1(k), \mathcal{B}_1(k), {\mathcal{B}}_2(k-1),H(\mt_{k})/t \in \dd y\right)   \right)-o(1)
\end{align*}
with
\begin{align*}
&    \nu_4^{\wk}(y)
\\
:=&
    \Pw_{\mt_{k}} \left( \sup_{x \in \Dt_{k} }  \lo_{X}(t(1-y),x) \leq t \tilde \alpha_t,\ H(\Lt_{k}) > t(1-y),\ H(\Lt_k)<H(\Lt_k^-) \right),
\end{align*}
Now, by computations similar  to the ones giving the upper bounds in \eqref{last0} and \eqref{last11}, we have
\begin{align*}
&
    \P(\mathcal{E}_3(N_t) \cap \mathcal{B}_1(N_t) \cap \mathcal{B}_2(n_t) \cap \mathcal{V}_t,\mathcal{Q} ) \\
& \geq
    \sum_{k=1}^{n_t} \int_{\varepsilon}^{1-\varepsilon} E\left( \nu_4^{\wk}(y) \dd\nu_5^{\wk}(y) \right) -o(1)
    = \sum_{k=1}^{n_t} \int_{\varepsilon}^{1-\varepsilon}  \nu_4(y) \dd\nu_5(y)-o(1).
\end{align*}
with $\nu_5^{\wk}(y):=\Pw \left( \mathcal{E}_3^1(k), \mathcal{B}_1(k), {\mathcal{B}}_2(k-1),\sum_{i=1}^{k-1}U_i/t \leq y\right) $,
$\nu_4(y):=E\big(\nu_4^{\wk}(y)\big)$ and $\nu_5(y):=E\big(\nu_5^{\wk}(y)\big)$.
The next step is to remove $\mathcal{B}_1(k)$ in the above expression. For that, we only have to prove that
\begin{align*}
    \sum_{k=1}^{n_t} \int_{\varepsilon}^{1-\varepsilon}
        E\left( \nu_4^{\wk}(y)  \Pw \left( \mathcal{E}_3^1(k), \bar{\mathcal{B}}_1(k), {\mathcal{B}}_2(k-1),\sum_{i=1}^{k-1}U_i/t \in \dd y\right)\right)
\end{align*}
is negligible, one can check that this quantity is smaller than
\begin{eqnarray*}
&&
    \sum_{k=1}^{n_t} \int_{\varepsilon}^{1-\varepsilon}
                    E\left[\Pw_{\mt_{k}} (H(\Lt_k)<H(\Lt_k^-),H(\Lt_{k}) > t(1-y))\right]
\\
&&
    \qquad\qquad\qquad\qquad\qquad
    \P \left(\bar{\mathcal{B}}_1(k),\mathcal{B}_2(k-1),\sum_{i=1}^{k-1}U_i/t \in \dd y\right)
\\
 & \leq &
    \sum_{k=1}^{n_t} \P\left(\sum_{i=1}^{k-1}U_i/t \leq 1, \sum_{i=1}^{k}U_i/t>1,\ \bar{\mathcal{B}}_1(k)\right)
\\
& \leq &
    \P\left(\bar{\mathcal{B}}_1( n_t)\right)
\leq
    C_2 v_t =o(1),
\end{eqnarray*}
where the last inequality comes from \eqref{TpsNeg}.
Therefore, collecting the above computations yields
\begin{align*}
    \P\left( \sup_{x \in \R} \loX(t,x) \leq \alpha \right )
& \geq
    \sum_{k=1}^{n_t} \int_{\varepsilon}^{1- \varepsilon} \nu_4(y) \dd\tilde \nu_5(y)   -o(1),
\end{align*}
with  $\tilde \nu_5(y):=e^{-\kappa \phi(t)} \sum_{k \leq n_t}   \P \left(\mathcal{E}^1_3(k),{\mathcal{B}}_2(k-1), \sum_{i=1}^{k-1}U_i/t    \leq y \right)$.
\\
\noindent We start with an estimation of the repartition function $\tilde \nu_5(y)$. Recall that like in the proof of Lemma \ref{lemX}, by the strong Markov property, the occupation time formula \eqref{1.2} and \eqref{1.3} the sequence $(U_j,\{\loX(H(\tL),x)-\loX(H(\tm),x),x \in \mathcal{D}_j \}, j \leq n_t)$ under $\mathcal{B}_2(n_t)$ is equal to a sequence $(H_j(\tL),\{\lo_j(H_j(\tL),x),x \in \mathcal{D}_j \}, j \leq n_t)$, with this time
\begin{eqnarray*}
    H_j(\tL)
& := &
    A^j(\tL) \int_{\tL^-}^{\tL} \bl{e^{-\tV(u)}}\mathcal{L}_{B^j}[\tau^{B^j}(1),A^j(u)/A^j(\tL)]\dd u,
\\
    \lo_j(H_j(\tL),x)
& :=  &
    A^j(\tL) e^{-\tV(x)}\mathcal{L}_{B^j}[\tau^{B^j}(1),A^j(x)/A^j(\tL)], \
\end{eqnarray*}
where $A^j(u)= \int_{\tm}^{u} e^{\tV(x)}\dd x$.
Using Remark \ref{RemEgaliteAvecouSansTilde}, Lemma \ref{CVs},  Fact \ref{Fact_Williams} {\bf(ii)},
and then \eqref{bessel4} and \eqref{bessel5}, we have for large $t$ for any $1\leq j\leq n_t$ since $\phi(t)=o(\log t)$,
\begin{equation}\label{D1tau}
\begin{split}
    P\big[\tilde \tau_j(\kappa r_t/8)\leq \tilde m_j+ r_t\leq \tilde \tau_j(r_t)\big]
\geq
    1-C_+ e^{ -(c_-) r_t},
\\
    P\big[\tilde \tau_j^-(r_t)\leq \tilde m_j- r_t\leq \tilde \tau_j^-(\kappa r_t/8)\big]
\geq
    1-C_+ e^{ -(c_-) r_t}.
\end{split}
\end{equation}
with $c_->0$. Therefore for any $j$,
$
    P\big(\mathcal{D}_j \subset \big[\tilde \tau^-_j(r_t),\tilde \tau_j(r_t) \big]\big) \geq 1-2C_+ e^{ -(c_-) r_t}
$.
Then on $\{\mathcal{D}_j \subset [\tilde \tau^-_j(r_t),\tilde \tau_j(r_t) ] \}$,
for any $x \in \mathcal{D}_j $,
$$
    \lo_j(H_j(\tL),x)
\leq
    A^j(\tL) \mathcal{L}_{B^j}[\tau^{B^j}(1),A^j(x)/A^j(\tL)].
$$
Also with probability $\geq 1-2C_+ e^{ -(c_-) r_t}$,
$\mathcal{D}_j \subset [\tilde \tau^-_j(r_t),\tilde \tau_j(r_t) ] $ so for any $x \in \mathcal{D}_j$,
\begin{align}
    A^j(\tilde \tau^-_j(r_t))
\leq
    A^j(x)
\leq
    A^j(\tilde \tau_j(r_t)).
\label{TrucpourV1b}
\end{align}
With Remark \ref{RemEgaliteAvecouSansTilde}, Lemma \ref{CVs}, Fact \ref{Fact_Williams}
and \eqref{4.6bb},
we obtain with a probability  larger than $1-e^{- (c_-) r_t }$,
\begin{eqnarray}
    -e^{-h_t/4}
& \leq &
    -e^{2 r_t} e^{-(1-1/2)h_t}
\leq
    \frac{A^j(\tilde \tau_j^-(r_t))}{A^j\big(\tL\big)}
\leq
    \frac{A^j(\tilde \tau_j(r_t))}{A^j\big(\tL\big)}
\nonumber\\
& \leq &
    e^{2 r_t } e^{-(1-1/2)h_t}
\leq
    e^{-h_t/4}.
\label{TrucpourV1}
\end{eqnarray}
Therefore, applying \eqref{Dev} (with $\delta=e^{-h_t/4}$ and $\varepsilon=\delta^{1/3}$),
we obtain with a probability larger than $1-e^{- (c_-) r_t}$,
\begin{align}
\sup_{x \in \mathcal{D}_j } {A^j\big(\tL\big)} \mathcal{L}_{B^j}\big(\tau^{B^j}(1),A^j(x)/A^j\big(\tL\big)\big)
\leq
    A^j\big(\tL\big)  \mathcal{L}_{B^j}\big(\tau^{B^j}(1),0\big)  \big(1+ e^{-h_t/12}\big) . \label{TrucpourV2}
\end{align}
Collecting the different estimates we then obtain,
\begin{align*}
    \tilde \nu_5(y)
\geq
    \P \left(\max_{ 1\leq j \leq k-1}  \lo_j\big(H_j\big(\tL\big),\tm\big) \leq  t  \bar \alpha_t
    , \ \sum_{j=1}^{k-1} \frac{H_j\big(\tL\big)}{t} {\leq y} \right) -C_+ e^{- (c_-) r_t},
 \end{align*}
with $\bar \alpha_t:=\tilde \alpha_t \big(1+ e^{-h_t/12}\big)^{-1}$.
We can then inverse the equality in law we have used above, and then obtain
\begin{align*}
    \tilde \nu_5(y)
\geq
    F_{\bar \alpha_t}(y)-C_+ e^{- (c_-) r_t},
 \end{align*}
with $F_{\bar \alpha_t}$ defined in Lemma \ref{Lem6.2}.
Then we can follow the same lines as for the upper bound
(especially  computations after \eqref{last0}), and obtain via Lemma \ref{Lem6.2}
and by choosing $C_0$ large enough in such a way that $(c_-) r_t/ \phi(t)= (c_-) C_0  > \kappa (1+ \delta)$:
 \begin{align*}
    \int_{\varepsilon}^{1-\varepsilon} \nu_4(y) \dd\tilde \nu_5(y)
 \geq
    \int_{\varepsilon}^{1-\varepsilon} \nu_4(y) \dd F_{\bar \alpha_t}^+(y)-o(n_t^{-1}).
\end{align*}

\noindent Remark also {that \eqref{TrucpourV2} implies} the  concentration of the local time at the $h_t$-minima:
with probability larger than $1-C_+ e^{- (c_-) r_t}$,
\begin{align}
    {\left| \sup_{y \in \mathcal{D}_j} \loX_j \big(H_j \big(\tL\big),y\big)
    - \lo_j \big(H_j\big(\tL\big),\tm\big) \right| \leq  e^{-h_t/12} \lo_j\big(H_j\big(\tL\big),\tm\big)}.
\label{maxaufond}
\end{align}

\noindent We now work on $\nu_4(y)$. By the second part of Lemma \ref{CVs}  it is equal to
\begin{align*}
&
    E\left(\Pw_{\mt_{1}} \left(\sup_{\bl{z \in \Dt_{1}} }  \lo_{X'}\bl{(t(1-y),z)} \leq t \tilde \alpha_t,
    H'\big(\Lt_{1}\big)>\bl{t(1-y)}, H'\big({\Lt_1}\big)<H'\big({\Lt^-_1}\big)  \right) \right)
\\
 =: &
    \tilde \nu_4(y),
\end{align*}
and by Lemma \ref{lem6.3b}, $\tilde \nu_4(y) \geq f_{\tilde \alpha_t}^-(y) -o(n_t^{-1})$. Therefore
\begin{align*}
    \int_{\varepsilon}^{1-\varepsilon} \nu_4(y) \dd\tilde \nu_5(y)
\geq
    \int_{\varepsilon}^{1-\varepsilon}f_{\tilde \alpha_t}^-(y)\dd F_{\bar \alpha_t}^+(y) -o(n_t^{-1}).
\end{align*}
From now on, the computations are very close from that of the upper bound (see \eqref{Fpfp} and  below) and we do not give more details. \hfill $\square$

\vspace{0.5cm}
\noindent {\bf Proof of Lemmata \ref{Lem6.2}, \ref{lem6.3b} and \ref{lem5.4}.}

\medskip

%
\noindent \\
\noindent {\bf Proof of Lemma \ref{Lem6.2}:}
This is a direct consequence {of Proposition \ref{Pro3.3}.}
\hfill $\square$

\bigskip

\noindent {\bf Proof of Lemma \ref{lem6.3b}:}
To obtain the result we use a similar method than in \cite{AndDiel}.
That is to say, we study the inverse of the local time at $\mt_1$, and use our knowledge about $H(\Lt_1)$. \bl{From the definitions of $f_{\gamma}$ and $\tilde f_{\gamma}$ we have easily $\tilde f_{\gamma}(x) \leq  f_{\gamma}(x)$ for all $x$.
So, to prove \eqref{tft}, we only need to prove the upper bound for ${f_{\gamma}}$ and the lower bound for ${\tilde f_{\gamma}}$.}
We fix $\varepsilon\in(0,1/2)$.


\noindent \textit{$\bullet$ Upper bound for $f_{\gamma}(x)$}.
Recall that $\sigma(u,\mt_1) = \inf \{s>0,\, \loX(s,\mt_1) \geq  u \}$, $u\geq 0$.
First, notice that for $0<x<1$, $f_\gamma(x)$ is equal to
\bl{\begin{align}
&
    E\left[  \P_{\mt_1}^{W_{\kappa}}\Big(  \loX(t(1-x),\mt_1) \leq \gamma t, H\big(\Lt_1\big)>t(1-x), H\big(\Lt_1\big)< H\big(\Lt_1^-\big) \Big)\right]
\nonumber \\
&  =
    E\left[  \P_{\mt_1}^{W_{\kappa}}\Big(\sigma(\gamma t,\mt_1) \geq t(1-x) , H\big(\Lt_1\big)>t(1-x), H\big(\Lt_1\big)< H\big(\Lt_1^-\big) \Big)\right]
\label{losi} \\
&  =
    E\left[  \P_{\mt_1}^{W_{\kappa}}\Big(H\big(\Lt_1\big) > \sigma(\gamma t,\mt_1) \geq t(1-x), H\big(\Lt_1\big)< H\big(\Lt_1^-\big) \Big)\right]
\label{losi1part} \\
&
    \qquad\qquad
    + E\left[  \P_{\mt_1}^{W_{\kappa}}\Big(\sigma(\gamma t,\mt_1) > H\big(\Lt_1\big)>t(1-x), H\big(\Lt_1\big)< H\big(\Lt_1^-\big) \Big)\right].
\label{losi2part}
\end{align}

Let us first study the expectation in \eqref{losi1part}.
On $\big\{H\big(\Lt_1\big) > \sigma(\gamma t,\mt_1), H\big(\Lt_1\big)< H\big(\Lt_1^-\big)\big\}$  under $\P_{\mt_1}^{W_{\kappa}}$,
$ X$ remains between $\Lt_1^-$ and $\Lt_1$ until time $\sigma(\gamma t,\mt_1)$ which is finite.
On this event and under $\P_{\mt_1}^{W_{\kappa}}$,
considering \eqref{1.2} and \eqref{1.3} as in (\cite{Shi}  p. 248),
the inverse of the local time can be written for $X$ starting at $\mt_1$ as }
\begin{equation} \label{eqSigmaX}
    \siX(\gamma t,\mt_1)
=
    \int_{\Lt_1^-}^{\Lt_1}e^{-\tilde V^{(1)}(z)}\lo_{B}\big(\sigma_{B}(\gamma t,0), A^1(z)\big)\dd z
=:
    I,
\end{equation}
where
$A^1(z)=\int_{\mt_1}^{z}e^{\tilde V^{(1)}(y)} \dd y$ and $B$ is a standard Brownian motion
independent of $W_{\kappa}$, such that $B$ starts at $A^1(\mt_1)=0$ and is killed when it first hits $A^1(\Lt_1)$.
In \eqref{eqSigmaX}, we integrate only between $\Lt_1^-$ and $\Lt_1$ because under $\P_{\mt_1}^{W_{\kappa}}$,
$e^{-\tilde V^{(1)}(z)}\lo_{B}\big(\sigma_{B}(\gamma t,0), A^1(z)\big)=\loX(\siX(\gamma t,\mt_1), z)=0$
for $z\notin \big[\Lt_1^-,\Lt_1\big]$ as explained after \eqref{losi2part}.
 We have
\[
    I
=
    {\gamma t}  \int_{\Lt_1^-}^{\Lt_1}e^{-\tilde V^{(1)}(z)}\lo_{\tilde B}\big(\sigma_{\tilde B}(1,0),\tilde a(z)\big)\dd z,
\]
with $\tilde a(z):=(\gamma t)^{-1}A^1(z)=(\gamma t)^{-1} \int_{\mt_1}^{z}e^{\tilde V^{(1)}(y)} \dd y$
and where $\tilde B := B((\gamma t)^2 .)/(\gamma t)$.
By scale invariance $\tilde B$ is also a standard Brownian motion that we still denote by $B$ in the sequel.
Also, recall that $\sigma_U(r,y):= \inf\{s>0, \ \lo_U(s,y) >r\}$ for $r>0$, $y\in \R$ is the inverse of the local time of the process $U$.
{Since we consider $X$ starting at $\mt_1$, we have
$H(\tilde L_1) = H(\tilde L_1)-H(\tilde m_1)=U_1$, for which Proposition \ref{Pro3.3} gives}
\begin{align}
    \E\Big(\P_{\mt_1}^{W_{\kappa}}\Big\{\big|H(\tilde L_1)-\mathcal{H}_1\big|\leq \varepsilon_t \mathcal{H}_1\Big\}\Big)
=
    \P\left(\left\{ |
            U_1
            -
            \mathcal{H}_1| \leq \varepsilon_t
            \mathcal{H}_1
            \right\}
        =:
            \mathcal{G}_1
            \right)
\geq
    1-e^{-D_1 h_t},
\label{R11}
\end{align}
with $\varepsilon_t:=e^{-d_1 h_t}$, if $\delta>0$ is chosen small enough. This will explain the appearance of $\mathcal{H}_1$
in $f_\gamma^{\pm}(x)$.
So, we now deal with $I$. Notice that $(\gamma t)^{-1} I $ can be split into two terms $(\gamma t)^{-1}I=I_1+I_2$, with
$$
    I_1
:=
    \int_{\tilde \tau_1^-(h_t/2)}^{\tilde \tau_1(h_t/2)}e^{-\tilde V^{(1)}(z)}\lo_{B}\big(\sigma_{B}(1,0),\tilde a(z)\big)\dd z,
$$
and $I_2:=(\gamma t)^{-1}I-I_1\geq 0$.
We now prove that the main contribution in $(\gamma t)^{-1}I $
comes from $I_1$ and obtain its approximation in probability.
Let $\varepsilon\in(0,1/100)$.
First, using the second part of Lemma \ref{CVs}, followed by
Remark \ref{RemEgaliteAvecouSansTilde},
Fact \ref{Fact_Williams} {\bf(ii)} (for which we need $i\geq 2$),
\eqref{4.6bb}
and finally the first part of Lemma \ref{CVs}, we get
\begin{eqnarray}
&&
    P\big[\big|A^1\big(\tilde \tau_1^-(h_t/2)\big)\big|\leq e^{h_t(1+ \varepsilon)/2},
          \big|A^1\big(\tilde \tau_1(h_t/2)\big)\big|\leq e^{h_t(1+ \varepsilon)/2 }
     \big]
\nonumber\\
& = &
    P\big[\big|A^2\big(\tilde \tau_2^-(h_t/2)\big)\big|\leq e^{h_t(1+ \varepsilon)/2},
          \big|A^2\big(\tilde \tau_2(h_t/2)\big)\big|\leq e^{h_t(1+ \varepsilon)/2}
     \big]
\nonumber\\
& \geq &
    1-2P\big[F^+(h_t/2)>e^{h_t(1+ \varepsilon)/2}\big]-P\big[\overline{{\mV}_t}\big]
\geq
    1-  C_+  e^{-\kappa \varepsilon h_t/4}.
\label{InegPourAhtsur2}
\end{eqnarray}
Therefore,
since
$
    \tilde a\big(\tilde \tau_1^-(h_t/2)\big)
\leq
    \tilde a(z)
\leq
    \tilde a\big(\tilde \tau_1(h_t/2)\big)
$
for all $z \in \big[ \tilde \tau_1^-(h_t/2),\tilde \tau_1(h_t/2)\big]$,
\begin{equation}\label{eqPetitesValeursatilde}
    P\big(\forall {z \in \big[ \tilde \tau_1^-(h_t/2),\tilde \tau_1(h_t/2)\big]},\   |\tilde a(z)| \leq e^{-(\log t) (1-3\varepsilon)/2}\big)
\geq
    1-  C_+  e^{-\kappa \varepsilon h_t/4 }.
\end{equation}
Also, using \eqref{5.96} and the second Ray-Knight theorem (see before \eqref{5.96}), we have
\begin{align}
    \P \left(\sup_{ |u| \leq  e^{-(\log t) (1-3\varepsilon)/2} }
    \Big|\lo_{B}(\sigma_{B}(1,0),u)-1\Big| \geq \widehat \varepsilon_t \right)
\leq
     e^{-t^{\varepsilon}/16}
\label{CVL_B} .
\end{align}
with
$\widehat \varepsilon_t:=t^{-(1-5\varepsilon)/4}$.
So we obtain
\begin{align}\label{InegI1}
    \E\Big[\P_{\mt_1}^{W_{\kappa}}\big(\big|I_1-  \tilde R_1 \big| \leq
    \widehat \varepsilon_t \tilde R_1  \big)\Big]
\geq
    1- C_+e^{-\kappa \varepsilon h_t/4 },
\end{align}
with $\tilde R_1:= \int_{\tilde \tau_1^-(h_t/2)}^{\tilde \tau_1(h_t/2)}e^{-\tilde V^{(1)}(z)}\dd z$.
We now prove that $I_2$ is negligible compared to the integral $\tilde R_1$ which appears in the previous equation,
and then  compared to $I_1$.
First thanks to \eqref{5.97} and the second Ray-Knight theorem, we have
$$
    \E\bigg[\P_{\mt_1}^{W_{\kappa}}\bigg(\sup_{z \in [{\Lt_1^-},{\Lt_1}] } \lo_B\big[\sigma_B(1,0),\tilde a(z)\big]>
    e^{\varepsilon \log t}\bigg)\bigg]
\leq
    2e^{-\varepsilon \log t}.
$$
So with probability larger than $1- 2e^{-\varepsilon \log t}$, we have
$$
    I_2 \leq e^{\varepsilon \log t} \left(\int_{\Lt_1^-}^{\tilde \tau_1^-(h_t/2)}e^{-\tilde V^{(1)}(z)}\dd z
    + \int_{\tilde \tau_1(h_t/2) }^{ \Lt_1} e^{-\tilde V^{(1)}(z)}\dd z \right)
=:
        e^{\varepsilon \log t} I_3 .
$$
By Lemma \ref{lemI3}, with a probability larger than $1-2e^{- (c_-) \varepsilon h_t}$ for large $t$,
\begin{align*}
I_3 \leq  C_+ h_t^2 e^{-(1- \varepsilon)h_t/2}.
\end{align*}
Also, by Lemma \ref{lemX},  with probability larger $1-e^{-(D_-)h_t}$, $\tilde R_1=R_1$
(which is the same $R_1$ as in \eqref{R11}), which law is given by the sum of two independent copies of $F^-(h_t/2)$.  So using \eqref{Fmoins}, with a probability larger than $ 1-2e^{-(D_-) h_t}$,
\begin{align*}
\tilde R_1 = R_1 \geq e^{-\varepsilon h_t/2}.
\end{align*}
We deduce from the last three inequalities that with a probability larger than $1- e^{-(c_-) \varepsilon h_t},$
\begin{equation}\label{InegI2}
    I_2
<
    R_1 e^{-(1-5 \varepsilon) h_t/2}
=
    \tilde R_1 e^{-(1-5 \varepsilon) h_t/2}.
\end{equation}
\bl{Finally, using $(\gamma t)^{-1}I=I_1+I_2$ together with \eqref{InegI1} and \eqref{InegI2}, we get
\begin{eqnarray}
&&
    E\left[\Pw_{\mt_1}\left(
                        \big|I - \gamma t R_1 \big|
                        \geq
                        2 t^{-(1-5\varepsilon)/4} (\gamma t) R_1,
                        \right.
     \right.
\nonumber\\
&&
    \left.
                        \left.
                        \qquad\qquad\quad
                        H\big(\Lt_1\big) > \sigma(\gamma t,\mt_1),\
                        H\big(\Lt_1\big)< H\big(\Lt_1^-\big)
                        \right)
    \right]
\nonumber \\
& \leq &
    C_+e^{-\varepsilon (c_-)  h_t }.
\label{invloct}
 \end{eqnarray}
We recall that by \eqref{eqSigmaX}, $\sigma(\gamma t,\mt_1) = I$ on
$\big\{H\big(\Lt_1\big) > \sigma(\gamma t,\mt_1), H\big(\Lt_1\big)< H\big(\Lt_1^-\big)\big\}$  under $\P_{\mt_1}^{W_{\kappa}}$.
Hence, combining \eqref{invloct} with \eqref{R11} gives for large $t$ for every $x\in[\varepsilon,1-\varepsilon]$,
\begin{align}
&
    \left \{H\big(\Lt_1\big) > \sigma(\gamma t,\mt_1) \geq t(1-x),
            H\big(\Lt_1\big)< H\big(\Lt_1^-\big) \right \}
\nonumber  \\
&
        \subset \left \{ \frac{1}{R_1} \leq \frac{\gamma}{1-x}(1 + \varepsilon_t'), \mathcal{H}_1 > t(1-x)(1 - \varepsilon_t'),
        H\big(\Lt_1\big) > \sigma(\gamma t,\mt_1)
        \right \}
        \cup \mathcal{E}_{\varepsilon}^1 \label{losi1partevt},
 \end{align}
where $\mathcal{E}_{\varepsilon}^1$ is such that
$E\big[\Pw_{\mt_1}(\mathcal{E}_{\varepsilon}^1 )\big]
\leq C_+e^{-(\varepsilon c_-)  h_t } + e^{-D_1 h_t}$ and where, as defined in the statement of the lemma,
$\varepsilon_t' = e^{-c_2 h_t}$ with $c_2>0$ chosen small enough.

Now, let us study \eqref{losi2part}. On the event inside the probability in \eqref{losi2part},
$\sigma(\gamma t,\mt_1)$ might be infinite.
We work under $\Pw_{\mt_1}$. There exists a Brownian motion $B$ such that, with $T^1$ playing under $\Pw_{\mt_1}$
the same role as $T$ does under $\P$ (see \eqref{1.2}),
$H\big(\tilde L_1\big) = T^1 \big(\tau^B\big(A^1\big(\tilde L_1\big)\big)\big)$
and  $\sigma(\gamma t,\mt_1)=T^1(\sigma_B(\gamma t, 0))$ (as in \eqref{eqSigmaX} and in \cite{Shi} p. 248). Also by \eqref{1.2},
notice for further use that under $\Pw_{\mt_1}$,
\begin{equation}\label{EqTempsLocalEnInverse}
    \mathcal{L}\big(\sigma(y t,\tilde m_1), z\big)
=
    e^{-\tilde V^{(1)} (z)} \mathcal{L}_B\big(\sigma_B(yt,0), A^1(z)\big),
\qquad z\in\R, y\in(0,+\infty).
\end{equation}
So, we have
\begin{eqnarray*}
    \sigma(\gamma t,\mt_1) > H\big(\Lt_1\big)
& \Leftrightarrow &
    \sigma_B(\gamma t, 0) > \tau^B\big(A^1\big(\tilde L_1\big)\big)
\\
& \Leftrightarrow &
    \mathcal{L}_B \big[\sigma_B(\gamma t, 0), 0 \big] = \gamma t
    > \mathcal{L}_B \big[\tau^B\big(A^1\big(\tilde L_1\big)\big), 0 \big].
\end{eqnarray*}
Now, note that, as in \eqref{eqDefej} in the proof of Lemma \ref{lemX},
$
    \mathcal{L}_B \big[\tau^B(A^1(\tilde L_1)), 0 \big]
=
     A^1(\tilde L_1) \mathcal{L}_{\tilde B} \big( \tau^{\tilde B}(1), 0 \big)
$,
where $\tilde B := B((A^1(\tilde L_1))^2 .)/A^1(\tilde L_1)$.
Also, by definition of ${\bf e}_1$ given in \eqref{eqDefej},
we have
$
    \mathcal{L}_{\tilde B} \big( \tau^{\tilde B}(1), 0 \big)
=
    {\bf e}_1
$.
As a consequence,
\[
    \sigma(\gamma t,\mt_1) > H\big(\Lt_1\big)
\Leftrightarrow
    \gamma t > A^1\big(\tilde L_1\big) {\bf e}_1
\Leftrightarrow
    \gamma t R_1 > A^1\big(\tilde L_1\big) {\bf e}_1 R_1.
\]
Then, according to \eqref{3.30}, we have $A^1\big(\tilde L_1\big) \geq \big(1-e^{-(d_-)h_t}\big) {S}_1$ with probability greater than $1- e^{-(D_-)h_t}$.
Moreover, according to $\eqref{R11}$ and to the fact that under $\Pw_{\mt_1}$ the diffusion $X$ starts at $\tilde m_1$,
we have $\mathcal{H}_1={\bf e}_1 S_1R_1 \geq (1 + \varepsilon_t)^{-1} H (\tilde L_1)$ with probability greater than $1- e^{-(D_-)h_t}$. As a consequence,
\begin{eqnarray}
    \sigma(\gamma t,\mt_1) > H(\Lt_1)
\Rightarrow
    \gamma t R_1 > \big(1-e^{-(d_-)h_t}\big) (1 + \varepsilon_t)^{-1} H(\Lt_1),
\label{casinfini}
\end{eqnarray}
except on an event which probability $E\big[\Pw_{\mt_1}(.)\big]$ is less than $2e^{-(D_-)h_t}$.
Combining this with \eqref{R11} we get
for large $t$ for every $x\in[\varepsilon,1-\varepsilon]$,
\begin{align}
&
    \left \{ \sigma(\gamma t,\mt_1) > H(\Lt_1)>t(1-x), H(\Lt_1)< H(\Lt_1^-) ) \right \} \nonumber \\
&
    \subset \left \{ \frac{1}{R_1} \leq \frac{\gamma}{1-x}(1+ \varepsilon_t'), \mathcal{H}_1 > t(1-x)(1 - \varepsilon_t'),
    \sigma(\gamma t,\mt_1) > H(\Lt_1)\right \}
    \cup
    \mathcal{E}_{\varepsilon}^2 \label{losi2partevt},
 \end{align}
where $\mathcal{E}_{\varepsilon}^2$ is such that
$
    E\big[\Pw_{\mt_1}(\mathcal{E}_{\varepsilon}^2 )]
\leq
    2e^{-(D_-)h_t} + e^{-D_1 h_t}
$
and where, as before, $\varepsilon_t' = e^{-c_2 h_t}$ with $c_2>0$  possibly smaller than before.

Combining \eqref{losi1part}, \eqref{losi2part} \eqref{losi1partevt} and \eqref{losi2partevt}
with the strong Markov property, we get
for large $t$ for every $x\in[\varepsilon,1-\varepsilon]$,
since $\phi(t)=o(\log t)$,
\begin{eqnarray*}
    f_{\gamma}(x)
& \leq &
    \P\left( \frac{1}{R_1} \leq \frac{\gamma}{1-x}(1+ \varepsilon_t'), \mathcal{H}_1 > t(1-x)(1 - \varepsilon_t') \right)
    + o(n_t^{-1}).
\\
& = &
    f_{\gamma}^+(x)
    + o(n_t^{-1}).
\end{eqnarray*}
}

\noindent \textit{$\bullet$ Lower bound for  $\tilde f_{\gamma}$.}
\bl{Let $\tilde \gamma := \gamma \big(1 + e^{-h_t/12}\big)^{-1}$ and
$y :=(1-x)/[R_1(1- 4\widehat\varepsilon_t)] $.
We have to distinguish the cases $H\big(\Lt_1\big) > \sigma(\tilde \gamma t,\mt_1)$
and $\sigma(\tilde \gamma t,\mt_1) > H\big(\Lt_1\big) $.
We work under $\Pw_{\mt_1}$.
On $\big\{y \leq \tilde \gamma, H\big(\Lt_1\big) > \sigma(\tilde \gamma t,\mt_1) \geq t(1-x), H\big(\Lt_1\big)< H\big(\Lt_1^-\big)\big\}$,
we can express the local time of $X$ at the inverse of its local time in $\tilde m_1$ at time $y t$ in terms of the standard Brownian motion
driving the diffusion}.
More precisely by \eqref{EqTempsLocalEnInverse} and by scale invariance,
there exists a Brownian motion $B$ such that for any $z\in \Dt_1$,
\begin{equation}
    \loX\big(\siX(y t,\mt_1),z\big)
=
    ( y t)e^{-\tilde V^{(1)}(z)} \lo_{B}\big(\sigma_{B}(1,0),\widehat a(z)\big)
\label{eqlaw}
\end{equation}
with $\widehat a(z):=(y t)^{-1} \int_{\mt_1}^{z}e^{\tilde V^{(1)}(u)} \dd u=A^1(z)/(yt)$.
Notice that by \eqref{eqDefF+G+},
$
    F^-(h_t/2)
\leq
    \tau^{\BP}(h_t/2)
$
in law, so
$
    P\big[R_1>8h_t/\kappa \big]
\leq
    2P\big[F^-(h_t/2)>4h_t/\kappa\big]
\leq
    2P\big[\tau^{\BP}(h_t/2)>4h_t/\kappa\big]
\leq
    e^{-(c_-)h_t}
$
for large $t$.
Moreover, we prove with the same method used to prove \eqref{D1tau}
that
$
    \tilde\tau^-(h_t/2)
\leq
    \tilde m_1-r_t
\leq
    \tilde m_1+r_t
\leq
    \tilde\tau(h_t/2)
$
with probability at least $1-C_+ e^{-(c_-)h_t}$. This and \eqref{InegPourAhtsur2} give
$
    -e^{h_t(1+\varepsilon)/2}
\leq
    A^1[\tilde\tau^-(h_t/2)]
\leq
    A^1(z)
\leq
    A^1[\tilde\tau(h_t/2)]
\leq
    e^{h_t(1+\varepsilon)/2}
$
for any $z \in \mathcal{D}_1$ with probability $\geq 1-e^{-(c_-)\varepsilon h_t}$.
So, for large $t$ for every $x\in[\varepsilon,1-\varepsilon]$,
$
    |\widehat a(z)|
\leq
        e^{h_t(1+\varepsilon)/2}R_1/(t(1-x))
\leq
    e^{-(\log t) (1-3\varepsilon)/2}
$
for these $z$ with such probability.
Hence with the same method we used to prove
\eqref{InegI1} from \eqref{eqPetitesValeursatilde} and \eqref{CVL_B}, we get
for large $t$ for every $x\in[\varepsilon,1-\varepsilon]$,
\begin{align*}
    E\left(\Pw_{\mt_1} \left(\sup_{ z \in \Dt_1}
            \Big| \lo_{B}(\sigma_{B}(1,0),\widehat a(z))-1\Big| \leq \widehat \varepsilon_t \right) \right)
\geq
    1-2e^{- (c_-) \varepsilon h_t}.
\end{align*}
The above inequality together with \eqref{eqlaw} imply
that for large $t$ for every $x\in[\varepsilon,1-\varepsilon]$,
\begin{eqnarray}
&&
    E\left(\Pw_{\mt_1} \left( \left\{ \exists z \in \Dt_1, \ \left| \loX(\siX(y t,\mt_1),z)- y t e^{-\tilde V^{(1)}(z)} \right| \geq 2 yt e^{-\tilde V^{(1)}(z)} \widehat \varepsilon_t, \right. \right. \right. \nonumber \\
&&
    \qquad\qquad\qquad  \left. \left. \left. y \leq \tilde \gamma,
    \ H\big(\Lt_1\big) > \sigma(\tilde \gamma t,\mt_1) \geq t(1-x),
    \ H\big(\Lt_1\big)< H\big(\Lt_1^-\big) \right\}   \right) \right)
\nonumber\\
& \leq  &
    2e^{- (c_-) \varepsilon h_t}. \label{tpsinv}
\end{eqnarray}
On $\big\{y \leq \tilde \gamma, H\big(\Lt_1\big) > \sigma(\tilde \gamma t,\mt_1) \geq t(1-x),
H\big(\Lt_1\big)< H\big(\Lt_1^-\big)\big\}$, if $ t(1-x) > \sigma(y t , \mt_1)$,
then $\sigma(y t , \mt_1)-y t R_1<-4t y R_1\widehat\varepsilon_t$,
and by \eqref{invloct} (applied with $\gamma$ replaced by $y$),
this has on the previous event a probability
$E\big(\Pw_{\mt_1}(.)\big)$ less than $C_+e^{-(c_-)\varepsilon h_t}$.
Thus on the previous event,
we have $ t(1-x) \leq \sigma(y t , \mt_1)$, except on a sub event of probability smaller than
$C_+e^{-(c_-) \varepsilon h_t}$.
This is true for every $x\in[\varepsilon,1-\varepsilon]$ for large $t$.

\bl{
Then since the local time is increasing in time,
we have on the previous event
for any $z \in \Dt_1$,
$
    \loX(t(1-x),z)
\leq
    \loX(\sigma(y t , \mt_1),z)
$,
which is  according to \eqref{tpsinv} less than
$
    y t e^{-V^{(1)}(z)} (1+2\widehat \varepsilon_t)
\leq
    y t  (1+2\widehat \varepsilon_t )
$
for every $z \in \Dt_1$ with probability $E\big(\Pw_{\mt_1}(.)\big)$ at least $1-2e^{- (c_-) \varepsilon h_t}$.
Combining this and the definition of our $y$ gives for large $t$,
for every $x\in[\varepsilon,1-\varepsilon]$,
\begin{eqnarray}
&&
    E\left(\Pw_{\mt_1} \left( \left\{\frac{\sup_{z \in \Dt_1}  \loX(t(1-x),z)}{t} >
            \frac{(1-x)}{R_1}   \frac{1+2\widehat \varepsilon_t }{1-4\widehat \varepsilon_t} \right\}=: \overline{\mathcal{G}_2} \right. \right.
\nonumber \\
&&
    \qquad\qquad\qquad
    \left. \left. \cap \left \{ y \leq \tilde \gamma, H\big(\Lt_1\big)
    > \sigma(\tilde \gamma t,\mt_1) \geq t(1-x), H\big(\Lt_1\big)< H\big(\Lt_1^-\big) \right \} \right) \right)
\nonumber\\
& \leq &
    (2+C_+) e^{-(c_-)\varepsilon  h_t}.
\label{L11}
\end{eqnarray}
As a consequence, for $t$ large enough so that
$1 + 2\widehat \varepsilon_t \leq 1 + e^{-h_t/12}$, we have
for every $x\in[\varepsilon,1-\varepsilon]$,
\begin{eqnarray}
&&
    \left \{ y \leq \tilde \gamma, H\big(\Lt_1\big) > \sigma(\tilde \gamma t,\mt_1) \geq t(1-x), H\big(\Lt_1\big)< H\big(\Lt_1^-\big) \right \}
\nonumber\\
& \subset &
    \left \{ \sup_{z \in \Dt_1}  \loX(t(1-x),z) \leq y(1+2\widehat \varepsilon_t)t\leq   \gamma t,
      \ H\big(\Lt_1\big) > \sigma(\tilde \gamma t,\mt_1) \right \} \cup \mathcal{E}_{\varepsilon}^3
      \phantom{blabla}
\label{mino5.3casfini}
\end{eqnarray}
by definition of $\tilde \gamma$,
where $\mathcal{E}_{\varepsilon}^3$ is such that $E\big(\Pw_{\mt_1}(\mathcal{E}_{\varepsilon}^3 )\big) \leq (2+C_+) e^{- (c_-)\varepsilon  h_t}$.

On the other hand,
from the definition of $\sigma(.,\mt_1)$, \eqref{maxaufond} and the definition of $\tilde \gamma$,
we have for large $t$ for every $x\in[\varepsilon,1-\varepsilon]$,
\begin{align}
&
    \left \{ y \leq \tilde \gamma, \sigma(\tilde \gamma t,\mt_1) > H\big(\Lt_1\big)>t(1-x),
            H\big(\Lt_1\big)< H\big(\Lt_1^-\big)  \right \}
\nonumber \\
\subset &
    \left \{  \loX\big(H\big(\Lt_1\big), \mt_1\big) \leq \tilde \gamma t, \sigma(\tilde \gamma t,\mt_1) > H\big(\Lt_1\big)>t(1-x),
                    H\big(\Lt_1\big)< H\big(\Lt_1^-\big)  \right \}
\nonumber \\
\subset &
    \left \{ \sup_{z \in \Dt_1}  \loX\big(H\big(\Lt_1\big),z\big) \leq \gamma t,
            \sigma(\tilde \gamma t,\mt_1) > H\big(\Lt_1\big)>t(1-x)
            \right \}
            \cup \mathcal{E}_{\varepsilon}^4
\nonumber \\
\subset &
    \left \{ \sup_{z \in \Dt_1}  \loX(t(1-x),z) \leq \gamma t, \ \sigma(\tilde \gamma t,\mt_1) > H\big(\Lt_1\big) \right \}
    \cup \mathcal{E}_{\varepsilon}^4
\label{losi2partevtbis},
\end{align}
where $\mathcal{E}_{\varepsilon}^4$ is the event where \eqref{maxaufond} fails,
it is such that $E\big(\Pw_{\mt_1}(\mathcal{E}_{\varepsilon}^4)\big) \leq C_+ e^{- (c_-) r_t}$.

Combining \eqref{mino5.3casfini} and \eqref{losi2partevtbis} we get
for large $t$ for every $x\in[\varepsilon,1-\varepsilon]$, under $\Pw_{\mt_1}$,
\[
    \left \{ y \leq \tilde \gamma, H(\Lt_1)>t(1-x), H(\Lt_1)< H(\Lt_1^-)  \right \} \subset \left \{ \sup_{z \in \Dt_1}  \loX(t(1-x),z) \leq \gamma t \right \} \cup \mathcal{E}_{\varepsilon}^5,
\]
where $\mathcal{E}_{\varepsilon}^5$ is such that
$
    E\big(\Pw_{\mt_1}(\mathcal{E}_{\varepsilon}^5 )\big)
\leq
    C_+ e^{- (c_-) r_t}
=
    C_+ e^{- (c_-) C_0\phi(t)}
=
    o(n_t^{-1})
$ as $t\to+\infty$ is we choose $C_0$ large enough.
Combining this with \eqref{R11}, \eqref{HLLM} and Proposition \ref{Pro3.3}, we obtain
for large $t$ for every $x\in[\varepsilon,1-\varepsilon]$,
\begin{eqnarray*}
    \tilde f_{\gamma}(x)
& \geq &
    \P \left(\frac{(1-x)}{R_1} \leq  \gamma (1-\varepsilon_t'),{\bf e}_1 S_1R_1>t  (1-x) (1+\varepsilon_t')  \right) -o(n_t^{-1})
\\
& = &
    f_{\gamma}^-(x)-o(n_t^{-1}),
\end{eqnarray*}
where the constant $c_2$ in the definition of $\varepsilon_t'=e^{-c_2 h_t}$ has been decreased if necessary.
This proves the lower bound for $\tilde f_{\gamma}(x)$ and then finishes the proof of the lemma.\hfill  $\square$
}

\bigskip

\noindent {\bf Proof of Lemma \ref{lem5.4}:}
Let $0<a<1/4$.
We start with \eqref{rec1}. By  Proposition \ref{Pro3.3}, the $\mathcal{H}_i$, $i\geq 1$ are i.i.d.,
so  $\bar{\mathcal{H}}_{k-1}$ and $\bar{\mathcal{H}}_{k}-\bar{\mathcal{H}}_{k-1}=\mathcal{H}_k$ are independent for $k\geq 1$. Thus for $t>0$,
\begin{eqnarray}
&&
    \sum_{1\leq k \leq n_t}\P\left[ \bar{\mathcal{H}}_{k} > 1-a/2,\ 1-2a < \bar{\mathcal{H}}_{k-1}\leq 1-3a/4 \right]
\nonumber\\
&= &
    \hspace{-1mm}
    \int_{1-2a}^{1-3a/4}d\mu_t(x)e^{\kappa \phi(t)}
    \P\left[ \mathcal{H}_{1} > 1-x-a/2 \right],
\end{eqnarray}
where the measure $\mu_t$ is defined by
$
    \int_0^xd\mu_t(y)
:=
    e^{- \kappa \phi(t)}\sum_{1\leq k \leq n_t}\P\left[  \bar{\mathcal{H}}_{k-1}\leq x \right]
$.
We know that $\mu_t$ converges vaguely as $t\to+\infty$ to the measure $\mu $ which has a density with respect to the Lebesgue measure
equal to $( \Gamma(\kappa) \CK)^{-1} x^{\kappa-1} \un_{x>0}$, with $\CK>0$ {(see  Lemma \ref{5.4})}.
Also thanks to Lemma \ref{lemproba}, $e^{\kappa \phi(t)}P\left[{\mathcal{H}}_{1}/t > x \right]$
converges uniformly on every compact subset of $(0,+ \infty)$ to $\CK x^{-\kappa}/ \Gamma(1- \kappa)$. Therefore,
\begin{align*}
&
    \lim_{t \rightarrow + \infty} \sum_{1\leq k \leq n_t}
    \P\left[ \bar{\mathcal{H}}_{k} > 1-a/2,\ 1-2a < \bar{\mathcal{H}}_{k-1}\leq 1-3a/4 \right] \\
&=
    \frac{1}{\Gamma(\kappa)\Gamma(1- \kappa)}\int_{1-2a}^{1-3a/4} x^{\kappa-1} (1-x- a/2)^{-\kappa}\dd x
\\
&\leq
        \textrm{const} \times a^{1-\kappa}.
\end{align*}
{For \eqref{NNeps}, we apply \eqref{6.5.2} with $r= \varepsilon\in(0,1/2)$ and $s=1- \varepsilon$, which gives
\begin{align*}
 &
    \lim_{t \rightarrow + \infty}  \P\left( \varepsilon t \leq H(m_{N_t}) \leq (1-\varepsilon) t  \right)
 \\
 & =
    1-\frac{\sin(\pi \kappa )}{ \pi}
   \left( \int_{0}^{\varepsilon} x^{ \kappa-1} (1-x)^{-\kappa} \textnormal{d}x+ \int_{1- \varepsilon}^{1} x^{ \kappa-1} (1-x)^{-\kappa} \textnormal{d}x \right)
 \\
 & \geq
    1-\frac{\sin(\pi \kappa )}{ \pi} \left( \frac{(1-\varepsilon)^{- \kappa}}{\kappa} \varepsilon^{\kappa}+\frac{(1- \varepsilon)^{\kappa-1}}{1-\kappa} \varepsilon^{1- \kappa} \right),
\end{align*}  which implies the result.  }
\hfill
$\square$
\\

\noindent {\bf Proof of Theorem \ref{ADV}:}
The proof of this theorem is a direct consequence of Propositions \ref{ProSupL} and \ref{CVY1Y2} and
of Lemmata \ref{continuityofJ} and  \ref{lemContJI}.
Notice that the proof of the upper bound does not use the proof of the lower bound,
but we use the upper bound for the proof of the lower bound.
In particular from the upper bound of Theorem \ref{ADV} (which makes use of the upper bound of Proposition \ref{ProSupL} but not of its lower bound),
we have
$
    \limsup_{t \rightarrow + \infty} \P(\lo^*(t) < 2 \tilde w_t)
\leq
    \P\big(\mathcal{Y}_1^{\natural}\big(\mathcal{Y}_2^{-1}(1)^-\big) \leq \varepsilon\big)
$
for any $\varepsilon>0$ as $\lim_{t \rightarrow + \infty}\tilde w_t/t  =0$.
From this, as $\mathcal{Y}_1^{\natural}\big(\mathcal{Y}_2^{-1}(1)^-\big)$
is positive, we obtain  $\lim_{t \rightarrow + \infty} \P(\lo^*(t) < 2 \tilde w_t)=0$,
which proves assertion \eqref{EqProuveeDebutPreuveTheorem}   at the beginning of the proof of the lower bound of Proposition \ref{ProSupL}.

\label{PageReferencePourPreuveEq518}

\bl{Thanks to Proposition \ref{ProSupL} and to the remark before this proposition,
we only need to study the convergence of $\mathcal{P}_1^{\pm}$ (the limit when $t$ goes to infinity and then the limit when $\varepsilon$ goes to $0$). The latter can be written in term of functionals of $(Y_1, Y_2)^t$ as follows.
Let
$\mathbb{Y}_t:=(Y_2^t)^{-1}(1-2\varepsilon)$; we have $\mathcal{N}_t^{2\varepsilon}e^{-\kappa\phi(t)}=\mathbb{Y}_t$, and
\begin{eqnarray*}
\mathcal{P}_1^{\pm}
& = &
    P\left[\big(1- Y_2^t (\mathbb{Y}^-_t)\big)\frac{Y_1^t (\mathbb{Y}_t) - Y_1^t (\mathbb{Y}_t^-)}{Y_2^t (\mathbb{Y}_t) - Y_2^t (\mathbb{Y}^-_t)}
    \leq \alpha_t^{\pm},
    \ (Y_1^t)^{\natural} (\mathbb{Y}_t^-) \leq \alpha_t^{\pm} \right]
\\
& = &
    P\Bigg[\big(1- \tilde K_{I, 1-2\varepsilon}^-((Y_1, Y_2)^t)\big)
    \frac{K_{I, 1-2\varepsilon}((Y_1, Y_2)^t) - K_{I, 1-2\varepsilon}^-((Y_1, Y_2)^t)} {\tilde K_{I, 1-2\varepsilon}((Y_1, Y_2)^t) - \tilde K_{I, 1-2\varepsilon}^-((Y_1, Y_2)^t)} \leq \alpha_t^{\pm},
\nonumber\\
&&
    \qquad\qquad\quad\qquad\qquad\qquad\qquad\qquad\qquad\qquad\quad\ \
    J_{I, 1-2\varepsilon}^-((Y_1, Y_2)^t) \leq \alpha_t^{\pm}
    \Bigg],
\end{eqnarray*}
with the notation $K_{I,a}$,   $\tilde K_{I,a},\dots$ introduced in \eqref{DefK} and before.
The hypotheses of Lemma 4.5 are: finite number of large jumps on compact intervals, strictly increasing, starting at $0$, and jumping over $1$ without reaching it.
These properties are naturally almost surely satisfied by a $\kappa$-stable subordinator so,
almost surely, the paths of $(\mathcal{Y}_1,\mathcal{Y}_2)$ satisfy the hypotheses of Lemma \ref{lemContJI}
(see e.g. \cite{Bertoin} III.2 p. 75).
Therefore they are points of continuity for $J_{I, 1-2\varepsilon}^-$, $K_{I, 1-2\varepsilon}^-$,
$K_{I, 1-2\varepsilon}$, $\tilde K_{I, 1-2\varepsilon}^-$ and $\tilde K_{I, 1-2\varepsilon}$.
Combining this continuity with Proposition \ref{CVY1Y2}, {continuous mapping theorem}, and replacing the functionals by their expressions, we obtain, when $t$ goes to infinity, the convergence of $\mathcal{P}_1^{\pm}$ to
\begin{align*}
&
    P\bigg[
        \big(1- \mathcal{Y}_2\big(\mathcal{Y}_2^{-1}(1-2\varepsilon)^-\big)\big)
        \frac{\mathcal{Y}_1\big(\mathcal{Y}_2^{-1}(1-2\varepsilon)\big) - \mathcal{Y}_1\big(\mathcal{Y}_2^{-1}(1-2\varepsilon)^-\big)} {\mathcal{Y}_2\big(\mathcal{Y}_2^{-1}(1-2\varepsilon)\big) - \mathcal{Y}_2\big(\mathcal{Y}_2^{-1}(1-2\varepsilon)^-\big)} \leq \alpha,
\\
&
    \qquad\qquad\qquad\qquad\qquad\qquad\qquad\qquad\qquad\qquad\qquad\qquad
        \mathcal{Y}_1^{\natural}\big(\mathcal{Y}_2^{-1}(1-2\varepsilon)^-\big) \leq \alpha
    \bigg].
\end{align*}
Then, note that almost surely $\mathcal{Y}_2\big(\mathcal{Y}_2^{-1}(1)^-\big) < 1$
so we have a.s. $\mathcal{Y}_2^{-1}(1-2\varepsilon) = \mathcal{Y}_2^{-1}(1)$ for all $\varepsilon$ small enough. We deduce that the  above expression converges to the repartition function of $\max(\mathcal{I}_1,\mathcal{I}_2)$ (see \eqref{defI} for  definitions of $\mathcal{I}_1$ and $\mathcal{I}_2$) when $\varepsilon$ goes to $0$, and this yields Theorem \ref{ADV}.
\hfill$\Box$}

\subsection{Favorite site (proof of Theorem \ref{ADV2})} \\
\bl{Thanks to Section \ref{secontion3}, we know precisely the nature of the contribution of each $h_t$-valley to the local time. The difficulty in proving Theorem \ref{ADV} was the need to consider only a part of the contribution of the last $h_t$-valley.
The proofs of the first two points \eqref{Th15eq1} and \eqref{Th15eq2} of Theorem \ref{ADV2} are thus easier to obtain, since they do not require to "cut" the contribution of any valley. Let us
prove
the first point \eqref{Th15eq1}
(the second one, \eqref{Th15eq2}, is obtained similarly). We have, using \eqref{InegRvTilde1},
\begin{align*}
&
    \mathbb{P} \left[ \loX^{*}(H(m_{N_t + 1})) \leq \alpha t \right]
\\
& \leq
    \mathbb{P} \left ( \loX^{*}\big(H\big(\tilde L_{N_t}\big)\big) \leq
    \alpha  t,
     \ \mathcal{Q}, \mV_t \right )
    + \mathbb{P} \left ( \overline{\mathcal{Q}} \right )
    + P\left(\overline{\mathcal V}_t\right) + \mathbb{P} \left ( \overline{\mathcal{B}_{3}(n_t)} \right )
\\
& \leq
    \mathbb{P} \bigg( \sup_{1 \leq j \leq N_t} \ell_j /t
    \leq (1 - \varepsilon_t)^{-1}\alpha,
    \ \mathcal{Q}, \mV_t \bigg)
    + \mathbb{P} \left ( \overline{\mathcal{Q}} \right )
    + o(1),
\end{align*}
where we fixed some $\varepsilon > 0$ and
$\mathcal{Q}:=\{\varepsilon t \leq  H(m_{N_t}) \leq (1-\varepsilon) t ,\ 1\leq N_t \leq n_t\}$
as
after \eqref{proba1}
(from there we see that $\lim_{\varepsilon \rightarrow 0} \lim_{t \rightarrow +\infty} \mathbb{P} \big( \overline{\mathcal{Q}}\big) = 0$).
In the last inequality we used Proposition \ref{Pro3.3}, Lemma \ref{CVs} and Lemma \ref{negloc1}.
To lighten notation, let
$
    \tilde \alpha_t
:=
    (1 - \varepsilon_t)^{-1}\alpha
$.
We have
\begin{eqnarray*}
&&
    \mathbb{P} \bigg( \sup_{1 \leq j \leq N_t} \ell_j /t \leq \tilde \alpha_t, \ \mathcal{Q}, \ \mV_t \bigg)
\\
& \leq &
    \mathbb{P} \bigg( \sup_{1 \leq j \leq N_t} \ell_j /t \leq \tilde \alpha_t, \ \mathcal{B}_{1}(n_t), \ \mathcal{Q}, \ \mV_t \bigg)
    + \mathbb{P} \left ( \overline{\mathcal{B}_{1}(n_t)} \right)
\\
& \leq &
    \mathbb{P} \bigg( \sup_{1 \leq j \leq N_t} \ell_j /t \leq \tilde \alpha_t, \ \bar{\mathcal{H}}_{N_t} \geq 1-\delta_t', \ \bar{\mathcal{H}}_{N_t-1}\leq 1-\varepsilon+ \delta_t', \ \mathcal{Q} \bigg)
    + o(1),
\end{eqnarray*}
with $\delta_t'=3\tilde v_t/t$ and where we used \eqref{TpsNeg} together with Proposition \ref{Pro3.3}. Partitioning on the values of $N_t$ we get that the above is less than
\[
    \sum_{1\leq k \leq n_t}
    \mathbb{P} \bigg( \sup_{1 \leq j \leq k} \ell_j /t \leq \tilde \alpha_t, \ \bar{\mathcal{H}}_{k} \geq 1-\delta_t', \ \bar{\mathcal{H}}_{k-1}\leq 1-\varepsilon+ \delta_t', \ \mathcal{Q} \bigg)
    + o(1).
\]
Since the sum $\sum_1$ defined in the proof of the upper bound of Proposition \ref{ProSupL} (see \eqref{sum1} and below) is smaller than $s(\varepsilon,t)$ satisfying $\lim_{\varepsilon \rightarrow 0} \lim_{t \rightarrow +\infty} s(\varepsilon,t) = 0$,
 we can intersect the event on the above probability with $\{ k = \mathcal{N}_{t}^{2\varepsilon} \}$ and get
\[
    \mathbb{P} \left[ \loX^{*}(H(m_{N_t + 1})) \leq \alpha t \right]
\leq
    \mathbb{P} \bigg( \sup_{1 \leq j \leq \mathcal{N}_{t}^{2 \varepsilon}} \ell_j /t \leq \tilde \alpha_t \bigg)
    + \mathbb{P} \left( \overline{\mathcal{Q}} \right)
    + s(\varepsilon,t)
    + o(1).
\]
Then, as in the proof of Theorem \ref{ADV} we have that $(\mathcal{Y}_1,\mathcal{Y}_2)$ almost surely satisfies the hypothesis of Lemma \ref{lemContJI}, and is therefore almost surely a point of continuity for $J_{I, 12-\varepsilon}$  defined just above \eqref{DefK}.
From this continuity, Proposition \ref{CVY1Y2} and the {continuous mapping theorem} we get
\[
    \sup_{1 \leq j \leq \mathcal{N}_{t}^{2\varepsilon}} \ell_j /t
=
    J_{I, 1-2\varepsilon} \left ( (Y_1, Y_2)^t \right )
\underset{t \rightarrow +\infty}{\cvloi}
    J_{I, 1-2\varepsilon} (\mathcal{Y}_1,\mathcal{Y}_2)
=
    \mathcal{Y}_1^{\natural} \big( \mathcal{Y}_2^{-1} (1-2\varepsilon)\big).
\]
Then, as in the proof of Theorem \ref{ADV} we have almost surely
$\mathcal{Y}_2^{-1}(1-2\varepsilon) = \mathcal{Y}_2^{-1}(1)$ for all $\varepsilon$ small enough
so $\mathcal{Y}_1^{\natural} \big( \mathcal{Y}_2^{-1} (1-2\varepsilon)\big)$
converges almost surely to $\mathcal{Y}_1^{\natural} \big( \mathcal{Y}_2^{-1} (1)\big)$ when $\varepsilon$ goes to $0$.
Thus, we get
\[
    \limsup_{t \rightarrow +\infty} \mathbb{P} \left[ \loX^{*}(H(m_{N_t + 1})) \leq \alpha t \right ]
\leq
    \mathbb{P} \left ( \mathcal{Y}_1^{\natural} \big( \mathcal{Y}_2^{-1} (1)\big) \leq \alpha \right ).
\]
A lower bound is proved similarly, so we get the following, proving \eqref{Th15eq1}:
\[
    \lim_{t \rightarrow +\infty} \mathbb{P} \left[ \loX^{*}(H(m_{N_t + 1})) \leq \alpha t \right]
=
    \mathbb{P} \left ( \mathcal{Y}_1^{\natural} \big( \mathcal{Y}_2^{-1} (1)\big) \leq \alpha \right ).
\]
}

To obtain the result \eqref{Th15eq3} for the favorite site, we first argue that we essentially need to obtain the asymptotic behavior of $N_t^*/N_t$, where $N_t^*:= \min\{j \geq 1, \loX(m_j,t) = \max_{1\leq k \leq N_t} {\lo(m_k,t)}\}$. Indeed, define for any $\varepsilon\in(0,1/2)$,
\begin{eqnarray*}
    \mathcal{K}_1
& := &
    \left\{ (1- \varepsilon)m_{N_t} \leq X(t) \leq (1+ \varepsilon)m_{N_t}\right\} ,
\\
    \mathcal{K}_2
& := &
    \left\{ (1- \varepsilon)m_{N_t^*} \leq F_t^* \leq (1+ \varepsilon)m_{N_t^*}\right\}.
\end{eqnarray*}
Then, we have, $ \lim_{t \rightarrow + \infty}\P(\mathcal{K}_1)=1$
by the localization result Theorem \ref{ththm} combined with the fact that $X(t)/t^\kappa$ converges in law under $\P$
to a positive limit as $t\to+\infty$ by \cite{KawazuTanaka}.

\bl{Let us now justify that $ \lim_{t \rightarrow + \infty}\P(\mathcal{K}_2)=1$.
According to \eqref{EqProuveeDebutPreuveTheorem} proved at the start of the proof of Theorem \ref{ADV},
to Lemma \ref{negloc2} and \eqref{Bas1}, we have
\[
    \mathbb{P} \left ( \sup_{x \in \R} \loX(t,x)  \geq 2 \tilde w_t, \mathcal{B}_{4}(n_t), N_t \leq n_t \right )
\underset{t \rightarrow + \infty}{\longrightarrow}
    1.
\]
Notice that on the event inside the above probability, for $t$ large enough so that $2 \tilde w_t \geq t e^{-2 \phi(t)}$,
we have $F_t^* \in \Dt_{N_t^*}$ (recall the definition of $\Dt_j$ in \eqref{Dj}).
Since $\Dt_{N_t^*}$ is centered at $m_{N_t^*}$ and its half-length is deterministic and equal to $r_t = C_0 \phi(t)$ we only need to justify that
\[ \mathbb{P} \left ( \varepsilon m_{N_t^*} \geq C_0 \phi(t) \right ) \underset{t \rightarrow + \infty}{\longrightarrow} 1. \]
We have $m_{N_t^*} \geq m_1$ and $\mathbb{P} (m_1 \geq C_0 \phi(t) / \varepsilon) \geq \mathbb{P} (\tilde m_1 \geq C_0 \phi(t) / \varepsilon) - o(1) $ by Lemma \ref{CVs}. So using \eqref{6.2.6}, we thus deduce that $ \lim_{t \rightarrow + \infty}\P(\mathcal{K}_2)=1$.
}

We can now write for $x>0$,
\begin{equation*}
    \P\left[ \frac{F_t^*}{X(t)} \leq x \right]
=
    \P\left[ \frac{F_t^*}{X(t)} \leq x, \mathcal{K}_1,\mathcal{K}_2 \right]+ v(\varepsilon,t)
\leq
    \P\left[ \frac{m_{N_t^*}}{m_{N_t}} \leq x \frac{1+ \varepsilon}{1- \varepsilon} \right]+ v(\varepsilon,t).
\end{equation*}
where $v(\varepsilon,t) \geq 0$,  satisfies $\lim_{\varepsilon \rightarrow 0} \lim_{t \rightarrow + \infty} v(\varepsilon,t)=0$.{
Similarly,  we have
$$
    \P\left[ \frac{F_t^*}{X(t)} \leq x \right]
\geq
    \P\left[ \frac{m_{N_t^*}}{m_{N_t}} \leq x \frac{1- \varepsilon}{1+ \varepsilon} \right]- v(\varepsilon,t).
$$
Hence, we obtain
\begin{align}
 \P\left[ \frac{m_{N_t^*}}{m_{N_t}} \leq x \frac{1- \varepsilon}{1+ \varepsilon} \right]- v(\varepsilon,t) \leq \P\left[ \frac{F_t^*}{X(t)} \leq x \right]  \leq  \P\left[ \frac{m_{N_t^*}}{m_{N_t}} \leq x \frac{1+ \varepsilon}{1- \varepsilon} \right]+ v(\varepsilon,t) \label{5.55}.
\end{align}
So, we observe that we only have to study the random variable $\frac{m_{N_t^*}}{m_{N_t}} $.
For that we first remark that $N_t^*$ and $N_t$ diverge when $t$ goes to infinity.
Indeed by Lemma \ref{6.7}, the correct normalisation for the convergence in law of $N_t$ is $e^{\kappa \phi(t)}$, so  $\P(N_t \geq e^{(1- \varepsilon)\kappa \phi(t)})=1-o(1)$.
For $N_t^*$, we first notice that the previous result for $N_t$ also gives for $t$ large, $\P(N_t \geq e^{(1- \varepsilon/2)\kappa \phi(t)})=1-o(1)$.
Therefore
$$
    \P\Big(N_t^* \leq e^{(1- \varepsilon)\kappa \phi(t)}\Big)
\leq
    \P\Big(  \max_{k \leq e^{(1- \varepsilon )\kappa \phi(t)}} {\lo(m_k,t)}
            \geq
            \max_{k< e^{(1- \varepsilon/2)\kappa \phi(t)}} {\lo(m_k,t)} \Big)
            +o(1).
$$
Now, since
$
    \lo(m_k,t)
=
    \lo\big(\tilde m_k,H\big(\tilde L_k\big)
        \wedge 
        \big(H\big(\tilde m_k\big)+H_{\tilde m_k\to \tilde L_k^-}\big)
        \big)
=:
    \widehat \ell_k
$ 
for $k<N_t$ on $\mV_t\cap\{N_t\leq n_t\}\cap \mathcal{B}_{2}(n_t)$ 
which has probability $1-o(1)$ by Lemmas \ref{CVs} and \ref{lemtps},  
\begin{eqnarray*}
&&
    \P\Big(
        \max_{k \leq e^{(1- \varepsilon )\kappa \phi(t)}}    {\lo(m_k,t)}
        \geq
        \max_{k < e^{(1- \varepsilon/2)\kappa \phi(t)}} {\lo(m_k,t)}
    \Big)
\\
& \leq &
    \P\Big(  \max_{k \leq e^{(1- \varepsilon )\kappa \phi(t)}} {\widehat \ell_k} 
    \geq \max_{k < e^{(1- \varepsilon/2)\kappa \phi(t)}} {\widehat \ell_k} \Big)
    +o(1),
\end{eqnarray*}
with $\big(\widehat \ell_k, k \leq e^{(1- \varepsilon/2)\kappa \phi(t)}\big)$ i.i.d. random variables under $\P$ 
by strong Markov property and the second part of Lemma \ref{CVs}, and 
with queue distributions given by \eqref{cvmesure7.1} and Proposition \ref{Pro3.3}. \\
It is then clear that for large $t$, 
$ 
    \P(  \max_{k \leq e^{(1- \varepsilon )\kappa \phi(t)}} {\widehat \ell_k} \geq \max_{k < e^{(1- \varepsilon/2)\kappa \phi(t)}} {\widehat \ell_k} )
=
    o(1)
$, 
    and we therefore obtain that $\P(N_t^* \geq e^{(1- \varepsilon)\kappa \phi(t)})=1-o(1)$. \\
Then, following the work of \cite{Faggionato}, we know that $(m_i-m_{i-1}, i \geq 2)$ are i.i.d. random variables with a known  Laplace transform (given by (2.19) in \cite{Faggionato}), this allows to compute the first and fourth moments  of $\Delta m_1:= m_2-m_{1}$ and therefore obtain after an  elementary but tedious computation that for large $t$,  $\E(\Delta m_1) \sim C_7 e^{ \kappa h_t} $ ($C_7>0$, see also (2.17) in \cite{Faggionato})
and $\E( (\Delta m_1-\E(\Delta m_1))^4) \sim C_8 e^{4 \kappa h_t}$ ($C_8>0$), which yields as $t\to+\infty$ and $k\to+\infty$,
$$
    \E\left[ \left( m_k/k - \E(\Delta m_1)  \right)^4\right] \sim C_8 e^{4 \kappa h_t}/k^2. 
$$
These facts allow us to write by a Markov inequality that
\begin{eqnarray*}
&&
    \P\left[ \left|m_{N_t}-\E(\Delta m_1) N_t \right|> \varepsilon  \E(\Delta m_1) N_t\right]
\\
& \leq &
    \sum_{j \geq e^{(1- \varepsilon)\kappa \phi(t)}}
    \P\Big[ \big|m_{j}-\E(\Delta m_1) j\big|> \varepsilon  \E(\Delta m_1) j \Big]+ o(1)
\\
 & \leq &
    \sum_{j \geq e^{(1- \varepsilon)\kappa \phi(t)}}  \frac{2C_8 (C_7)^{-4}}{ \varepsilon ^4 j^2}+ o(1)
\\
& \leq &
    C_+ \varepsilon ^{-4} e^{-(1- \varepsilon)\kappa \phi(t)}+o(1).
\end{eqnarray*}
This yields that $\{\left|m_{N_t}-\E(\Delta m_1) N_t \right| \leq \varepsilon  \E(\Delta m_1) N_t\}$ as well as (with a similar computation) $\{\left|m_{N_t^*}-\E(\Delta m_1) N_t^* \right| \leq \varepsilon  \E(\Delta m_1) N_t^*\}$ are realized  with a probability close to one. \\
{Now including these events in the probability in \eqref{5.55}, eventually enlarging $v(\varepsilon,t)$ we get
\begin{align*}
  \P\left[ \frac{{N_t^*}}{{N_t}} \leq x \frac{(1- \varepsilon)^2}{(1+ \varepsilon)^2} \right]- v(\varepsilon,t) \leq \P\left[ \frac{F_t^*}{X(t)} \leq x \right] & \leq  \P\left[ \frac{{N_t^*}}{{N_t}} \leq x \frac{(1+ \varepsilon)^2}{(1- \varepsilon)^2} \right]+ v(\varepsilon,t).
\end{align*}
Notice that the random variables involved now ($N_t^*$ and $N_t$) only depend 
of what happens in the bottom of the $h_t$-valleys, and we have to deal with
\begin{align*}
    \P\left[ \frac{{N_t^*}}{{N_t}} \leq y \right]
= 
    \P\big[ {{N_t^*}}={{N_t}}\big] \un_{\{y=1\}}
    +\P\left[ \frac{{N_t^*}}{{N_t}} \leq y,\ {{N_t^*}}<{{N_t}}\right] \un_{\{y\leq 1\}}
    +\un_{\{y> 1\}}, 
\end{align*}
for any $y>0$.
We are now interested in the limit when $t$ goes to infinity of the above two probabilities. We first use the same  lines as for the proof of Section 5.1, that is to say we give a lower and an upper bound of this probability involving the i.i.d. sequences $(\ell_j,j)$ and $(\mathcal{H}_j,j)$. In the same way we have obtained Proposition \ref{ProSupL}, we then have for any $\varepsilon>0$ and large $t$,
\begin{align*}
\tilde{\mathcal{P}} -v(\varepsilon,t) \leq \P\left({{N_t^*}}={{N_t}}  \right) \leq \tilde{\mathcal{P}}+v(\varepsilon,t)
\end{align*}
with
$$\tilde{\mathcal{P}}:= \P\left[(1- \bar{\mathcal{H}}_{ \mathcal{N}^{2\varepsilon}_t-1})\frac{\bar{\ell}_{ \mathcal{N}^{2\varepsilon}_t}-\bar{\ell}_{ \mathcal{N}^{2\varepsilon}_t-1} } {(\bar{\mathcal{H}}_{ \mathcal{N}^{2\varepsilon}_t}-\bar{\mathcal{H}}_{ \mathcal{N}^{2\varepsilon}_t-1})} > {\max_{ 1\leq j \leq  \mathcal{N}^{2\varepsilon}_t-1} \frac{\ell_j}{t}  } \right],$$
recall that 
$\bar{\mathcal{H}}_k= Y_2(ke^{-\kappa \phi(t)})= \frac{1}{t}\sum_{i=1}^k\mathcal{H}_i$, 
$\bar{\ell}_k=Y_1(ke^{-\kappa \phi(t)})=\frac{1}{t}\sum_{i=1}^k\ell_i$, 
$\mathcal{N}^{2\varepsilon}_t :=\inf\{m \geq 1, \bar{\mathcal{H}}_{m} >1- 2\varepsilon \}$, 
and $v$ is a positive function such that $\lim_{t \rightarrow + \infty} v(\varepsilon,t) \leq \textrm{const} \times \varepsilon^{ \kappa  \wedge (1- \kappa)}$ 
with an eventually larger  const than in Proposition \ref{ProSupL}. 
In the same way, for any $y>0$, $\varepsilon>0$ and $t$ large enough,}
\begin{align*}
    \bar{\mathcal{P}}_1^{-} -v(\varepsilon,t)  
\leq 
    \P\left[ \frac{{N_t^*}}{{N_t}} \leq y,\ {{N_t^*}}<{{N_t}}\right] \un_{y\leq 1}   
\leq 
    \bar{\mathcal{P}}_1^{+} +v(\varepsilon,t),
\end{align*}
where
\begin{align*}
&
    \tilde{\mathcal{P}}^{\pm}_1
\\
& :=
    \P\left[\mathcal{N}_t^*/ \mathcal{N}^{2\varepsilon}_t\leq y \pm \varepsilon,  (1- \bar{\mathcal{H}}_{ \mathcal{N}^{2\varepsilon}_t-1})\frac{\bar{\ell}_{ \mathcal{N}^{2\varepsilon}_t}-\bar{\ell}_{ \mathcal{N}^{2\varepsilon}_t-1} } {(\bar{\mathcal{H}}_{ \mathcal{N}^{2\varepsilon}_t}-\bar{\mathcal{H}}_{ \mathcal{N}^{2\varepsilon}_t-1})} \leq  {\max_{ 1\leq j \leq  \mathcal{N}^{2\varepsilon}_t-1} \frac{\ell_j}{t}  } \right]\un_{y\leq 1},
\end{align*}
with $\mathcal{N}_t^*:= \min\{j \geq 1, \ell_j = \max_{k \leq \mathcal{N}^{2\varepsilon}_t} {\ell_k}\}$.
This together with Lemma \ref{lem4.7} yields for large $t$,
$$
\big|
    \P\left[ {{N_t^*}}={{N_t}}\right]
    -\P\left[  \mathcal{I}_1 < \mathcal{I}_2 \right]
\big| 
\leq 
    \lim_{t \rightarrow + \infty} v(\varepsilon,t)
    +o(1)
$$
and
\begin{eqnarray*}
&&
    \left|
        \P\left[ \frac{{N_t^*}e^{- \kappa \phi(t)}}{{N_t}e^{- \kappa \phi(t)}} \leq y, {{N_t^*}}<{{N_t}}\right]
        -
        \P\left[ \frac{{F^*(\mathcal{Y}_1,\mathcal{Y}_2) }}{\mathcal{Y}_2^{-1}(1)} \leq y,\ \mathcal{I}_1 \geq  \mathcal{I}_2 \right]
    \right|
\\
& \leq &
    \lim_{t \rightarrow + \infty} v(\varepsilon,t),
    +o(1),
\end{eqnarray*}
where $F^*$ is defined at the beginning of Section \ref{secCont}. Replacing $y$ by $x \frac{(1- \varepsilon)^2}{(1+ \varepsilon)^2}$ for the lower bound and by $x \frac{(1+ \varepsilon)^2}{(1- \varepsilon)^2}$ for the upper bound and taking the limit when $t$ goes to infinity and then $\varepsilon \rightarrow 0$ we obtain for $0<x<1$, 
 \begin{align*}
    \lim_{t \rightarrow + \infty }\P\left[ \frac{{N_t^*}}{{N_t}} \leq x \right] 
 =  
    \P\left[ \frac{{F^*(\mathcal{Y}_1,\mathcal{Y}_2) }}{\mathcal{Y}_2^{-1}(1)} \leq x,\ \mathcal{I}_1 \geq  \mathcal{I}_2 \right].
\end{align*} }

 To finish the proof of the last result of Theorem \ref{ADV2} we finally have to prove Lemma \ref{5.5b} below.

\medskip

\begin{lemma} \label{5.5b}
The random variable $\frac{F^*(\mathcal{Y}_1,\mathcal{Y}_2)}{\mathcal{Y}_2^{-1}(1)}$ follows a uniform law $U_{[0,1]}$
and is independent of the couple $\left (\mathcal{I}_1 ,  \mathcal{I}_2 \right )$.
\end{lemma}

\medskip

\noindent{\bf Proof:}
\bl{
For any $s>0$, let $\mathcal{G}_1(s):=\inf\{ u \leq s,\ \mathcal{Y}_1(u) - \mathcal{Y}_1(u-)= \mathcal{Y}_1^{\sharp}(s)  \}.$ The fact that for every $s>0$, $\mathcal{G}_1(s)/s$ follows a uniform distribution is basic. Since the independence that we seek is specific we give some details.

The process of the jumps of $(\mathcal{Y}_1, \mathcal{Y}_2)$ in $[0, s]$ is a Poisson point process in $[0, s] \times (\mathbb{R}_+)^2$ (the coordinate in $[0, s]$ for the instant when the jump occurs and the other coordinate for the jump) with intensity measure $\lambda \times \nu$ where $\lambda$ is the Lebesgue measure on $[0, s]$ and $\nu$, as defined in the introduction, is the L\'evy measure of $(\mathcal{Y}_1, \mathcal{Y}_2)$. Let us give a particular construction of the process $(\mathcal{Y}_1, \mathcal{Y}_2)$ on $[0, s]$:

Let $(P_n)_{n \geq 1}$ be a countable partition of $(\mathbb{R}_+)^2$ by Borelian sets such that $\forall n \geq 1, \ 0 < \nu(P_n) < +\infty$. For each $n$ we define an \textit{i.i.d.} sequence $(S^n_k)_{k \geq 1}$ of random variables in $(\mathbb{R}_+)^2$, an \textit{i.i.d.} sequence $(U^n_k)_{k \geq 1}$ of random variables in $[0, s]$ and a random variable $T_n$ such that
\begin{itemize}
\item
$\displaystyle{
    \forall n \geq 1, \ S^n_1 \sim \nu(. \cap P_n) / \nu(P_n), \ U^n_1 \sim U_{[0,s]}, \ T_n \sim \mathcal{P} (\nu(P_n)),
}$
\item For any $n \geq 1$, the variables $(S^n_k)_{k \geq 1}$, $(U^n_k)_{k \geq 1}$ and $T_n$ are independent,
\item The triplets $\left ( (S^n_k)_{k \geq 1}, (U^n_k)_{k \geq 1}, T_n \right )_{n \geq 1}$ are independent.
\end{itemize}
We know that the random set
\[ \mathcal{S}_n := \left \{ (U^n_k, S^n_k), \ n \geq 1, \ 1 \leq k \leq T_n \right \} \]
is a Poisson point process in $[0, s] \times (\mathbb{R}_+)^2$ with intensity measure $\lambda \times \nu$. Since $(\mathcal{Y}_1, \mathcal{Y}_2)$ is pure jump it is equal in law to the process $(\mathcal{Z}_1, \mathcal{Z}_2)$ defined by
\[
    \forall r \in [0, s],
\qquad
    (\mathcal{Z}_1, \mathcal{Z}_2)(r) = \sum_{n \geq 1, 1 \leq k \leq T_n} S^n_k \ \mathds{1}_{U^n_k \leq r}.
\]
In particular
\begin{align*}
    \mathcal{Z}^{\sharp}_1(s)
& =
    \max \{ \pi_1 (S^n_k), \ n \geq 1, \ 1 \leq k \leq T_n \}, \\
    \mathcal{G}^\mathcal{Z}_1(s)
& = \inf \left \{ U^n_k, \ n \geq 1, \ 1 \leq k \leq T_n, \ \pi_1 (S^n_k) = \mathcal{Z}^{\sharp}_1(s) \right \}, \\
    \mathcal{Z}_1(s)
& = \sum_{n \geq 1, 1 \leq k \leq T_n} \pi_1 (S^n_k),
\qquad
    \mathcal{Z}_2(s) = \sum_{n \geq 1, 1 \leq k \leq T_n} \pi_2 (S^n_k).
\end{align*}
We thus have that $\mathcal{G}_1(s)/s \egloi U_{[0,1]}$ and it is independent from $(\mathcal{Y}^{\sharp}_1(s), \mathcal{Y}_1(s), \mathcal{Y}_2(s))$ and from the sigma-field $\sigma ( (\mathcal{Y}_1, \mathcal{Y}_2)
(t+s) - (\mathcal{Y}_1, \mathcal{Y}_2)(s), \ t \geq 0)$.

We now have to replace $s$ by $\mathcal{Y}_2^{-1}(1)$. For that we can consider for example the dyadic approximations of  $\mathcal{Y}_2^{-1}(1)$,
that is,
$(t_n := \max \left \{ k \in \mathbb{N}, \frac{k}{2^n} < \mathcal{Y}_2^{-1}(1) \right \},n)$. Then, partitioning on the values of $t_n$, using the independence we just proved and the fact that $\mathcal{G}_1(s)/s$ follows a uniform distribution on $[0, 1]$ we get that $\mathcal{G}_1(t_n)/t_n$ follows a uniform distribution on $[0, 1]$ and is independent from
\begin{eqnarray}
   \big(
        (\mathcal{Y}^{\sharp}_1(t_n),
        \ \mathcal{Y}_2(t_n),
        \ \mathcal{Y}_1(t_n + 2^{-n}) - \mathcal{Y}_1(t_n),
        \ \mathcal{Y}_2(t_n + 2^{-n}) - \mathcal{Y}_2(t_n)
   \big).
\label{favptunif}
\end{eqnarray}
We let $n$ goes to infinity, $t_n$ converges almost surely to $\mathcal{Y}_2^{-1}(1)$ from below. As a consequence, $\mathcal{G}_1(t_n)/t_n$ converges almost surely to $\frac{F^*(\mathcal{Y}_1,\mathcal{Y}_2)}{\mathcal{Y}_2^{-1}(1)}$
while the quadruple in \eqref{favptunif} converges almost surely to
\begin{align*}
&
    \big(
        \mathcal{Y}^{\sharp}_1(\mathcal{Y}_2^{-1}(1)-),
        \ \mathcal{Y}_2(\mathcal{Y}_2^{-1}(1)-),
\\
&
        \qquad\qquad\qquad\ \
        \ \mathcal{Y}_1(\mathcal{Y}_2^{-1}(1)) - \mathcal{Y}_1(\mathcal{Y}_2^{-1}(1)-),
        \ \mathcal{Y}_2(\mathcal{Y}_2^{-1}(1)) - \mathcal{Y}_2(\mathcal{Y}_2^{-1}(1)-)
    \big).
\end{align*}
As a consequence, $\frac{F^*(\mathcal{Y}_1,\mathcal{Y}_2)}{\mathcal{Y}_2^{-1}(1)}$ follows a uniform distribution on $[0,1]$ and is independent from the above quadruple for which $\left (\mathcal{I}_1 ,  \mathcal{I}_2 \right )$ is a measurable function, this yields the lemma.
}
\hfill$\Box$

{
\section{Results and additional arguments from the paper \cite{AndDev} }\label{Section6}

\subsection{Some estimates on the diffusion $X$}

The first lemma below gives the right normalisation in law of the number of $h_t$-valleys visited by $X$ before time $t$.

\medskip

\begin{lemma}[number of visited $h_t$-valleys] \label{6.7}
Assume that $ 0<\kappa<1$.
Then,  under the annealed law $\P$,
${N_t} e^{-\kappa \phi(t)}\to_{t\to+\infty}\mathcal{N}$ in law.
The law of $\mathcal{N}$  is determined by its Laplace transform:
\begin{equation}
\label{eqTransfoLaplaceLimiteNbValleesVisitees}
    \forall u>0,
\qquad
    \E\left(e^{-u\mathcal{N}}\right)
=
    \sum_{j=0}^{+ \infty}
    \frac{1}{\Gamma(\kappa j+1) }
    \left(\frac{ -u}{\CK}\right)^j,
\end{equation}
where $\CK$ is a positive constant. Moreover $\P(N_t > n_t) \leq e^{-\phi(t)}$.
\end{lemma}

\bigskip

\noindent{\bf Proof:}
The convergence in distribution is exactly Proposition 1.6 of \cite{AndDev}. For the second fact we have $\P(N_t \geq  n_t) \leq \P(\tilde N_t \geq n_t)+\P(\overline{\mathcal V}_t) \leq \P(\tilde N_t \geq n_t) +e^{[-\k /2+o(1)]h_t}$ by Lemma \ref{CVs}, with $\tilde N_t:= \max\{j\geq 1,\ \tm \leq \sup_{s \leq t} X(s)\}$. Then equation (5.3) in \cite{AndDev} gives $ \P(\tilde N_t \geq n_t) \leq \exp(-2 \phi(t)) $, which yields the result.  \hfill$\Box$

\noindent \\ The lemma below deals with the renewal structure we speak about on the introduction,
and the consequence on the hitting time $H(m_{N_t})$  of the ultimate $h_t$-valley visited by $X$ before time $t$.

\medskip

\begin{lemma} \label{5.4} Assume $0<\k<1$ and $0<\delta<\inf\{2/27, \k^2/2\}$.
For $t>0$, let $\mu_t$ be the positive measure on $\R_+$ such that
$$
    \forall x\geq 0,
\qquad
    \mu_t([0,x])
:=
    e^{-\kappa \phi(t)}
    \sum_{j=1}^{n_t}\P\left(\bar{\mathcal{H}}_j \leq x\right).
$$

Recall that for any  $k$, $\bar{\mathcal{H}}_k:= \sum_{j=1}^{k} {\mathcal{H}}_j/t $, and $\mathcal{H}_1=R_1S_1 \bf{e}_1$ is defined in Proposition \ref{Pro3.3}.
Then, $(\mu_t)_t$ converges vaguely as $t\to+\infty$ to 
$\mu$ defined by
$$
\dd \mu(x):= (\CK \Gamma(\kappa))^{-1}  x^{\kappa-1}\un_{(0,+\infty)}(x) \dd x,
$$
with $\CK$ is the same constant as in Lemma \ref{6.7}. For $0\leq r <s\leq 1$,
\begin{align}
    \lim_{t \rightarrow + \infty}
    \P\bigg(1-s\leq  \frac{H(\mf_{N_t})}{t} \leq 1-r\bigg)
& =  \frac{\sin(\pi \kappa )}{ \pi}
    \int_{1-s}^{1-r} x^{ \kappa-1} (1-x)^{-\kappa} \textnormal{d}x
\label{6.5.2}.
\end{align}
\end{lemma}

\bigskip

\noindent{\bf Proof:}
\noindent The first part of the above lemma is very close to Lemma 5.1 of \cite{AndDev}, indeed Proposition \ref{Pro3.3} gives the proximity between the random variables $(U_i,i \leq n_t)$ and the random variables $(\mathcal{H}_i, i \leq n_t)$, moreover an important preliminary result in \cite{AndDev} (Proposition 4.1) states that $e^{\kappa \phi(t) }(1-\E(e^{- \lambda U_1 /t}))= \mathcal{C}_ \kappa \lambda^{\kappa}+o(1)$ for large $t$. So we also know that \begin{align} e^{\kappa \phi(t) }(1-\E(e^{- \lambda \mathcal{H}_1 /t}))= \mathcal{C}_ \kappa \lambda^{\kappa}+o(1) \label{6.91b},\end{align} notice that this result could also be deduced from \eqref{cvmesure9.1} with the help of a Tauberian theorem. Then by independence of the random variables $\mathcal H_j$ and the fact that they are i.d., for any  $\lambda>0$
\begin{eqnarray*}
    \int_0^{+ \infty} e^{- \lambda x} \dd \mu_t(x)
& = &
     \frac{1}{e^{\kappa \phi(t)}}\sum_{j=1}^{n_t}\left(\E\Big(e^{-\lambda\ \frac{\mathcal{H}_1}{t}}\Big) \right)^j
\end{eqnarray*}
By  \eqref{6.91b} as  $n_t e^{-\k\phi(t)}\to_{t\to+\infty}+\infty$, $[\E\left(e^{- \lambda \mathcal{H}_1/t}\right)]^{n_t+1}=o(1)$. Hence, we get as $t\to +\infty$,
again by \ref{6.91b}
\begin{eqnarray*}
    \int_0^{+ \infty} e^{- \lambda x} \dd \mu_t(x)
& = &
    \frac{e^{-\kappa \phi(t)}(1
    +o(1))}{1-\E\left(e^{- \lambda  \mathcal{H}_1/t}\right)}
    +o(1)
 =
   \frac{1}{\mathcal{C}_ \kappa \lambda^{\kappa}}
   +o(1)
\\
& = &
    \int_0^{+ \infty} \frac{e^{-\lambda x}x^{\kappa-1}}{\mathcal{C}_ \kappa \Gamma(\kappa)} \textnormal{d}x+o(1),
\end{eqnarray*}
which gives the vague convergence of measure $(\mu_t)_t.$
Also \eqref{6.5.2} is equation (1.2) of Corollary 1.5 in \cite{AndDev}.
\hfill$\Box$

\bigskip

{
\noindent In Lemma \ref{Lemma63} below,  we approximate $\tilde h_j$, the exit time of $h_t$-valley number $j$
(if $X$ leaves it on the right),
by a product of $3$ simpler random variables.
To this aim, we recall that with the notation of Lemma \ref{lemX} and of its proof,
for each $1\leq j \leq n_t$,
$\tilde {R}_j= \int_{\tilde \tau_j^-(h_t/2)}^{\tilde \tau_j(h_t/2)}e^{-\tV(x)} {\dd x}$,
and $A^j(u)=\int_{\tm}^{u} {e^{\tV(x)}}\dd x$, $u\in\R$.
Moreover, for some independent Brownian motions $ B^j$, $1\leq j \leq n_t$, independent of $\wk$,
\begin{eqnarray*}
    \tilde  h_j
& = &
    \int_{\tL^-}^{\tL} {e^{-\tV(u)}}\mathcal{L}_{B^j}[\tau^{B^j}(A^j(\tL)),A^j(u)]\dd u,
\\
    {\bf e}_j
& = &
    \mathcal{L}_{B^j}\big[\tau^{B^j}(A^j(\tL)),0\big]/A^j(\tL).
\end{eqnarray*}

\medskip

\begin{lemma}\label{Lemma63}
Let $0<\varepsilon<\inf\{2/27, \k^2/2\}$.
For large $t$, we have for every $1\leq j \leq n_t$,
\begin{equation}\label{EqApproxhj}
    \P\left(  \left | \tilde h_j-A^j(\tL) {{\bf e}_j} \tilde {R}_j \right|
        >
        2e^{-(1-3\varepsilon)h_t/6} A^j(\tL) {{\bf e}_j}  \tilde {R}_j
    \right)
\leq
    C_+ e^{-(c_-) \varepsilon h_t}.
\end{equation}
\end{lemma}

\bigskip

\noindent{\bf Proof:}
We first notice that $\big(\tilde h_j,A^j(\tL), {{\bf e}_j}, \tilde {R}_j\big)$
is measurable with respect to the $\sigma$-field generated by
$\big(\tilde{V}^{(j)}(x+\tilde {L}^+_{j-1}),\ 0\leq x \leq \tilde {L}^+_j-\tilde {L}^+_{j-1} \big)$
and $B^j$, so, thanks to the second fact of Lemma \ref{CVs},
its law under $\P$ does not depend on $j$.
Thus, the left hand side of \eqref{EqApproxhj} does not depend on $j$. Hence
we just have to prove \eqref{EqApproxhj} for $j=2$.

This is actually already proved in \cite{AndDev}, for which it is an important step.
Indeed in this paper \cite{AndDev},
our $A^j$, $\tilde B^2$ and $\tilde h_2$
are denoted respectively by $\tilde A_j$, $B$ and $\bU$, as defined in (\cite{AndDev}, eq. (3.17) and (3.18)),
and our $\tilde {R}_2$  and ${\bf e}_2$
by $\Im$ and ${\bf e}_1$, as defined in (\cite{AndDev}, after eq. (4.17)).
Hence our \eqref{EqApproxhj} for $j=2$ is exactly (\cite{AndDev}, Lemma 4.7), which proves our lemma.

The proof of (\cite{AndDev}, Lemma 4.7) is quite technical, however we can give a simple heuristic
in order for the present paper to be more self-contained.
The idea of the proof of (\cite{AndDev}, Lemma 4.7) is that, loosely speaking, for $u$ close to $\tilde m_j$,
that is for $u\in\big[\tilde \tau_j^-(h_t/2), \tilde \tau_j(h_t/2)\big]$,
$
\mathcal{L}_{B^j}[\tau^{B^j}(A^j(\tL)),A^j(u)]
$
is nearly
$
\mathcal{L}_{B^j}[\tau^{B^j}(A^j(\tL)),0]
=A^j(\tilde L_j){\bf e_j}
$,
whereas for $u$ far from $\tilde m_j$, that is for $u\in[\tilde L_j^-, \tilde L_j]$
but $u\notin [\tilde \tau_j^-(h_t/2), \tilde \tau_j(h_t/2)]$,
$e^{-\tilde V^{(j)}(x)}$ is "nearly" $0$, with large probability.
Finally, combining these heuristics gives
$\tilde h_j\approx A^j(\tilde L_j){\bf e_j} \tilde R_j$.
\hfill$\Box$
}

\bigskip

The following lemma is used to prove Lemma \ref{lemX} and uses the notation of this lemma,
and where the independent r.v. $G^+(h_t/2,h_t)$, $F_1^+(h_t)$, $F_2^-(h_t/2)$ and $F_3^-(h_t/2)$
defined before Proposition \ref{Pro3.3}.

\bigskip

\begin{lemma} \label{LemmaConstructionS2}
Assume $0<\delta<\inf\{2/27, \kappa^2/2\}$.
For large $t$, possibly on an enlarged probability space,
there exists $R_2\egloi F^{-}_2(h_t/2)+F^{-}_3(h_t/2)$
and $S_2\egloi F^{+}_1(h_t)+G^{+}(h_t/2, h_t)$,
such that  $R_2$, $S_2$ and ${\bf e}_2$ are independent
and
\begin{equation}\label{EqApproxIntbySCase2}
    P\Bigg(
        \left\{
            \left| \int_{\tilde m_2}^{\tilde L_2} {e^{\tilde V^{(2)}(x)}}\dd x - {S}_2\right| \leq e^{-(d_-)h_t} {S}_2,\,
            \tilde R_2=R_2
        \right\}
    \Bigg)
\geq
    1- e^{-(D_-)h_t},
\end{equation}
where $D_->0$.
\end{lemma}

\bigskip

\noindent{\bf Proof:}
Due to (\cite{AndDev} Lemma 4.5) with its notation, we have
$\Ip_0:=\int_{m_2}^{\tau_2(h_t)} e^{V^{(2)}(x)}\dd x\egloi F^+(h_t)$,
$\Ip_2:=\int_{\tau_2(h_t)}^{L_2} e^{V^{(2)}(x)}\dd x\egloi G^+(h_t/2, h_t)$,
$\Im_1:=\int_{m_2}^{\tau_2(h_t/2)} e^{-V^{(2)}(x)}\dd x\egloi F^-(h_t/2)$
and finally
$\Im_2:=\int_{\tau_2^-(h_t/2)}^{m_2} e^{-V^{(2)}(x)}\dd x\egloi F^-(h_t/2)$
with
$L_2:=\inf\{x>\tau_2(h_t),\, V^{(2)}(x)=h_t/2\}$.
The problem is that $\Ip_0$ is not independent of $\Im_1$, so we would like to replace it by some $\Ip_1\egloi\Ip_0$ of it
with better independence properties.
It is proved in (\cite{AndDev}, at the top of page 32) that for large $t$,
possibly in an enlarged probability space, there exists $\Ip_1$
such that
$
    |\Ip_0-\Ip_1|
\leq
    e^{-(1-3\delta)h_t/2} \Ip_1
$
with probability greater than $1-4e^{-\kappa \delta h_t/2}$ and
where $\Ip_1\egloi F^+(h_t)$ by (\cite{AndDev}, eq. (4.35)).

Let $S_2:=\Ip_1+\Ip_2\geq \Ip_1$.
Notice that on $\mV_t$, by Remark \ref{RemEgaliteAvecouSansTilde},
$\tilde R_2 =\Im_1+\Im_2 =:R_2$
and
$
    \int_{\tilde m_2}^{\tilde L_2} {e^{\tilde V^{(2)}(x)}}\dd x
=
    \int_{m_2}^{L_2} {e^{V^{(2)}(x)}}\dd x
=
    \Ip_0+\Ip_2
$.
The two previous inequalities give
$
    \big|\int_{\tilde m_2}^{\tilde L_2} {e^{\tilde V^{(2)}(x)}}\dd x-S_2\big|
=
    \big|\Ip_0-\Ip_1\big|
\leq
    e^{-(1-3\delta)h_t/2} S_2
$ and $\tilde R_2 =R_2$
with probability at least $1- 5e^{-\kappa \delta h_t/2}$ thanks to Lemma \ref{CVs}.
This proves \eqref{EqApproxIntbySCase2}.

Moreover, by (\cite{AndDev}, Prop. 4.4 (i)), $\Ip_1$, $\Ip_2$, $\Im_1$, $\Im_2$
and ${\bf e}_2$ (which is denoted by ${\bf e}_1$ in \cite{AndDev})
are independent. So, ${\bf e}_2$, $S_2=\Ip_1+\Ip_2$ and $R_2=\Im_1+\Im_2$
are independent, and
$R_2\egloi F^{-}_2(h_t/2)+F^{-}_3(h_t/2)$
and $S_2\egloi F^{+}_1(h_t)+G^{+}(h_t/2, h_t)$.
\hfill $\Box$

\medskip

\noindent The last lemma of this section tells that with large probability, the diffusion $X$ leaves
every $h_t$-valley $[\tilde L_j^-, \tilde L_j]$, $1\leq j \leq n_t$ from its right.
Recall that $B^j$ is defined after \eqref{DefAj}.

\bigskip

\begin{lemma}\label{LemmaProbaRetourEnmi}
For large $t$, there exists $c_->0$ such that
\begin{equation}
    \P \bigg[ \cap_{j=1}^{n_t} \Big\{ \max_{u<\tL^-}\mathcal{L}_{B^j}[\tau^{B^j}(A^j(\tL)),A^j(u)]=0 \Big\}   \bigg]
\geq
    1-e^{- (c_-) h_t} \label{6.3.2}.
    \end{equation}
\end{lemma}

\bigskip

\noindent{\bf Proof:}
 \eqref{6.3.2} is essentially Lemma 3.2 in \cite{AndDev}:  \\
 Indeed, recall the definition of $\mathcal{A}_j:= \{ \max_{u<\tL^-}\mathcal{L}_{B^j}[\tau^{B^j}(A^j(\tL)),A^j(u)]=0 \}$, we have $ \cap_{j=1}^{n_t}\mathcal{A}_j = \cap_{j=1}^{n_t} \{ H_j(\tL) < \{H_j(\tL^-) \}$, with, for any $ \tL^- \leq x \leq \tL $, $H_j(x)= \inf\{ s>0, B_j(s)=x\}$, with $B_j$  a Brownian motion. Therefore $\Pw(\mathcal{A}_j)$ is equal to the probability $\Pw(\overline{\mathcal{E}}_j)$ of Lemma 3.2 in \cite{AndDev}. It is proved in this lemma see (3.10) that for large $t$, $P(\mathcal{B}:=\{ \Pw(\overline{\mathcal{E}}_j) \leq e^{- (\kappa /2) h_t}\} ) \geq 1-3e^{-\kappa \delta h_t}$, so we obtain \eqref{6.3.2} as $\P(\overline{\mathcal{E}}_j) \leq E(\Pw(\overline{\mathcal{E}}_j) \un_{\mathcal{B}}  ) + P(\overline{\mathcal{B}})  \leq e^{- c_- h_t}/n_t$, for $c_->0$ small enough.
\hfill$\Box$

\subsection{Some estimates on the potential $W_\kappa$ and its functionals}

\noindent \\
\noindent We start this section with the Laplace transform of the important functional $\mathcal{R}_{\kappa}$ :

\medskip

\begin{lemma} \label{LTR}
Recall that $0<\k<1$.
For any $\gamma>0$,
\begin{equation}\label{eqTransfoLaplaceRk}
    E\big(e^{-\gamma \mathcal{R}_{\kappa}}\big)
=
    \left(\frac{ (2\gamma)^{{\k/2}}}{\k\Gamma(\k)I_\k(2\sqrt{2\gamma})}\right)^2.
\end{equation}
Moreover, $\mathcal{R}_{\kappa}$ admits moments of any positive order.
\end{lemma}

\medskip

\noindent  {\bf Proof:}
$\int_0^{+ \infty} e^{- W_{\kappa}^{\uparrow}(u)} \dd u$ is the limit in law under $P$ of
$\int_0^{\tau^{W_{\kappa}^{\uparrow}}(x)} e^{- W_{\kappa}^{\uparrow}(u)} \dd u$ as $x\to+\infty$.  This limit
is given by (\cite{AndDev}, Lemma 4.2), which proves \eqref{eqTransfoLaplaceRk}.
Note that in (\cite{AndDev}, Lemma 4.2), $W_{\kappa}^{\uparrow}$ is denoted by $R$, and $\int_0^{\tau^{W_{\kappa}^{\uparrow}}(x)} e^{- W_{\kappa}^{\uparrow}(u)} \dd u$
is denoted respectively by $F^-(x)$.
Moreover the Laplace transform of $\mathcal{R}_{\kappa}$ is of class $C^{\infty}$ on a neighborhood of $0$ since $x\mapsto x^\k/I_\k(x)$ is $C^{\infty}$ on such a neighborhood
(see e.g. \cite{Borodin} p. 638).
Therefore $\mathcal{R}_{\kappa}$ admits moments of any positive order.
\hfill$\Box$

\bigskip

The following Lemma is a series of estimates concerning the different coordinates of valleys.

\medskip

\begin{lemma} For $t$ large enough, for every $1\leq i \leq n_t$,
\begin{align}
     & P(0<M_0<m_1)\leq C_+  h_t e^{-\k h_t}, \label{6.2.2} \\
    & P(\tilde \tau_{i+1}^*(h_t)\neq\tilde \tau_{i+1}(h_t) ) \leq    C_+ h_t e^{-\k h_t}, \label{6.2.1}\\
    & P\Big(\inf_{[\tilde \tau_i^-(h_t^+), \tilde \tau_i^-(h_t)]}\tilde V^{(i)} < h_t/2 \Big) \leq e^{-\kappa h_t /8}, \label{6.2.3}  \\
    & P(\tilde L_i^+ -\tilde L_i^- \geq 40 h_t^+/ \kappa ) \leq  e^{- \kappa h_t / 8},
    \label{6.2.5} \\
    & P(\tilde \tau_{i}(h)-\tilde m_i \geq 8h/ \kappa ) \leq C_+ e^{- \kappa h/ (2 \sqrt 2)}, \qquad 0 \leq h \leq h_t,  \label{6.2.4} \\
    & P(\tilde m_1 \leq r ) \leq e^r \exp\big(\big( \kappa/2- \sqrt{2+ \kappa^2/4}\big) h_t^+\big)=o(1), \qquad \forall r=o(h_t^+).  \label{6.2.6}
\end{align}
\end{lemma}

\bigskip

\noindent  {\bf Proof:}
\eqref{6.2.2} follows from eq. (2.8) of \cite{AndDev}; \eqref{6.2.1} is eq. (3.41) of \cite{AndDev}. \eqref{6.2.3} and \eqref{6.2.5} are respectively  eq. (2.34) and (2.32) of Lemma 2.7 of \cite{AndDev}.
Moreover, \eqref{6.2.4} is eq. (2.22) of the same reference. For \eqref{6.2.6}, we know from definitions in \eqref{eqDefTaui1} that $ \tilde m_1   \geq \tilde L_1^{\sharp} = \tau^{W_{\kappa}}(-h_t^+)$, where $ \tau^{W_{\kappa}}(-h_t^+)$ is the first positive time the drifted Brownian motion $W_{\kappa}$ reaches $-h_t$. Using a Markov inequality together with (2.0.1) page 295 of \cite{Borodin} we obtain $P( \tau^{W_{\kappa}}(-h_t^+)  \leq r) = P( e^{-\tau^{W_{\kappa}}(-h_t^+)}  \geq e^{-r})  \leq  e^r e^{( \kappa/2- \sqrt{2+ \kappa^2/4}) h_t^+}$, which is exactly \eqref{6.2.6}.
\hfill$\Box$

\bigskip

The lemma below deals with two functionals involving coordinates far from the bottom $\tilde m_1$
of the first visited $h_t$-valley $[\tilde L_1^-, \tilde L_1]$.

\bigskip

\begin{lemma} \label{lemI3} There exists $c_->0$ such that for any $\varepsilon>0$ and $t$ large enough,
\begin{align*}
P\left( \int_{\tilde \tau_1(h_t/2) }^{ \Lt_1} e^{-\tilde V^{(1)}(x)}\dd x
\leq C_+ h_t^2 e^{-(1- \varepsilon)h_t/2} \right) \geq 1-e^{- (c_-) \varepsilon h_t},
\\
P\left( \int_{ \Lt_1^-}^{\tilde \tau_1^-(h_t/2) } e^{-\tilde V^{(1)}(x)}\dd x  \leq C_+ h_t^2 e^{-(1- \varepsilon)h_t/2} \right)
\geq
    1-e^{- (c_-) \varepsilon h_t}.
\end{align*}
\end{lemma}

\bigskip

\noindent{\bf Proof:}
The proof is inspired from steps 1 and 2 of Lemma 4.7 of \cite{AndDev}. For the first integral, let
$$
    \B_1^{}
:=
    \big\{\inf\nolimits_{[\tts1,\tilde{\tau}_1(h_t)]}\tilde V^{(1)}>(1-\e) h_t/2\big\},
\qquad
    \B_2^{}
:=
    \big\{\tilde L_1^{+}-\tilde L_1^- \leq 40 h_t^+ / \kappa\big\}.
$$
We have on
$
\B_1^{}\cap \B_2^{}$,
\begin{equation}\label{eqApproxJ2Bis}
    \int_{\tts1}^{\tilde L_1} e^{-\tilde V^{(1)}(u)}\dd u
\leq    e^{- (1- \e) h_t/2} \big[\tilde L_1-\tts1 \big]
\leq
    \frac{40 h_t^+ h_t}{\kappa} e^{- (1- \e) h_t/2}.
\end{equation}
Now, Fact \ref{Fact_Williams}, equation  (\ref{3.10}) with $\alpha=1/2$, $\gamma=(1-\e)/2$ and $\omega=1$,
and Lemma \ref{CVs} give
$$
    P\big(\overline{\B_1^{}} \big)
\leq
    P\big[\inf\nolimits_{[\tau_1(h_t/2),\tau_1(h_t)]} V^{(1)}\leq(1-\e) h_t/2,\mV_t\big]
    + P(\overline{\mV}_t)
\leq
    3e^{-\k \e h_t/2}.
$$
Moreover,
$
    P\big(\overline{\B_2^{}}\big)
\leq
    e^{-\k h_t/8}
\leq
    e^{-\k \e h_t/2}
$
by \eqref{6.2.5} since we can take $\e<1/4$.
The second inequality, can be proved similarly.
\hfill$\Box$

\bigskip

\begin{lemma}  Recall that for $h>0$,
$
\beta_0(h):=
    E
    \left(
        \int_{0}^{\tau_1^*(h)}
        e^{\wk(u)}\dd u
    \right)
$, with  $\tau_1^*(h) :=
    \inf\{u\geq 0,\ \wk(u)-\inf\nolimits_{[0,u]}\wk\geq h\}$.
  For large $h$,
\begin{equation}
    \beta_0(h)
\leq
    C_+ e^{(1-\k)h}. \label{6.1.1}
\end{equation}
\end{lemma}

\bigskip

\noindent{\bf Proof:}
\eqref{6.1.1} is (\cite{AndDev}, eq. (3.38)),
since in \cite{AndDev}, $\beta_0(h)$ is defined at the top of page 23
and $\tau_1^*(h)$ in its Lemma 3.6.
\hfill$\Box$
}

\section{Appendix}\label{SectAppendix}

\subsection{Some estimates for Brownian motion, Bessel processes, $\BP$ and their functionals}

\noindent We provide in this section some known formulas for some processes that appear in our study.
The first lemma is about Laplace transforms of the
exponential functionals defined in \eqref{eqDefF+G+} and \eqref{eqDefF+G+Bis}.
Its proof can be found in (\cite{AndDev}, Lemma 4.2).
Recall that $C_+$ (respectively $c_-$) is a positive constant that is as large (resp. small) as needed.


\bigskip

\begin{lemma} \label{lem4.5}
There exist $C_9>0$, $M>0$ and $\eta_1\in(0,1)$ such that $\forall y>M, \forall \gamma\in(0, \eta_1]$,
\begin{align}
\left|E \left(e^{-\gamma F^+(y)/e^y}\right)-[1-2\gamma/(\k+1)]\right|
    & \leq
    C_9\max(e^{-\k y},\gamma^{3/2})
    ,\label{Ff2}\\
\left|E \left(e^{-\gamma G^+(y/2,y)/e^y}\right)-[1-\Gamma(1-\kappa)(2\gamma)^{\k}/\Gamma(1+\k)]\right|
    & \leq
    C_9\max(\gamma^\k e^{-\k y/2},\gamma)
    .\label{Gpf2}
\end{align}
Moreover, there exists $C_{10}>0$ such that for all $y >0$, $E \left( {F^+(y)}/{e^y}  \right)\leq C_{10}$.
\end{lemma}

\bigskip

Recall that $\BP$ is a $(-\kappa/2)$-drifted Brownian motion $W_{\kappa}$ Doob-conditioned to stay positive
(see above \eqref{eqDefF+G+}). We have,

\medskip

\begin{lemma}\label{Lemma72}
Let $0< \gamma< \alpha  <\omega $. For all $h$ large enough, we have
\begin{align}
    P^{\alpha h}\left(\tau^{\BP^{}}(\gamma  h)<\tau^{\BP^{}} (\omega h)  \right)
& \leq
    2 e^{-\kappa (\alpha-\gamma) h }
    \label{3.10},
\\
    P\left(\tau^{\BP^{}}( \omega h)-\tau^{\BP^{}} (\alpha h) \leq 1  \right)
& \leq
    4 e^{- [(\omega-\alpha) h]^2/3},
    \label{3.10b}
\\
   P\big(\tau^{\BP}(h)>8h/\k \big)
& \leq
    C_+ e^{- \kappa h /(2\sqrt{2})},
    \label{bessel4}
\\
    P\big(\tau^{\BP}(h) \leq h \big)
& \leq
    C_+ e^{-(c_-) h},
\label{bessel5}\\
    P\big(\tau^{\BP^{}}( \gamma h) \leq 1  \big)
& \leq
    C_+e^{- (c_-)[\gamma h]^2},
    \label{3.10c}
\end{align}
where $P^{\alpha h}$ denotes the law of $\BP$ starting from $\alpha h$. Moreover the first inequality is still true if $\omega$ is a function of $h$ such that $\lim_{h \rightarrow + \infty} \omega(h)= + \infty$.
\end{lemma}

\medskip

\noindent{\bf Proof:}
The first 3 inequalities come from (\cite{AndDev}, Lemma 2.6).
The fact that, in \eqref{3.10}, $\omega$ can actually be taken as a function of $h$ comes directly from eq. (2.31) of \cite{AndDev},
which shows that the right hand side of \eqref{3.10} is equivalent to $e^{-\kappa (\alpha-\gamma) h }$ as $h\to+\infty$
if $w=w(h)\to_{h\to+\infty}+\infty$.
\eqref{3.10c} is a consequence of \eqref{3.10b} with $\omega=\gamma$ and $\alpha=\gamma/2$.
We turn to \eqref{bessel5}.
By (\cite{AndDev}, eq. (2.7) and Fact 2.1, coming from \cite{Faggionato}),
$
    E\big( e^{- \alpha \tau^{\BP}(h)} \big)
\sim_{h \rightarrow + \infty}
    \textrm{const}.
    e^{h(\kappa/2-\sqrt{2 \alpha+ \kappa^2/4})}
$,
in particular for $\alpha=1-\kappa$.
Then a Markov inequality for $P\big( e^{-\alpha \tau^{\BP}(h)}>e^{-\alpha h} \big)$
proves \eqref{bessel5} since $1-\kappa/2 -\sqrt{2(1-\kappa)+\kappa^2/4}<0$.
\hfill$\Box$


\bigskip

We also need the following lemma, focusing only on some exponential functionals.

\medskip

\begin{lemma} \label{EstimR} Recall that $F^{\pm}$ and $G^+$ are defined in \eqref{eqDefF+G+}
and \eqref{eqDefF+G+Bis}. For all $0<\zeta \leq 1$ and $0<\e<1$,
for $h$ large enough,
\begin{align}
  &  P\left[e^{(1- \e)\zeta h }\leq F^{+}(\zeta h)  \leq e^{(1+ \e) \zeta  h }\right]
\geq
    1-4 e^{-   \kappa \e\zeta h/2}, \label{4.6bb}
\\
&
    P\big[F^-(h)  \geq e^{-\varepsilon h}\big]
\geq
    1-e^{-(c_-) \varepsilon^2 h^2}, \label{Fmoins}
\\
&
    P\big[G^+(\alpha h,h)  \leq  b(h) e^{h}\big]
\geq
    1- C_+ [b(h)]^{- \kappa},
\qquad
    0<\alpha<1, \, b(h)>0. \label{G}
\end{align}
\end{lemma}

\bigskip

\noindent{\bf Proof:}
By Markov inequality and the last line of Lemma \ref{lem4.5},
$$
    P\big[ F^{+}(\zeta h) > e^{(1+ \e)\zeta h }\big]
\leq
    C_{10} e^{-\varepsilon \zeta h }
\leq
    e^{-\kappa\e\zeta h /2}
$$
for large $h$.
For the lower bound, we have by (\cite{AndDev}, eq. (2.29)) for large $h$,
$$
    P\big[F^+(\zeta h)\geq e^{(1-\e)\zeta h}\big]
\geq
    1-3e^{-\kappa\e\zeta h /2}.
$$
These two inequalities prove \eqref{4.6bb}.
For \eqref{Fmoins},
first $F^-(h)  \geq e^{-\varepsilon h} \tau^{\BP^{}}(\varepsilon h)$,
and using \eqref{3.10c}, $\tau^{\BP^{}}(\varepsilon h)  \geq 1$
with a probability larger than $1-e^{-(c_-) \varepsilon^2 h^2 }$, which proves \eqref{Fmoins}.
Finally,  notice that in law
$
    G^+(\alpha h,h)
\leq
    e^h \int_0^{+ \infty} e^{\wk(x)}\dd x
=
    e^h
    A_{\infty}
$.
By \cite{Dufresne}, $2/A_\infty$ is a gamma variable of parameter $(\k,1)$,
and so has a density equal to  $e^{-x}x^{\k-1}\un_{\R_+}(x)/\Gamma(\k)$, which leads to \eqref{G}.
\hfill$\Box$

\bigskip

The following lemma is exactly Lemma 4.3 in \cite{AndDev} which  proof can be found in that paper.

\medskip

\begin{lemma}\label{LemmaContinuiteEn0}
Let $(B(s),\, s\in\R)$ be a standard two-sided Brownian motion. For every $0<\e<1$, $0<\delta<1$ and $x>0$,
\begin{align}
&
    \P\bigg(\sup_{u\in[-\delta, \delta]}\big|\mathcal{L}_B\big(\tau^B(1), u\big)-\mathcal{L}_B\big(\tau^B(1), 0\big)\big|
                >\e \mathcal{L}_B\big(\tau^B(1), 0\big)\bigg)
\leq
    C_{+}\frac{\d^{1/6}}{\e^{2 /5}}, \label{Dev} \\
&
    \P\bigg(\sup_{u\in[0,1]} \mathcal{L}_B\big(\tau^B(1), u\big) \geq x \bigg)
\leq
    4e^{-x/2}, \label{Diel} \\
&
    \P\left(\sup_{u\leq 0} \mathcal{L}_B\big(\tau^B(1), u\big) \geq x \right)
\leq
    4/x. \label{Diel2}
\end{align}
\end{lemma}

\bigskip

The next lemma says that with large probability, a $2$-dimensional squared Bessel Process is bounded by some deterministic function.
This lemma may be of independent interest.

\bigskip

\begin{lemma}\label{lemmaTechniqueMajorationBessel2}
Let $(\BQ_2(u),\ u\geq 0)$ be a Bessel process of dimension $2$, starting from $0$,
and two functions $a(.)$ and $k(.)$ from $(0,+\infty)$ to $(0,+\infty)$, having limit $+\infty$ on $+\infty$.
We have for large $t$,
\[  P\Big(
        \forall u\in(0,k(t)],\
            \BQ_2^2(u)
        \leq
            2 e  \big[a(t)+4\log\log[e k(t)/u] \big] u
    \Big)
\geq
    1-C_+\exp[-a(t)/2]. \]
\end{lemma}

\bigskip

\noindent{\bf Proof:}
We consider for $t>0$ and $i\in\N$,
$$
    \B_{1,i}
:=
    \bigg\{
        \sup_{[k(t)/e^{i+1}, k(t)/e^i]}\BQ_2^2
        \leq
            2 \frac{ k(t)}{e^{i}}[a(t)+4\log(i+1)]
    \bigg\},
\qquad
    \B_{2}
:=
    \bigcap_{i=0}^\infty
    \B_{1,i}.
$$
We recall that there exist two standard independent Brownian motions $(B_1(u),\ u\geq 0)$ and $(B_2(u),\ u\geq 0)$
such that
$(\BQ_2^2(u), \ u\geq 0)$ is equal in law to $(B_1^2(u)+B_2^2(u),\ u\geq 0)$.
So for $i\in\N$,
\begin{eqnarray*}
    P\big(\overline{\B}_{1,i}\big)
& \leq &
    2
    P\big(
        \sup\nolimits_{[k(t)/e^{i+1}, k(t)/e^i]}B_1^2
        >
         k(t)e^{-i}[a(t)+4\log(i+1)]
    \big)
\\
& \leq &
    4
    P\Big(
        \sup\nolimits_{[0, k(t)/e^i]}B_1
        >
        \sqrt{ k(t)e^{-i}[a(t)+4\log(i+1)]}
    \Big)
\\
& = &
    4
    P\left(
        |B_1(1)|
        >
        \sqrt{a(t)+4\log(i+1) }
    \right)
\\
& \leq &
    8\exp[- a(t)/2-2\log(i+1)]
\end{eqnarray*}
for large $t$ so that $a(t)\geq 1$, by scaling, and since
$B_1\egloi -B_1$,
$\sup_{[0,1]} B_1\egloi |B_1(1)|$ and
$P(B_1(1)\geq x)\leq e^{-x^2/2}$ for $x\geq 1$.
Consequently for large $t$,
\begin{equation}\label{InegProbaE2LemmeInegBessel}
    P\big(\overline{\B}_{2}\big)
 \leq
    \sum_{i=0}^\infty
    P\big(\overline{\B}_{1,i}\big)
\leq
    8\exp[-a(t)/2]
    \sum_{i=0}^\infty \frac{1}{(i+1)^2}
=
    C_+\exp[-a(t)/2].
\end{equation}
Now,
let $0< u\leq k(t)$. There exists $i\in\N$ such that $k(t)/e^{i+1}< u \leq k(t)/e^i$.
We have,
$e^i\leq k(t)/u$,
so
$e^{i+1}\leq e k(t)/u$
and then
$
    \log(i+1)
\leq
    \log\log [e k(t)/u]
$.
Consequently on $\B_2$,
$$
    \BQ_2^2(u)
\leq
    2 \big(k(t)/e^{i}\big)[a(t)+4\log(i+1)]
\leq
    2 e u \big[a(t)+4\log\log[e k(t)/u] \big].
$$
This, combined with \eqref{InegProbaE2LemmeInegBessel},
proves the lemma.
\hfill$\Box$

We also need some estimates on the local time of $B$ at a given coordinate $y \in \R$ at the inverse of the local time of $B$ at 0. Recall that $\sigma_B(r,y)= \inf\{s>0, \ \lo_B(s,y) >r\}$ for $r>0$, $y\in \R$. By the second Ray-Knight Theorem, the processes
$(\lo_B(\sigma_B(r,0),y),y \in \R_+)$ and $(\lo_B(\sigma_B(r,0),-y),y \in \R_+)$  are two independent squared  Bessel processes
 of dimension $0$ starting at $r$.
The following lemma is proved in (\cite{Talet1}, Lemma 3.1;
the results are stated for a Bessel process but are actually true for a {\it squared} Bessel process; see also \cite{Diel},  Lemma 2.3).

\bigskip

\begin{lemma} \label{LemBessel0} We denote by $(Q_0(y), \, y\geq 0)$ the square of a $0$-dimensional Bessel process  starting at $1$.  Let $M>0$, $u>0$ and $v>0$.
Then,
\begin{align}
&
    P\left(\sup_{ 0 \leq y \leq v } \big|Q_0(y)-1\big| \geq u \right)
\leq
    4 \frac{\sqrt{(1+ u)v}}{u}\exp\left[-u^2/(8(1+u) v) \right],
\label{5.96} \\
&
    P\left(\sup_{y \geq 0} Q_0(y) \geq M \right)
=
    1/M.
\label{5.97}
\end{align}
\end{lemma}

\bigskip

\section*{Acknowledgements}

We are grateful to an anonymous referee for his or her very careful reading of the paper,
and for comments that helped us improve the clarity and the presentation of the paper.


\bibliographystyle{alea3}
\bibliography{thbiblioE_ALEA}

\end{document}